\input amstex
\documentstyle{amsppt}
\magnification=\magstephalf   
\hcorrection{0.55truein}\vcorrection{0.125truein}
\parindent=1pc             
\indenti=2pc               
\pagewidth{30pc}           
\pageheight{47pc}          
\refstyle{B}
\let\bs\backslash

\let\hra\hookrightarrow

\define\hyf{\snug@-\ignorespaces}

\let\isom\cong
\redefine\cong{\equiv}

\redefine\Im{\operatorname{Im}}
\let\lmt\longmapsto
\let\lra\longrightarrow

\let\ol\overline

\define\simto{\overset\sim\to\longrightarrow}
\let\tensor\otimes

\predefine\imaginary{\Im}
\predefine\Par{\P}
\redefine\pmod#1{~(\text{\rm mod}~#1)}
\predefine\sect{\S}


\loadeusm 	      
\let\script\eusm

\loadeusb             

\loadbold            


\define\?{\llap{\bf ??\qquad}\ignorespaces}

\define\0{{\bold 0}}
\define\1{{\bold 1}}
\define\fra{{\frak a}}

\define\atild#1{\tilde a_{\bold #1}}

\define\alphabar{\bar{\alpha}}

\define\alphatild{\tilde{\alpha}}
\define\Alb{\operatorname{Alb}}
\define\AJ{\operatorname{AJ}}
\define\alt#1{A_{#1}\,}

\define\boldb{{\boldkey b}}
\let\b\boldb

\define\bigdot{{\sssize{\bullet}}}
\define\betatild{\tilde{\beta}}
\define\cref#1{{\bf (???)}[{\it #1}]}
\define\CC{{\Bbb C}}
\define\CH{\operatorname{CH}}
\define\Chow#1#2{\CH^{#1}(#2,\QQ)}
\define\Chowhom#1#2{\CH^{#1}_{\text{hom}}(#2,\QQ)}
\define\Chowalb#1#2{\CH^{#1}_{\Alb}(#2,\QQ)}
\define\Chowaj#1#2{\CH^{#1}_{\AJ}(#2,\QQ)}

\define\chitild{\widetilde{\chi}}
\define\chibar{\ol{\chi}}

\define\CKd{Chow-K\"unneth decomposition}
\define\cl{\gamma}

\define\Corr{\operatorname{Corr}}

\define\bsce{{\boldkey{e}}}
\define\e{\bsce}

\define\ebar{{\bar{\bsce}}}

\define\scE{E}
\define\E{{E}}

\define\Estar{\E^{*}}
\define\Ebar{\overline{\E}\vphantom{E}}
\define\Einf{\E^{\infty}}
\define\Etild{\widetilde{\E}\vphantom{E}}
\define\Ebarinf{\Ebar^\infty}

\define\Etwo{\lup2\E}
\define\Etwostar{\lup2\E^*}
\define\Etwobar{\lup2{\vphantom{E}}\Ebar}
\define\Etwotild{\lup2{\vphantom{E}}\Etild}
\define\Etwobarinf{\Etwobar^\infty}
\define\Etwotildinf{\Etwotild^\infty}

\define\Ektild{\lup k{\vphantom{E}}\Etild}

\define\End{\operatorname{End}}
\define\frG{{\frak G}}
\define\G{\frG}

\define\Gal{{\operatorname{Gal}}}
\define\tGamma{{\tr{\Gamma}}}
\define\graph#1{[\Gamma_{#1}]}
\define\tgraph#1{[\tGamma_{#1}]}
\define\graphof#1{[\Gamma(#1)]}
\define\tgraphof#1{[\tGamma(#1)]}

\define\Gm{{\Bbb{G}_{\text{m}}}{\vphantom{\Bbb{G}}}}

\define\frH{{\frak H}}
\define\Hdot{H_{\bigdot}}

\define\Hom{\operatorname{Hom}}

\define\id{{\operatorname{id}}}
\define\inv{^{-1}}

\define\Jac{\operatorname{Jac}}
\define\Ker{\operatorname{Ker}}

\define\LL{{\Bbb L}}
\define\lambdabar{\ol{\lambda}\vphantom{\lambda}}
\define\lambdatild{\tilde{\lambda}\vphantom{\lambda}}

\define\lsub#1{\mathstrut_{#1}\kern-0.05em}
\define\lsup#1{\mathstrut^{#1}\kern-0.07em}
\let\lup\lsup
\define\m{{\boldkey m}}
\define\scM{{\script{M}}}
\define\M{{M}}
\define\Mbar{{\ol{M}}}
\define\Minf{{M^{\infty}}}
\define\mubar{\bar\mu\vphantom{\mu}}
\define\mutild{\tilde\mu\vphantom{\mu}}

\define\NS{{\operatorname{NS}}}

\define\phistar{\phi^{*}}

\define\phibar{\ol{\phi}\vphantom{\phi}}
\define\phitild{\widetilde{\phi}\vphantom{\phi}}
\define\phitwo{\lup2\phi}
\define\phitwostar{\lup2\phistar}
\define\phitwobar{\lup2\vphantom{\phi}\phibar}
\define\phitwotild{\lup2\vphantom{\phi}\phitild}

\define\psisharp{\psi^{\#}}
\define\psitwosharp{\lup2\psisharp}
\redefine\P{{\script P}}
\define\PP{{\Bbb P}}
\define\pbar{\ol{p}}

\define\tildepi{\widetilde\pi}
\define\barpi{\ol{\pi}}
\define\pibar#1{\ol{\pi}_{#1}(\Ebar/\Mbar)}
\define\pibari{\pibar{i}}
\define\pibarf{\pibar{f}}
\define\pibarinf{\pibar{\infty}}
\define\pibarrel#1{\ol{\pi}_{#1}^{\text{rel}}(\Ebar/\Mbar)}
\define\pitild#1{\widetilde{\pi}_{#1}}
\define\pitildi{\pitild{i_1,i_2}}
\define\pitildf{\pitild{f}}
\define\pitildinf{\pitild{\infty}}
\define\pican#1{\pi_{#1}^{\text{can}}}

\define\Pic{\operatorname{Pic}}
\define\pr{\operatorname{pr}}

\define\QQ{{\Bbb Q}}

\define\Qdot{{\QQ_{\sssize\bullet}}}
\define\Qed{\hbox to 0.5em{ }\nobreak\hfill\hbox{$\square$}}

\define\rel{^{\text{rel}}}

\define\frS{{\frak S}}
\redefine\S{\frS}

\define\SL{{\operatorname{SL}}}

\define\Span{\operatorname{Span}}
\define\spec{\operatorname{Spec}}

\define\sbangle#1{_{\langle #1 \rangle}}
\define\spangle#1{^{\langle #1 \rangle}}
\define\sparen#1{^{(#1)}}
\define\ustar#1{#1^*}
\define\lstar#1{#1_*}
\define\Sym{\operatorname{Sym}}
\define\sym#1{S_{#1}\,}

\define\taubar{\bar\tau}
\define\tautild{\tilde\tau}
\define\thetabar{\ol{\vartheta}}
\define\thetatild{\tilde\vartheta}

\define\tr#1{{\vphantom{#1}^{\text{t}}\kern-0.05em#1}}
\define\twedge{{\tsize{\bigwedge}}}
\define\U{{\script{U}}}
\define\V{{\script{V}}}

\define\scW{{\script W}}
\define\W#1{\lup#1{\scW}}

\define\ZZ{{\Bbb Z}}
\define\ZmodN{\ZZ/N\ZZ}
\define\ZmodNsq{(\ZmodN)^2}
\catcode`\@=11
\def\hookrightarrowfill{$\m@th\mathord\lhook\mkern-3mu%
 \mathord-\mkern-6mu%
  \cleaders\hbox{$\mkern-2mu\mathord-\mkern-2mu$}\hfill
  \mkern-6mu\mathord\rightarrow$}
\def\hookleftarrowfill{$\m@th\mathord\leftarrow\mkern-6mu%
  \cleaders\hbox{$\mkern-2mu\mathord-\mkern-2mu$}\hfill
  \mkern-6mu\mathord-\mkern-3mu\mathord\rhook$}
\atdef@ C#1C#2C{\ampersand@\setbox\z@\hbox{$\ssize
 \;{#1}\;\;$}\setbox\@ne\hbox{$\ssize\;{#2}\;\;$}\setbox\tw@
 \hbox{$#2$}\ifCD@
 \global\bigaw@\minCDaw@\else\global\bigaw@\minaw@\fi
 \ifdim\wd\z@>\bigaw@\global\bigaw@\wd\z@\fi
 \ifdim\wd\@ne>\bigaw@\global\bigaw@\wd\@ne\fi
 \ifCD@\hskip.5em\fi
 \ifdim\wd\tw@>\z@
 \mathrel{\mathop{\hbox
to\bigaw@{\hookrightarrowfill}}\limits^{#1}_{#2}}\else
 \mathrel{\mathop{\hbox to\bigaw@{\hookrightarrowfill}}\limits^{#1}}\fi
 \ifCD@\hskip.5em\fi\ampersand@}
\atdef@ D#1D#2D{\ampersand@\setbox\z@\hbox{$\ssize
 \;\;{#1}\;$}\setbox\@ne\hbox{$\ssize\;\;{#2}\;$}\setbox\tw@
 \hbox{$#2$}\ifCD@
 \global\bigaw@\minCDaw@\else\global\bigaw@\minaw@\fi
 \ifdim\wd\z@>\bigaw@\global\bigaw@\wd\z@\fi
 \ifdim\wd\@ne>\bigaw@\global\bigaw@\wd\@ne\fi
 \ifCD@\hskip.5em\fi
 \ifdim\wd\tw@>\z@
 \mathrel{\mathop{\hbox
to\bigaw@{\hookleftarrowfill}}\limits^{#1}_{#2}}\else
 \mathrel{\mathop{\hbox to\bigaw@{\hookleftarrowfill}}\limits^{#1}}\fi
 \ifCD@\hskip.5em\fi\ampersand@}
\catcode`\@=13
\define\today{\number\day\space   \ifcase\month\or January\or
  February\or March\or April\or May\or June\or July\or August\or
  September\or October\or November\or December\fi \space\number\year}
\topmatter
\title  Chow Motives of Elliptic \\
    Modular Surfaces and Threefolds  \endtitle
\rightheadtext{Chow Motives of Elliptic Modular Threefolds}
\author B.B.~Gordon$^1$ and J.P.~Murre$^2$ \endauthor
\leftheadtext{Gordon and Murre}
\affil  University of Oklahoma and Leiden University  \endaffil
\address  Department of Mathematics, University of Oklahoma, 601~Elm,
  Room~423, \newline Norman, OK\enspace 73019  U.S.A. \endaddress
\email  {\tt bgordon\@ou.edu}  \endemail
\address  Department of Mathematics and Computer Science,
  Leiden University, \newline Neils  Bohrweg 1, P.O. Box 9512,
  2300~RA Leiden, The Netherlands  \endaddress
\email  {\tt murre\@rulwinw.LeidenUniv.NL}  \endemail
\thanks \strut$^1$~Partially supported by NSA Grant MDA904--92--H--3093,
  NATO Research Grant CRG931416, The Thomas Stieltjes Institute for
  Mathematics (The Netherlands), The University of Leiden, and NSA
  Grant~MDA904--95--1--0004.
\newline\indent
  \strut$^2$~Partially supported by NATO Research Grant CRG931416,
  The Clarence J. Karcher Endowment (University of Oklahoma), and The
  University of Oklahoma
\endthanks
\abstract
   The main result of this paper is the proof for elliptic modular
threefolds of some conjectures formulated by the second-named author
and shown by Jannsen to be equivalent to a conjecture of
Beilinson on the filtration on the Chow groups of smooth projective
varieties.  These conjectures are known to be true for surfaces in
general, but for elliptic modular surfaces we obtain more precise
results which are then used in the proof of the conjectures for
elliptic modular threefolds.

 Let $\phi:E\to M$ be the universal elliptic curve with level@-$N$
structure, whose smooth completion is an elliptic modular surface
{}~$\Ebar$.  An elliptic modular threefold is a desingularization
$\Etwotild$ of the fibre product $\Ebar\times_\Mbar \Ebar$.
  The first main result is that there exists a decomposition of the
diagonal $\Delta(\Etwotild)$ modulo rational equivalence as a sum of
mutually orthogonal idempotent correspondences $\pi_i$ which lift the
K\"unneth components of the diagonal modulo homological equivalence.
These correspondences act on the Chow groups of ~$\Etwotild$, and
secondly we show that $\pi_i\cdot \CH^j(\Etwotild) =0$ for $i<j$ or
$i>2j$; the implication of this is that there is a filtration on
$\CH^j(\Etwotild)$ that has $j$~steps, as predicted by the general
conjectures.  The third main result is that the first step of this
filtration, the kernel of $\pi_{2j}$ acting on $\CH^j(\Etwotild)$,
coincides with the kernel of the cycle class map from $\CH^j(\Etwotild)$
into the cohomology $H^{2j}(\Etwotild)$; which is to say that there is
a natural, geometric description for this step of the filtration.  We
also identify $F^2\CH^3(\Etwotild)$ as the Albanese kernel.  As a
by-product of our methods we also obtain some information about the Chow
groups of the Chow motives for modular forms $\W k$ defined by Scholl,
for $k=1$ and~2, for example that $\CH^{2}(\W 1) =
\CH^2_{\text{Alb}}(\Ebar)$, and that $\CH^3(\W2) = F^3\CH^3(\Etwotild)$
lives at the deepest level of the filtration, within the Albanese
kernel.
  \endabstract
\endtopmatter
\document
\tolerance=1000            
\mathsurround1pt           
\parskip=\smallskipamount    
\newpage                   

\head \bf Introduction \endhead
\noindent
  Let $X$ be a smooth projective variety of dimension ~$d$.  The
second-named author has conjectured that as an element of the Chow group
$\CH^d(X\times X)\tensor\QQ$ the diagonal $\Delta(X)$ can be decomposed
as a sum $\Delta(X) = \sum_{i=0}^{2d} \pi_i$ of mutually orthogonal
idempotent correspondences modulo rational equivalence which lift the
K\"unneth components of the diagonal \cite{Murre, 1993}.  These
Chow-K\"unneth components of the diagonal, which are not in general
canonical, act on the Chow groups of $X$ to give a decomposition of the
form $\CH^j(X)\tensor \QQ = \bigoplus_{i=0}^{2d} \pi_i \cdot
(\CH^j(X)\tensor\QQ)$.  Then conjecturally $\pi_i \cdot \CH^j(X)\tensor
\QQ =0$ for $i<j$ or $i>2j$; and when this is the case, the filtration
defined by $F^\nu\CH^j(X)\tensor\QQ := \bigoplus_{i=j}^{2j-\nu} \pi_i
\cdot (\CH^j(X)\tensor \QQ)$ has precisely $j$~steps.  A third conjecture
asserts that this filtration is independent of the choice of projectors
{}~$\pi_i$; and as the first step in this direction, a fourth conjecture
proposes that $F^1\CH^j(X)\tensor\QQ$ is precisely the kernel of the
cycle class map into cohomology.  U.~Jannsen has shown that these
conjectures of the second-named author \cite{op.~cit} together are
equivalent to conjectures of Beilinson on the existence of a canonical
filtration on the Chow groups of smooth projective varieties
\cite{Jannsen, 1994}; see also \cite{Bloch, 1980} \cite{Beilinson,
1987}.  The class of varieties for which the conjectures are known to be
true is still very small:  For curves, it is elementary (compare
\cite{Manin, 1968}, \cite{Kleiman, 1972}); for surfaces, see \cite{Murre,
1990}; for products of surfaces and curves, see \cite{Murre, 1993, ~II};
for uniruled threefolds, see \cite{del Angel and M\"uller-Stach, 1996};
and the existence of a Chow-K\"unneth decomposition is known for abelian
varieties \cite{Shermenev, 1974} \cite{Deninger and Murre, 1991}
\cite{K\"unnemann, 1994}.

  The main result of the present paper is the proof of the conjectures
(except for some points concerning the canonicity of the filtration) for
elliptic modular threefolds.  To describe these, let $N\ge 3$ be an
integer (which we suppress from the notation), let $M:=M_N$ be the
modular curve parameterizing elliptic curves with full level-$N$
structure, and let let $\phi :E \to M$ be the universal elliptic curve
(with full level-$N$ structure) over ~$M$.  Then the smooth completion
$\phibar : \Ebar\to\Mbar$ of $E$ over the compactification $\Mbar$ of
{}~$M$ obtained by adjoining the cusps is an elliptic modular surface
\cite{Shioda, 1972}.  The fibre product $\Etwobar := \Ebar\times_\Mbar
\Ebar$ over $\Mbar$ has only rational double points for singularities,
and by blowing these up we get the nonsingular elliptic modular threefold
{}~$\Etwotild$ that is the main focus of our attention; such threefolds
have also been studied in \cite{Schoen, 1986}.  For the fibre
products $\Ebar\times_\Mbar \cdots \times_\Mbar\Ebar$ ($k\ge 1$ times)
there is a natural desingularization $\Ektild$ due to \cite{Deligne,
1969}, but see also \cite{Scholl, 1990}; the first-named author of
the present paper has
looked at the cohomology and the Hodge structure of these ~$\Ektild$, and
verified the generalized Hodge conjecture for them \cite{Gordon, 1993}.

  To prove the conjectures for elliptic modular surfaces and threefolds,
we begin by constructing projectors for $\Ebar$ that extend the canonical
relative projectors that are known for $E$ as elliptic curve scheme over
{}~$M$ \cite{Deninger and Murre, 1991} \cite{K\"unnemann, 1994}.  Using
these projectors, we construct a finer Chow motive decomposition of
$\Ebar$ than the Chow-K\"unneth decomposition, and thus obtain more
precise results about the filtration of the Chow groups of $\Ebar$ than
can be proved for surfaces in general \cite{Murre, 1990}.  Then using the
fibre product structure of $\Etwobar$ we can construct projectors for
$\Etwotild$ which extend the relative tensor products (over ~$M$) of the
canonical relative projectors for ~$E/M$.  The projectors we get
this way, together with a detailed knowledge of the cohomology of
$\Etwotild$, then give us a Chow-K\"unneth decomposition from
which we are able to deduce the other conjectures as well.  We expect
that these methods can be generalized to give a Chow-K\"unneth
decomposition for $\Ektild$ for any $k\ge 1$; however, this becomes
technically more complicated, and we intend to return to it later.

For $k\ge1$ Scholl has constructed Chow motives $\W k :=
\vphantom{\script W}_N^k\kern-0.06em\script W$ supported on ~$\Ektild$,
and he has shown that their cohomology groups are the parabolic
cohomology groups attached to cusp forms of weight $(k+2)$ and
level@-$N$ \cite{Scholl, 1990}.  Not surprisingly, we also encounter
these motives, for $k=1$ and~2, and we recover the same results about
their cohomology.  But then we also study their Chow groups:  For $\W1$
we show that (modulo torsion) $\CH(\W1) = \CH^2(\W1) =
\CH^2_{\Alb}(\Ebar)$, the kernel in $\CH^2(\Ebar)$ of the Albanese map;
for $\W2$ we find that (modulo torsion) it has only two Chow groups,
namely $\CH^2(\W2)$, which is related to the intermediate Jacobian
$J^2(\Etwotild)$, and $\CH^3(\W2)$, which we find lies in the deepest
level of the filtration on $\CH^3(\Etwotild)$.

  The paper is organized as follows.  In section one we recall the
definitions and some facts about Chow motives, and give the precise
statements of the conjectures.  In section two we collect together some
of the facts we need about elliptic modular surfaces and threefolds.  In
section three we construct projectors which extend the canonical
relative relative projectors for $E/M$ to the fibre variety $\Ebar$ over
{}~$\Mbar$, and we also construct projectors for $\Etwotild$ that
extend the tensor products of these canonical relative projectors.  We
also need some extra projectors to account for the degenerate fibres over
the cusps, and we introduce these in section three as well.  Section four
is the technical center of the paper, for there we identify the
motives defined by the projectors of section three with motives
supported on varieties of lower dimension and the ``Scholl motives''
{}~$\W k$, for $k=1,2$.  In section five we study the cohomology of the
motives from section four, and obtain Chow-K\"unneth decompositions for
$\Ebar$ and ~$\Etwotild$.  Finally in section six we use the
Chow-K\"unneth decompositions from section five to study the Chow groups
of $\Ebar$ and ~$\Etwotild$ and obtain the desired results about the
filtrations on those Chow groups.  The main results are stated precisely
in Theorems ~4.2, 5.1 and ~6.2, and each section has a small introduction
of its own.

  It is a pleasure to thank A.~Besser, M.~Hanamura and A.J.~Scholl for
valuable conversations related to this project.

\head \bf 1. Chow motives and the conjectures. \endhead
\noindent
  Let $k$ be a field, let $\rho: S\to \spec k$ be a smooth, connected,
quasi-projective scheme, and let $\V(S)$ be the category of projective
$S$-schemes $\lambda: X\to S$ with $\lambda$ smooth.  When
$S=\spec k$ we write $\V(k)$ for $\V(\spec k)$.  Let $\CH^j(X)$ denote
the Chow group of codimension ~$j$ algebraic cycles on ~$X$
modulo rational equivalence, and $\Chow j X := \CH^j(X)\tensor_\ZZ
\QQ$.  In those cases where we need to consider the Chow group of a
singular
variety
{}~$V$ we write $\CH_i(V)$ for the Chow group in the sense of
\cite{Fulton,
1984} of dimension ~$i$ algebraic cycles on ~$X$
modulo rational equivalence, and $\CH_i(V,\QQ) := \CH_i(V)\tensor_\ZZ
\QQ$.   For a cycle $Z$ on $X$ we write $[Z]$ for its class in
$\CH^j(X,\QQ)$ or $\CH_i(X,\QQ)$.

\subhead 1.1. The category of Chow motives \endsubhead
  To establish some general notations and fix ideas, we briefly recall
some basic definitions and properties for the category of Chow
motives over ~$S$, specifically allowing the possibility that $S=\spec
k$.  For more details see \cite{Scholl, 1994} and \cite{Deninger and
Murre, 1991}.

\subhead 1.1.1. Definition of the category of Chow motives
\endsubhead
  Let $X$ and $Y$ in $\V(S)$, and for convenience we assume that $X$ is
connected and of relative dimension $d_S(X)$ over ~$S$.  Then the group
of relative correspondences of degree ~$r$ from $X$ to $Y$ is
$$
\Corr^r_S(X,Y) := \Chow{d_S(X) +r}{X\times_S Y} .
$$
 There is also the usual bilinear composition
$$ (\alpha,\beta) \lmt \beta\circ\alpha :=
 \pr_{13*}(\pr_{12}^*(\alpha)\cdot\pr_{23}^*(\beta)) ,
$$
 where $\pr_{ij}:X_1\times X_2\times X_3 \to X_i\times X_j$ is the
projection and the intersection product is taken in $\CH^*(X_1\times_S
X_2\times_S X_3,\QQ)$.  Then the category $\scM(S)$ of Chow motives over
$S$ can be defined by:  Objects are triples $(X, p, m)$, where
$X$ is in $\V(S)$, and $m\in\ZZ$, and $p\circ p = p \in \Corr^0_S(X,X)$
is an idempotent (projector); and morphisms are given by
$$\spreadlines{1\jot} \align
\Hom_{\scM(S)}((X/S,p,m),(Y/S,q,n)) :&= q\circ \Corr^{n-m}_S(X,Y) \circ p
\\
  &=q\circ \Chow{d_S(X) +n-m}{X\times_S Y} \circ p .
\endalign
$$
 When $m=0$ we usually write $(X/S,p)$ for $(X/S,p,0)$.

\subsubhead 1.1.2. Examples \endsubsubhead
  (a) ~There is a unit object in $\scM(S)$, namely $\1_S :=
(S,\id_S)$.  More generally, when $X$ in $\V(S)$ (connected) has a
rational section $e:S\to X$ (or still more generally, a relative
zero-cycle of degree one), then $\1_S \simeq (X/S,(e(S)\times_S X))$.

\smallpagebreak
 (b) ~By definition, the Lefschetz motive is $\LL_S := (S,\id_S,-1)$, and
more generally we let $\LL_S^d := \LL_S\tensor_S \dots \tensor_S \LL_S =
(S,\id_S,-d)$.  When $d_S(X) =d$ and there exists a rational section
$e:S\to X$, then also $\LL_S^d \simeq (X/S, (X\times_S e(S))$.  In
particular, when $S=\spec k$ and $X$ is any curve with a rational point
{}~$e$, then $\LL \simeq (X, X\times\{e\})$.

\smallpagebreak
 (c) ~For simplicity, let $S=\spec k$, and let $\alpha\in\Chow p X$ and
$\beta \in \Chow q X$, where $p+q=\dim X$, and suppose the intersection
multiplicity $(\alpha\cdot\beta)=1$.  Then $\alpha\times\beta\in \Chow d
{X\times X}$ is a projector, and  $(X, \alpha\times\beta) \simeq \LL^q$
in ~$\scM(k)$.  In fact, $\alpha\in \Corr^{-q}(X,\spec k)$ and
$\beta\in\Corr^q(\spec k,X)$, and these induce inverse isomorphisms.

\subsubhead 1.1.3. The tensor product in $\scM(S)$ \endsubsubhead
  There is a tensor product in $\scM(S)$, induced by the direct
product in ~$\V(S)$.   First, for $\alpha \in \Chow * {X_1}$ and $\beta
\in \Chow * {X_2}$ let
$$
 \alpha \times_S \beta := (\pr_1^*(\alpha)\cdot \pr_2^*(\beta)) \in
\Chow * {X\times_S Y} .
$$
 Next, for correspondences $\phi \in \Chow *{X_1\times_S X_2}$ and
$\psi\in \Chow*{X_3\times_S X_4}$, let
$$
 \phi\tensor_S \psi := t^*(\phi\times_S\psi) \in \Chow*{(X_1\times_S
X_3)\times_S (X_2\times_S X_4)} ,
$$
 where
$$
 t: (X_1\times_S X_3)\times_S (X_2\times_S X_4) \lra (X_1\times_S
X_2)\times_S (X_3\times_S X_4)
$$
 permutes the factors.  This determines the tensor product on morphisms,
and then the tensor product of two objects is given by
$$
 (X/S,p,n) \tensor (Y/S,q,m) := ((X\times_S Y)/S, p\tensor_S q, m+n) .
$$

\subsubhead 1.1.4. The direct sum in $\scM(S)$ \endsubsubhead
  There is also a direct sum in $\scM(S)$, induced by taking disjoint
union in ~$\V(S)$.  When $m=n$ it is defined by
$$
  (X/S,p, m) \oplus (Y/S, q, m) := (X\sqcup Y, p\oplus q, m) .
$$
 If $m < n$, say, then rewrite
$$
  (X/S,p,m) \isom (X/S,p,n)\tensor \LL_S^{n-m} =
(X\times_S(\PP^1_S)^{n-m}, p', n)
$$
 for some projector $p'$, and then the direct sum is defined by
$$
  (X/S,p,m)\oplus(Y/S, q, n) :=((X\times_S(\PP^1_S)^{n-m} \sqcup Y)/S,
p'\oplus q, n) .
$$

\subsubhead 1.1.5. $\scM(S)$ is pseudoabelian \endsubsubhead
 With these definitions it can be shown that $\scM(S)$ is a
$\QQ$-linear pseudoabelian tensor category.
An additive category is said to be pseudoabelian iff
for every object $M$ every idempotent $g\in\End_{\scM(S)}(M)$ has an
image, or equivalently a kernel, and the canonical map
$$
 (\Im(g) \oplus \Im(\id - g))\simto M
$$
 is an isomorphism.  See \cite{Jannsen, 1992} or \cite{Scholl, 1994,
Cor.3.5} to see that $\scM(S)$ is not in general an abelian category.

\subsubhead 1.1.6. The functor $\V(S) \to \scM(S)$ \endsubsubhead
  There is a natural contravariant functor from $\V(S)$ to $\scM(S)$,
given by associating to a morphism $f: X\to Y$ of smooth projective
$S$-schemes the class of the transpose of its graph, $\tgraph f \in
\Chow{d_S(Y)}{Y\times_S X}$, and associating to $X$ in $\V(S)$ the
object $(X,[\Delta(X)])$, where $\Delta(X)$ denotes the diagonal in
$X\times_S X$.  When $S = \spec k$ we write
$$
 h(X) :=  (X,[\Delta(X)]).
$$

\subsubhead 1.1.7. Some formulas \endsubsubhead
  For later use we note that for $X$, $Y$, $Z$ in $\V(S)$, and
$f:X\to Y$, $f' :Y\to X$, $g:Y\to Z$, $g':Z\to Y$, and
$\alpha \in \CH(X\times_S Y,\QQ)$ and $\beta \in \CH(Y\times_S Z,\QQ)$,
$$\spreadlines{\medskipamount}\alignat2
 [\Gamma_g]\circ \alpha &= (\id_X \times_S g)_*(\alpha)
 &[\tGamma_{g'}]\circ\alpha &= (\id_X\times_S g')^*(\alpha)  \\
 \beta\circ [\Gamma_f] &= (f\times_S\id_Z)^*(\beta)
 &\beta\circ [\tGamma_{f'}] &= (f'\times_S \id_z)_*(\beta)  \\
 [\Gamma_g]\circ [\Gamma_{f}] &= [\Gamma_{g\circ f}]
 &\qquad [\tGamma_{g'}] \circ [\tGamma_{f'}] &= [\tGamma_{f'\circ g'}]
\endalignat
$$
 see \cite{Deninger and Murre, 1991, 1.2.1}.  In particular, if
$f_1:X\to Y$, $f_2:X\to Y$, then
$$
 \graph{f_2} \circ \tgraph{f_1} = (f_1\times f_2)_*([\Delta(X)]) .
$$

\subsubhead 1.1.8. Remark on the relation between relative and absolute
motives \endsubsubhead
 When $S$ is projective the covariant functor $\V(S)\to
\V(k)$ taking $\lambda:X\to S$ to $\rho\circ\lambda: X\to \spec k$
induces a natural covariant functor $\Psi:\scM(S) \to \scM(k)$ which
makes the diagram
$$\CD
\V(S) @> \rho_* >> \V(k) \\
@VVV  @VV h V \\
\scM(S) @>> \Psi > \scM(k)
\endCD
$$
 commute.  Namely, let $i : X\times_S Y  \hra X\times Y$ be the inclusion,
and consider the morphism
$$
 i_* :\Chow{d_S(X)}{X\times_S Y} \to \Chow{\dim X}{X\times Y}.
$$
 Then it is easy to see that the codimensions work out so that a relative
correspondence of degree zero maps to an absolute correspondence of
degree zero, and it can also be checked that composition of relative
correspondences agrees with composition of absolute correspondences under
this
``pushing forward.''  Thus when $Y=X$ relative projectors map to absolute
projectors, and in this way we get a functor $\Psi$ as claimed.
Although this remark does not precisely apply to the situation of this
paper, it has been useful as part of the philosophy behind our methods;
see
also 3.2.7--3.2.8 and 3.3.7--3.3.9 below.

\subsubhead 1.1.9. The Chow groups of a Chow motive  \endsubsubhead
 Recall that in general a correspondence
$\gamma\in\Chow {}{X_1\times_k X_2}$ acts on a cycle class $[Z] \in\Chow
{}
{X_1}$ by
$$
 \gamma([Z]) := (\pr_2)_* (\pr_1^*(Z)\cdot \gamma) .
$$
  Then the Chow groups of $(X,p,m)$ in $\scM(k)$ are defined by
$$\spreadlines{1\jot} \align
\Chow j {(X,p,m)} :  &= p(\Chow{j+m} X ) \\
&= \Hom_{\scM(k)}(\LL_k^j,(X,p,m))
\endalign
$$
 and we let
$$
 \CH((X,p,m),\QQ) := \bigoplus_{j\in\ZZ} \Chow j {(X,p,m)} .
$$

\subsubhead 1.1.10. The cohomology groups of a Chow motive
\endsubsubhead
  In principle the cohomology groups of a Chow motive can be defined
with respect to any Weil cohomology theory, cf\. \cite{Kleiman, 1968,
1994}, but in this paper we will only consider Betti and \'etale
cohomology.  For a smooth, projective scheme $X$ over ~$k$, we write
$\Hdot^i(X,\Qdot)$ to signify either the Betti cohomology of
$X(\CC)^{\text{an}}$ with coefficients in $\QQ_B :=\QQ$, if $k$ comes
with an embedding into ~$\CC$, or the \'etale cohomology of
$X\times_{\spec k} \spec k^{\text{sep}}$ with coefficients in
$\QQ_\ell$; after taking Tate twists into account, cf\. \cite{Deligne,
1982, \sect 1}, we have
$$\spreadlines{1\jot}
 \Hdot^i(X,\Qdot (r)) :=\cases H^i_B(X(\CC)^{\text{an}}, \QQ_B(r)) , \\
   H^i_{\text{\'et}}(X\times_{\spec k} \spec k^{\text{sep}}, \QQ_\ell(r))
{}.
\endcases
$$
 Then the cohomology groups of $(X,p,m)$ in $\scM(k)$ are defined by
$$
\Hdot^i((X,p,m),\Qdot) := p( \Hdot^{i+2m}(X,\Qdot(m)) ).
$$
 Note that the $i^{\text{th}}$ cohomology group of $(X,p,m)$ has
weight~$i$, and for instance $\Hdot^i((X,p,m),\Qdot) \ne
\Hdot^i((X,p),\Qdot(m))$.  Let
$$
 \Hdot ((X,p,m),\Qdot) := \bigoplus_{i\in\ZZ} \Hdot^i((X,p,m),\Qdot) .
$$

\subhead  1.2. The conjectures \endsubhead
  Continuing to establish general terminology, as well as some of the
underlying motivation, we briefly recall the conjectures from
\cite{Murre, 1993} about the Chow groups of smooth projective varieties.
For more details and a summary of what is known, see {\sl op.~cit.};
for the relationship with conjectures of Beilinson, see \cite{Jannsen,
1994}.

\subhead 1.2.1.  Definition of Chow-K\"unneth decomposition
\endsubhead
 Let $X$ be a smooth projective variety of dimension~$d$.
A {\sl Chow-K\"unneth decomposition of $X$} is a
collection of mutually orthogonal projectors $\pi_0(X), \ldots ,
\pi_{2d}(X)$ in $\Chow d {X\times X} = \Corr^0(X,X)$ such that
$$
\sum_{i=0}^{2d} \pi_i(X) = [\Delta(X)] ,
$$
 and
$$
 \pi_i(X)(\Hdot (X,\Qdot)) = \Hdot^i(X,\Qdot) .
$$
 When a Chow-K\"unneth decomposition of $X$ exists, we let
$$
 h^i(X) := (X,\pi_i(X)).
$$

\proclaim{Conjecture A}
 For any smooth projective variety $X$ there exists a Chow-K\"unneth
decomposition of ~$X$.
\endproclaim

\proclaim{Conjecture B}
 Let $X$ be a smooth projective variety, and assume that there exists a
Chow-K\"unneth decomposition of ~$X$.  Then
$$
 \CH^j(h^i(X),\QQ):= \pi_i(X)(\Chow j X ) = 0 \quad  \text{for }i<j
\text{ or }i>2j .
$$
\endproclaim

\subsubhead 1.2.2. A filtration on the Chow groups of ~$X$ \endsubsubhead
  Let $X$ be a smooth projective variety, and assume that there exists a
Chow-K\"unneth decomposition of ~$X$ such that $\CH^j(h^i(X),\QQ)=0$ for
$i<j$ or $i>2j$.  Then there is a $j$@-step filtration on $\Chow j X$
defined by
$$\spreadlines{1\jot}\aligned
F^\nu\CH^j(X,\QQ) :&= \Ker\{\pi_{2j-\nu+1} \big|
F^{\nu-1}\CH^j(X,\QQ) \}  \\
  &= \bigoplus_{i=j}^{2j-\nu} \CH^j(h^i(X),\QQ) ,
\endaligned
$$
 for $0\le \nu\le j$.

\proclaim{Conjecture C}
 Let $X$ be a smooth projective variety, and assume that there exists a
Chow-K\"unneth decomposition of ~$X$ such that $\CH^j(h^i(X),\QQ)=0$ for
$i<j$ or $i>2j$.  Then the filtration $F^\nu\CH^j(X,\QQ)$ is independent
of the choice of Chow-K\"unneth projectors $\pi_i(X)$.
\endproclaim

\subsubhead 1.2.3. The cycle class map \endsubsubhead
  Let
$$
 \Chowhom j X :=\Ker (\cl : \Chow j X \to \Hdot^{2j}(X,\Qdot(j))) ,
$$
 where $\gamma$ is the cycle class map.  Then it follows from the
commutative diagram
$$\CD
 \Chow j X @> \pi_{2j}(X) >> \Chow j X \\
@V \cl VV  @VV \cl V \\
 \Hdot^{2j}(X,\Qdot(j)) @> \sim > \pi_{2j}(X) > \Hdot^{2j}(X,\Qdot(j))
\endCD
$$
 that
$$
 F^1\Chow j X := \Ker(\pi_{2j}(X)\big| \Chow j X) \subseteq
\CH_{\text{hom}}^j(X,\QQ) .
$$

\proclaim{Conjecture D}
  Let $X$ be a smooth projective variety, and assume that there exists a
Chow-K\"unneth decomposition of ~$X$ such that $\CH^j(h^i(X),\QQ)=0$ for
$i<j$ or $i>2j$.  Then
$$
F^1\Chow j X = \CH_{\text{hom}}^j(X,\QQ) \quad \text{for }1\le j\le
\dim(X).
$$
\endproclaim

\subsubhead 1.2.4. A generalization \endsubsubhead
 Suppose $M=(X,p,m)\in\scM(k)$ is a Chow motive, with $X$
equidimensional of dimension ~$d$.  Then one can define a
Chow-K\"unneth decomposition of $M$ as a collection of mutually
orthogonal projectors $\pi_i(M)\in \End_{\scM(k)}(M):= \Corr^0(M,M)$,
with $-2m\le i\le 2d-2m$,
such that $\sum_{i} \pi_i(M) = \id_M  =p$ and
$\pi_i(M)(\Hdot (M,\Qdot)) = \Hdot^i(M,\Qdot)$.  It might sometimes
be useful, as it is for us below in sections four through six, to
decompose a variety as a sum of submotives in some other way than a
Chow-K\"unneth decomposition, and then verify the conjectures on the
various submotives, in the sense of the following lemma.

\proclaim{Lemma 1.2.5}
Suppose $M \simeq M_1 \oplus M_2$ in $\scM(k)$.
\roster
\item If a Chow-K\"unneth decomposition of ~$M_1$ exists and a
Chow-K\"unneth decomposition of ~$M_2$ exists, then a Chow-K\"unneth
decomposition of ~$M$ exists.
\item If in addition $\pi_i(M_t)(\Chow j {M_t}) =0$ whenever $i<j$ or
$i>2j$ for both $t=1$ and $t=2$, then with the induced Chow-K\"unneth
decomposition $\pi_i(M)(\Chow j {M}) =0$ whenever $i<j$ or $i>2j$.
\item If in addition
$$
\Ker\big(\pi_{2j}(M_t)\big| \Chow j {M_t}\big) = \Ker \{\cl:\Chow j
{M_t} \to \Hdot^{2j}(M_t,\Qdot(j))\}
$$
 for both $t=1$ and $t=2$, then with the induced Chow-K\"unneth
decomposition
$$
\Ker\big(\pi_{2j}(M)\big| \Chow j {M}\big) = \Ker \{\cl:\Chow j {M} \to
\Hdot^{2j}(M,\Qdot(j))\} .
$$
\endroster
\endproclaim

\demo{Proof}
 \therosteritem1
 Let $M_1 = (X_1, p_1, m_1)$ and $M_2 = ( X_2, p_1, m_2)$, an suppose
first for
simplicity that $m_1 =m_2 =:m$, say, so that by definition~1.4
$$
M \simeq (X_1\sqcup X_2, p_1\oplus p_2, m).
$$
   Then the inclusions $j_1$ and $j_2$ of $X_1$ and $X_2$ respectively
into $X_1\sqcup X_2$ induce orthogonal central idempotents, say $e_1$
and $e_2$, whose sum is the identity in $\End_{\scM(k)}(M_1\oplus M_2)$.
Therefore $M_t \isom (X_1\sqcup X_2, e_t^*(p_1\oplus p_2), m)$, $t=1,2$.
So if $\id_{M_t} = \sum_i \pi_i(M_t)$ is a Chow-K\"unneth decomposition
for
$M_t$, then (up to isomorphism)
$$
 \id_{(M_1\oplus M_2} = \sum_i \big( e_{1*}\pi_i(M_1) + e_{2*}\pi_i(M_2)
\big)
$$
  is a Chow-K\"unneth decomposition for $M$.  In case $m_1 < m_2$, say,
then as
in 1.4 we have
$$
 M_1 \simeq M_1' := (X_1 \times (\PP^1)^q, p_1' , m_2),
$$
 for a suitable choice of ~$p_1'$ and $q:= m_2-m_1$.  Then $M\simeq M_1'
\oplus
M_2$, so it suffices to know that the existence of a Chow-K\"unneth
decomposition for $M_1$ implies the existence of a Chow-K\"unneth
decomposition
for ~$M_1'$.  However, the isomorphism $M_1 \simeq M_1'$ can be used to
transform a Chow-K\"unneth decomposition of $M_1$ into a Chow-K\"unneth
decomposition of $M_1'$ with $\pi_i(M_1') \simeq \pi_{i-2q}(M_1)$.

 \therosteritem2
  As in ~\therosteritem1, first suppose $m_1=m_2=:m$.  Then from 1.1.9 we
see
that
$$\align
 \CH^j(M,\QQ) &\simeq \Hom_{\scM(k)}(\LL^j, M_1\oplus M_2)\\
 &\isom \Hom_{\scM(k)}(\LL^j, M_1) \oplus \Hom_{\scM(k)}(\LL^j, M_1) \\
 & = \CH^j(M_1,\QQ) \oplus \CH^j(M_2,\QQ) .
\endalign
$$
 Thus if $\pi_i(M_t)(\Chow j {M_t}) =0$ whenever $i<j$ or $i>2j$ for both
$t=1$ and $t=2$, then the same must be true for $M$ as well.  Now if $m_1
<
m_2$, say, then we need to know that $\pi_i(M_1')(\Chow j{M_1'}) =0$
whenever
$i<j$ or $i>2j$, with $M_1'$ as above.  So consider the diagram
$$\CD
 \Chow j{M_1'} @= \Chow{j-q}{M_1} \\
 @V \pi_i(M_1') VV  @VV \pi_{i-2q}(M_1) V \\
 \Chow j{M_1'} @= \Chow{j-q}{M_1}
\endCD
$$
 where the equalities follow from ~1.1.9.  Since $q>0$, if $i<j$ then
$i-2q <
j-q$ and if $i>2j$ then $i-2q>2(j-q)$, so $M_1'$ satisfies the hypothesis
of
{}~\therosteritem2 whenever $M_1$ does, as required.

 \therosteritem3
 When $m_1 =m_2$, then similarly as above we see that the
cycle class map $\cl:\CH^j(M,\QQ) \to \Hdot^{2j}(M,\Qdot(j))$ is the
direct sum of the two cycle class maps $\cl:\CH^j(M_t,\QQ) \to
\Hdot^{2j}(M_t,\Qdot(j))$, for $t=1,2$, and the claim follows directly.
And if
$m_1 < m_2$, then the diagram in ~\therosteritem2 above with $i=2j$ can be
combined with the diagram in~1.2.3 to show that $M_1'$ satisfies the
hypothesis
of ~\therosteritem3 whenever $M_1$ does, as required.  This completes the
proof
of the lemma.
\Qed
\enddemo

\head \bf 2.  Elliptic modular surfaces and threefolds \endhead
\noindent
   We review the geometric structure of elliptic modular surfaces and
threefolds with level-$N$ structure.  To begin, we fix an integer $N \ge
3$ once and for all, and a ground field ~$K$ in which $2N$ is invertible
and which contains $N^{\text{th}}$ roots of unity.  When there is no
danger of confusion we will drop $N$ or $K$ from the notation.

\subhead 2.1. The elliptic modular curve \endsubhead
    Let $M:=M_N$ be the elliptic modular curve over $K$ that represents
the functor which to a $K$-scheme ~$S$ associates the set of isomorphism
classes of elliptic curves $\scE/S$ with level-$N$ structure, where a
level-$N$ structure consists of an isomorphism
$$
 \alpha: \ZmodNsq \times S \simto \scE[N]/S
$$
 of group schemes over~$S$, compare \cite{Deligne and Rapoport, 1973,
Ch.IV} or \cite{Katz and Mazur, 1985, Ch.III}.  If $K$ is a subfield of
$\CC$, the analytic space $\M^{\text{an}}(\CC)$ associated to $\M$ is
isomorphic to $\Gamma(N)\bs\frH$, where $\Gamma(N)\subset\SL_2(\ZZ)$ is
the subgroup of matrices congruent to the identity modulo~$N$.  A smooth
completion of $M$
$$
 j: \M \hra \Mbar
$$
 is obtained by adjoining a finite set of cusps
$$
\Minf := \Mbar - \M
$$
 which parameterize generalized elliptic curves.

\subhead 2.2. The elliptic modular surface  \endsubhead
    Since $N\ge 3$ there
exists a universal elliptic curve with level~$N$ structure $\phi:\E\to
\M$.  Then the universal generalized elliptic curve with level-$N$
structure $\phibar :\Ebar\to\Mbar$ is the canonical minimal smooth
completion of $\phi : \E\to \M$ \cite{Shioda, 1972}, \cite{Deligne and
Rapoport, 1973}.  Let
$$
 \alphabar :  \ZmodNsq \times \Mbar \simto \Ebar .
$$
 denote the extension of the level-$N$ structure to $\Ebar$.
 The N\'eron model $\Estar$ of $\E$ over $\Mbar$ is the open subscheme
where $\phibar$ is smooth.  The following diagram summarizes the
notation.
$$
\CD
\E @CCC \Estar @CCC \Ebar @DDD \Einf \\
 @V \phi VV  @V \phistar VV  @VV \phibar V  @VVV \\
 \M @CC j C \Mbar @= \Mbar @DDD \Minf \\
\endCD
$$

\subsubhead 2.2.1. Description of $\Einf$  \endsubsubhead
  For $c\in\Minf$, the fibres
$$
 \Ebar_c :=\phibar^{-1}(c) \simeq \ZmodN \times \PP^1
$$
 are standard N\'eron $N$-gons, where we can number the components by
letting $\theta_c(m) \simeq \PP^1$ be the component containing
$\alphabar((m,n),c)$, for $(m,n) \in \ZmodNsq$.  Note that for fixed
$m$, as $n$ varies the $\alphabar((m,n),c)$ all lie in the same
component, and may be identified with $N^{\text{th}}$ roots of unity
when $\theta_c(m)$ minus its intersections with $\theta_c(m-1)$ and
$\theta_c(m+1)$ is identified with $\Gm$.  Sometimes we refer
to $\theta_c(0)$ as the identity component.  In this notation the
intersection relations among the components of ~$\Einf$ are
$$\spreadlines{1\jot}
 (\theta_c(m)\cdot \theta_{c'}(m')) = \cases
  -2 &\text{if  } c=c' \text{ and } m=m' \\
 \hphantom{-} 1 &\text{if  } c=c' \text{ and } m-m' =\pm 1 \\
 \hphantom{-} 0 &\text{otherwise}
\endcases
$$
  \cite{Kodaira, 1963, III}, \cite{Shioda, 1972}, \cite{Ash et al., 1975,
I.4}.  In particular, the rank of intersection matrix for the
components of the fibre over a cusp is ~$(N-1)$.

\subsubhead Remark 2.2.2 \endsubsubhead
 It follows from \cite{Shioda, 1972, Thm.1.1} that a basis for
$\NS(\Ebar)\tensor \QQ$ is given by the zero-section $\ebar
:=\alphabar((0,0),\Mbar)$, a regular fibre, and the components of the
cusp fibres other than the identity component.

\subhead 2.3.  The elliptic modular threefold \endsubhead
 Consider the fibre products
$$\spreadlines{2\jot} \align
\phitwo : \Etwo := \E\times_\M  \E &\lra \M \\
\phitwostar : \Etwostar := \Estar \times_\Mbar  \Estar &\lra \Mbar\\
\phitwobar : \Etwobar := \Ebar\times_\Mbar  \Ebar &\lra \Mbar\\
 \Etwobarinf := \Einf \times_\Minf  \Einf &\lra \Minf .
\endalign
$$
  Then $\Etwobar$ is not smooth:  Using the local coordinates of
\cite{Deligne, 1969, Lemme~5.5} or \cite{Scholl, 1990, \sect2}, compare
also
\cite{Schoen, 1986}, one can check that the points over $c\in\Minf$ that
are a product of two double points of $\Ebar_c$ are rational double
points in ~$\Etwobar$.  If we let $\Etwobarinf_0 =
\Etwobar^{\text{sing}} \subset \Etwobarinf$ denote the reduced subscheme
of $\Etwobar$ consisting of all these points, for all $c\in\Minf$, then
applying \cite{Deligne, 1969, Lemmes ~5.4, 5.5} or \cite{Scholl, 1990,
Prop.2.1.1, Thm.3.1.0(i)} gives us the following description of the
desingularization ~$\Etwotild$ of ~$\Etwobar$.

\proclaim{Proposition 2.3.1}
  Let
$$
 \beta : \Etwotild \lra \Etwobar
$$
 be the blowing-up of $\Etwobar$ along $\Etwobarinf_0$.  Then
$\Etwotild$ is nonsingular.  Further, let
$$
 \Etwotildinf := (\phitwobar \circ \beta)\inv(\Minf)
$$
 be the union of the resulting fibres over $\Minf$.  Then $\Etwotildinf$
consists of $2N^2\cdot\#(\Minf)$ components, half of which are
quadric surfaces (isomorphic to $V(xy-zw) \subset \PP^3$) that are
the components of the exceptional divisor, and half of which are
isomorphic to $\PP^1\times\PP^1$ with four (smooth) points blown ~up,
these being the proper transforms with respect to ~$\beta$ of the
components of $\Etwobarinf$.  In particular, all the components of
$\Etwotildinf$ are rational surfaces.
\endproclaim

\subsubhead Remark 2.3.2 \endsubsubhead
  In fact $\#(\Minf) = \frac 1 2  N^2 \prod_{p|N} (1- p^{-2})$
\cite{Miyake, 1989}, though this will play no explicit role for us.

\subsubhead 2.3.3. Notation \endsubsubhead
 As a matter of notation, let
$$
\phitwotild := \phitwobar \circ \beta : \Etwotild \lra \Mbar
$$
 be the fibre structure map.   The following diagram then summarizes the
rest of the notation.
$$\CD
 \Etwo @CCC \Etwostar @CCC \Etwotild  @DDD \Etwotildinf \\
 @|  @|  @V \beta VV  @VVV\\
 \Etwo @CCC \Etwostar @CCC \Etwobar  @DDD \Etwobarinf \\
@V \phitwo VV  @V \phitwostar VV @V \phitwobar VV  @VVV \\
 \M @C j CC \Mbar @= \Mbar  @DDD \Minf
\endCD
$$

\subsubhead 2.3.4. Indexing the components of $\Etwotildinf$
\endsubsubhead
   For use later (in 3.3.11 and 4.4.1) we also index the components
$\Theta_c$ of the cusp fibres $\Etwotild_c$, for $c\in\Minf$.  According
to Proposition ~2.3.1, half the components are the proper transforms
of the components $\theta_c(m)\times_{\{c\}} \theta_c(n)$ of
$\Etwobar_c$, so these $\Theta_c(m,n)$ are naturally indexed by pairs
$(m,n)\in\ZmodNsq$.  The remaining components come from blowing up
points which can be described as the (fibre) product (over ~$c$ in
$\Minf$) of the point where $\theta_c(m)$ intersects $\theta_c(m+1)$ with
the point where $\theta_c(n)$ intersects $\theta_c(n+1)$, as $m$ and $n$
run over $\ZmodN$.  Then the correct incidence relations and symmetries
are best described if we call the blowing-up of this point
$\Theta_c(m+\frac12,n+\frac12)$, indexed by a pair of half-integers
mod~$N\ZZ$; compare \cite{Deligne and Rapoport, 1973, \sect VII.1}.

\head \bf 3. Construction of projectors for $\Ebar$ and $\Etwotild$
\endhead
\nobreak
\subhead 3.1. Introduction to the construction \endsubhead
  When $S$ is a smooth, connected, quasi-projective base scheme and
$A\to S$ is an abelian scheme, then there exist canonical, mutually
orthogonal relative projectors $\pican i (A/S)$ in
$\Chow{d_S(A)}{A\times_S A}$ whose sum is the diagonal \cite{Shermenev,
1974}, \cite{Deninger and Murre, 1991}, \cite{K\"unneman, 1994}.  These
are characterized by the property that
$$
 \tgraph{\mu(n)} \circ \pican i (A/S) = \pican i (A/S) \circ
\tgraph{\mu(n)} = n^i\, \pican i (A/S) ,
\tag 1
$$
 where $\mu(n): A\to A$ is the multiplication by $n$ endomorphism of
$A/S$.  In particular, when $S$ is a point, these $\pican i(A)$ define
a Chow-K\"unneth decomposition of ~$A$.

  In our situation, even though $\Ebar\to\Mbar$ and $\Etwotild\to\Mbar$
are not abelian schemes, one of the underlying ideas for the projectors
we define in this section is to extend in a suitable
sense the canonical relative projectors for $\E/\M$ and $\Etwo/M$ to
projectors for $\Ebar$ and $\Etwotild$, respectively.  The idea is that
$E\times_M E$ naturally embeds in $\Ebar \times \Ebar$, and this
embedding factors through the natural embedding of $\Ebar\times_\Mbar
\Ebar$ into $\Ebar\times\Ebar$; and likewise the natural embedding of
$\Etwo\times_M\Etwo$ into $\Etwotild\times\Etwotild$ factors through
$\Etwotild\times_\Mbar \Etwotild$.  Then what we would like to do is
``push forward'' $\pican i(\E/\M)$ from $\Chow 1 {E\times_M E}$ to $\Chow
2{\Ebar\times\Ebar}$, and similarly ``push forward'' $\pican i(\Etwo/M)$
from $\Chow 2{\Etwo\times_M\Etwo}$ to $\Chow
3{\Etwotild\times\Etwotild}$.  The trouble is, there is no natural push
forward for this situation, and the only alternative seems to be to
choose an explicit cycle to represent the rational equivalence class
$\pican i(\E/\M)$, then take its closure in
$\Ebar\times_\Mbar \Ebar$, and then push that forward to a cycle
on $\Ebar\times\Ebar$; and likewise for $\Etwotild\times\Etwotild$.
Conceptually this is what we do, but as a matter of logical presentation
it seems preferable to begin by describing explicit cycles supported on
$\Ebar\times_\Mbar\Ebar$, and then show that they have nice properties as
mutually orthogonal projectors in $\Chow 2{\Ebar\times\Ebar}$.  The point
is that for technical reasons these cycles on $\Ebar\times_\Mbar\Ebar$
are not simply the closures of the obvious ``natural'' representatives for
$\pican i(\E/\M)$, as described in \cite{K\"unnemann, 1994} for
example.  So it requires some work to show that their restrictions to
$\E\times_M\E$ do indeed represent the canonical relative projectors for
$E/M$, see Proposition~3.2.8.

  Then for $\Etwotild$ we use the previously-defined cycles on
$\Ebar\times_\Mbar\Ebar$ in the definition of explicit cycles supported on
$\Etwotild\times_\Mbar \Etwotild$ that behave nicely as mutually
orthogonal projectors in $\Chow 3 {\Etwotild\times\Etwotild}$.  One of
the ideas underlying this construction is to take advantage of the
relative product structure of $\Etwo/M$, for in this way we get nine
projectors corresponding (in the sense of 3.3.7--3.3.9 below) to
the $\pican i (\E/\M)
\tensor_M \pican j (\E/\M)$, for $0\le i,j\le 2$, in $\Chow
2{\Etwo\times_M\Etwo}$, rather than just the five that correspond to the
$\pican i (\Etwo/\M)$, for $0\le i\le 4$.

\subhead 3.2. Extending canonical relative projectors to ~$\Ebar$
\endsubhead
\par\nobreak
\subsubhead 3.2.1. The zero-section and its transpose \endsubsubhead
  Let $\e := \alpha((0,0), M)$ be the zero-section as curve in ~$E$.
Then
$$\spreadlines{1\jot}
{\aligned
 [\E\times_\M \e] &= \pican 2(\E/\M), \\
 [\e\times_\M \E] &= \pican 0 (\E/\M) = \tr{\vphantom{\pi}} \pican
2(\E/\M)
\endaligned} \qquad \text{\sl in } \Chow 1{\E\times_\M\E}
$$
 \cite{K\"unneman, 1994, 4.1.2(iv)}.  Now let $\ebar :=
\alphabar((0,0),\Mbar)$ be the zero-section as curve in ~$\Ebar$.  Then
$$\spreadlines{1\jot}
{\aligned
 \pbar_2 &:= [\Ebar\times_\Mbar\ebar], \\
 \pbar_0 &:= [\ebar\times_\Mbar \Ebar] = \tr{\vphantom{p}}\pbar_2
\endaligned} \qquad  \text{\sl in } \Chow 2{\Ebar\times\Ebar}
$$
 are projectors, but unexpectedly, they are not orthogonal, as the next
lemma
explains.   In order to formulate this lemma precisely, and also for later
purposes (see~3.2.7), we consider the inclusions
$$
\psi: E\times_M E @C \psi_1 CC \Ebar\times_\Mbar \Ebar @C \psi_2 CC
\Ebar\times\Ebar .
$$

\proclaim{Lemma 3.2.2}
 In $\Chow 2{\Ebar\times\Ebar}$
\roster
\smallskip
\item $\pbar_0\circ\pbar_0 =\pbar_0$ and $\pbar_2\circ\pbar_2 =\pbar_2$;
\smallskip
\item $\pbar_2\circ\pbar_0 =0$;
\smallskip
\item $\pbar_0\circ\pbar_2 = (\psi_2)_*(\phibar\times_\Mbar\phibar)^{*}
    \phibar_*[\ebar \cdot \ebar] \ne 0$, where $[\ebar \cdot \ebar]$
denotes the self-intersection cycle in $\Chow 2 {\Ebar}$.
\endroster
\endproclaim

\demo{Proof}
 Let $\mubar(0) := \alphabar((0,0),\phibar(\bigdot)) :\Ebar\to
\Ebar$ be the morphism given by projection onto the zero-section.
(The notation is meant to suggest ``multiplication by zero,'' extending
to $\Ebar$ of the fibre-wise group homomorphism that maps everything to
the identity element.)  Then $\pbar_2$ and $\pbar_0$ correspond to the
graph and transposed graph of $\mubar(0)$, respectively,
$$
 \pbar_2 = \graph{\mubar(0)}, \qquad\quad \pbar_0 = \tgraph{\mubar(0)}.
\tag 4
$$
 Then \therosteritem1 follows because $\mubar(0)\circ\mubar(0)=
\mubar(0)$, and \therosteritem2 because by ~1.1.7
$$
 \graph{\mubar(0)}\circ\tgraph{\mubar(0)} = (\mubar(0)\times\mubar(0))_*
([\Delta(\Ebar)])  ,
$$
 which vanishes in $\Chow 2{\Ebar\times\Ebar}$ for dimension reasons.
As for ~\therosteritem3, we verify this by direct computation.  In order
to
have proper intersection for this computation, we move the graph
$\graph{\mubar(0)}$ on $\Ebar\times\Ebar$ by first moving the divisor
$\ebar$
in its linear equivalence class on $\Ebar$ to a divisor $\ebar'$
intersecting
$\ebar$ properly on $\Ebar$ (and moreover, for simplicity, also such that
over
a cusp $\ebar'$ passes through neither the crossing points of the
components of
that fibre nor through the intersection of $\ebar$ with the fibre).  Also
note
that the cycle class we finally get is the class of a cycle supported on
the
singular variety $\Ebar\times_\Mbar\Ebar$ and therefore we have to go via
$\CH_1(\Ebar\times_\Mbar\Ebar, \QQ)$ (in the sense of \cite{Fulton,
1988}).
The nonvanishing is a
consequence of the fact that the self-intersection number
$(\ebar\cdot\ebar) = -(p_a +1) <0$ \cite{Kodaira, 1963, p.15},
\cite{Shioda, 1972, p.25}.
\Qed
\enddemo

\subhead 3.2.3. Definition of $\pibar0$ and $\pibar2$ \endsubhead
 If we now let, in $\Chow 2{\Ebar\times\Ebar}$,
$$\spreadlines{1\jot} \align
 \pibar0 &:= \pbar_0 - \tfrac 1 2 \pbar_0 \circ \pbar_2 =
\tgraph{\mubar(0)} -\tfrac12 \tgraph{\mubar(0)}\circ\graph{\mubar(0)}
\\
 \pibar2 &:= \pbar_2 - \tfrac 1 2 \pbar_0 \circ \pbar_2 =
\graph{\mubar(0)} -\tfrac12 \tgraph{\mubar(0)}\circ\graph{\mubar(0)}
 = \tr{\vphantom{\pi}}\pibar0
\endalign
$$
 then it follows from Lemma~3.2.2 that these are orthogonal
projectors.  For use below we also choose a zero cycle $\frak a$ on
$\Mbar$ representing $\phibar_*[\ebar\cdot\ebar]$, i.e.,
$$
 [\fra] = \phibar_*[\ebar\cdot\ebar] \in \Chow 1\Mbar,
$$
 and observe that by
doing so we get representative cycles for $\pibar0$ and $\pibar2$
supported on $\Ebar\times_\Mbar\Ebar$.  Also
for later reference note that the ``correction term''
$$
 \tfrac 1 2 \pbar_0 \circ \pbar_2 = \tfrac12
\tgraph{\mubar(0)}\circ\graph{\mubar(0)} =
\tfrac12(\psi_2)_*(\phibar\times_\Mbar\phibar)^*([\frak a])
$$
 is nilpotent of order~2 in $\Chow 2{\Ebar\times\Ebar}$.

\subsubhead 3.2.4. Automorphism correspondences on $\Ebar$
\endsubsubhead
  Following \cite{Scholl, 1990} we consider a group of automorphisms
acting on ~$\Ebar$.   Firstly, for $\b\in (\ZZ/N\ZZ)^2$ translation by
$\alpha(\b,z)$ in each fibre $\phi\inv(z)$ defines an automorphism
$\tau(\b): E\to E$ of ~$E$.  Since this depends only on the group
structure of $E/M$, it extends first to an automorphism $\tau^*(\b)$ of
$\Estar$, and then by Zariski's Main Theorem \cite{Hartshorne, 1977,
V.5.2, p.410}, since the invertibility of $\tau^*(\b)$ away from the
isolated points of $\Ebar -\Estar$ precludes the total transform of any
of these points in the closure of the graph of $\tau^*(\b)$ having
dimension one or more, to an automorphism $\taubar(\b) : \Ebar \to
\Ebar$.  In this way we get a group action of $(\ZZ/N\ZZ)^2$ on ~$\Ebar$.
  By the same reasoning, the fibrewise inversion map is an automorphism
of $\E$ that extends first to an automorphism of $\Estar$ and then to an
automorphism $\mubar(-1):\Ebar\to\Ebar$ of ~$\Ebar$, and together with
the identity map this gives a group action of $\mu_2$ on $\Ebar$.  These
two group actions together give a group action of the semidirect
product
$$
 \G := \ZmodNsq \rtimes \mu_2
$$
 on $\Ebar$, which can be extended $\QQ$-linearly to define an action of
the group ring $\QQ[\G]$ on ~$\Ebar$.  In particular, by associating to
a group element $g\in\G$ the class of its graph ~$\graph g$
(respectively, transposed graph ~$\tgraph g$) in $\Chow 2
{\Ebar\times\Ebar}$, we get a $\QQ$-algebra homomorphism
$$
  \QQ[\G] \lra \Chow 2 {\Ebar\times\Ebar}
$$
 (respectively, antihomomorphism $\QQ[\G]^{\text{opp}} \to \Chow 2
{\Ebar\times\Ebar}$) from the group ring of $\G$ into the ring of
degree-zero correspondences on ~$\Ebar$.  Further, since the group
actions operate fibrewise, these correspondences are supported
on $\Ebar\times_\Mbar\Ebar$.  We remark also that for automorphisms of
$\Ebar$ such as those defined by the action of $g\in\G$,
$$
 \tgraph g = \graph {g\inv} .
$$

\subhead 3.2.5. Definition of $\pibar1$ \endsubhead
 We take for $\pibar1$ the projector $\Pi_\varepsilon$ defined in
\cite{Scholl, 1990, 1.1.2} for $k=1$, which may be described as
follows:  Let $\varepsilon = \varepsilon_1$ be the character of ~$\G$
defined by the product of the trivial character on $\ZmodNsq$ and the
sign character on $\mu_2$; then one description of $\pibar 1$ is
$$
 \pibar 1 := \pibar {\varepsilon} = \frac 1 {2N^2} \sum_{g\in\G}
\varepsilon(g)\inv  \graph g .
$$
 As the homomorphic image of an idempotent in $\QQ[\G]$, it follows that
$\pibar1$ is a projector in $\Chow 2 {\Ebar\times\Ebar}$, and it is also
clear that $\tr{\vphantom{\pi}}\pibar1 = \pibar1$.  Another description
of $\pibar1$ comes from observing that
$$
 \lambdabar :=  \frac 1 2 \big( \graph{\mubar(1)} - \graph{\mubar(-1)}
\big), \qquad
 \thetabar := \frac 1 {N^2} \sum_{b\in(\ZZ/N\ZZ)^2} \graph{\taubar(b)}
$$
 are homomorphic images of commuting idempotents in $\QQ[\G]$, and then
$$
 \pibar1 = \lambdabar \circ \thetabar = \thetabar \circ \lambdabar.
$$

\proclaim{Proposition 3.2.6}
  The $\pibari$, for $i=0,1,2$, are mutually orthogonal projectors in
$\CH^2(\Ebar\times\Ebar,\QQ)$.
\endproclaim

\demo{Proof}
  We have already seen the idempotency of each $\pibari$, and the
orthogonality of $\pibar0$ and $\pibar2$, so it only remains to check
that $\pibar1$ is orthogonal to the other two.  To see this, we can use
3.2.3 that
$$
\pibar2 := \graph{\mubar(0)} - \frac 1 2
\tgraph{\mubar(0)}\circ\graph{\mubar(0)} = \tr{\vphantom{\pi}}\pibar0 .
$$
  Then from the observation that
$$
 \mubar(0)\circ \mubar(\pm1) = \mubar(\pm1)\circ \mubar(0) =\mubar(0)
$$
 and the formulas ~1.1.7 it follows immediately that
$\lambdabar$ is orthogonal to both $\pibar2$ and $\pibar0$, and thus
$\pibar1$ is as well.
\Qed
\enddemo

\subhead 3.2.7. Notations and definitions related to cycles on
$\Ebar\times_\Mbar\Ebar$ \endsubhead
 Suppose $\alpha : \Ebar\to\Ebar$ is a morphism such that $\alpha$
respects the fibre structure of $\Ebar\to\Mbar$.  Then the graph
$\Gamma_\alpha$ of $\alpha$ is supported on $\Ebar\times_\Mbar\Ebar$.  In
order to emphasize this we may write $\Gamma_\alpha^{\text{rel}}$ for the
graph of $\alpha$ as a cycle on $\Ebar\times_\Mbar\Ebar$, and $\graph
\alpha ^{\text{rel}} := [\Gamma_\alpha^{\text{rel}}]$ for its class in
$\CH_2(\Ebar\times_\Mbar\Ebar,\QQ)$ (the Chow group in the sense of
\cite{Fulton, 1984}, since $\Ebar\times_\Mbar\Ebar$ is singular).  Now
consider
again the inclusions
$$
\psi: E\times_M E @C \psi_1 CC \Ebar\times_\Mbar \Ebar @C \psi_2 CC
\Ebar\times\Ebar .
$$
  Then $\graph \alpha = (\psi_2)_*(\graph \alpha^{\text{rel}})$.
By abuse of notation define
$$
 \psisharp(\graph\alpha) := \psi_1^*(\graph\alpha^{\text{rel}}) \qquad
\text{in
} \Chow 1{E\times_M E} .
$$
  Then
$\psisharp(\graph\alpha)$ is just the class in $\Chow 1 {E\times_M E}$ of
the graph of the restriction of ~$\alpha$ to ~$E$.  Therefore if the
morphism $\beta:\Ebar\to\Ebar$ also respects the fibre structure of
$\Ebar\to\Mbar$ then we have
$$
 \psisharp(\graph{\alpha\circ\beta}) = \psisharp(\graph{\alpha}) \circ
\psisharp(\graph{\beta}) ,
$$
 and, if we allow the same notations and definitions for the transpose of
a graph, also
$$
 \psisharp(\tgraph{\alpha\circ\beta}) = \psisharp(\tgraph{\beta}) \circ
\psisharp(\tgraph{\alpha}) .
$$

  Now if we extend these notations and definitions by linearity and apply
them to the cycles and projectors defined in 3.2.5, then we have
$$
 \pibarrel1 = \frac 1 {2N^2} \sum_{g\in\G} \varepsilon(g)^{-1} \graph g
\rel
$$
 in $\CH_2(\Ebar\times_\Mbar\Ebar,\QQ)$, and $\lambdabar\rel$ and
$\thetabar\rel$ may be defined similarly.  We may also apply these
notations and definitions to the cycles and projectors in ~3.2.3,
and let
$$\spreadlines{1\jot} \align
\pibarrel0 &:= \tgraph{\mubar(0)}\rel -\tfrac12
[(\phibar\times_\Mbar\phibar)^*(\fra)] \\
\pibarrel2 &:= \graph{\mubar(0)}\rel -\tfrac12
[(\phibar\times_\Mbar\phibar)^*(\fra)] .
\endalign
$$
 Then for $i=0,1,2$,
$$
 \pibar i = (\psi_2)_*\pibarrel i ,
$$
 and thus we have elements
$$
 \psisharp\pibar i := \psi_1^*\pibarrel i \in \Chow 1{E\times_M E} ,
$$
 and we get $\psisharp(\lambdabar)$ and $\psisharp(\thetabar)$
similarly.

\proclaim{Proposition 3.2.8}
 With the notation as above,  in $\CH^1(E\times_M E,\QQ)$
$$
 \psisharp\pibari = \pican i(E/M)
$$
  for $0\le i\le 2$.
\endproclaim

\demo{Proof}
 Consider $i=1$ first.  Then it follows immediately from the
considerations in ~3.2.7 that
$$
 \psisharp(\pibar1) = \psisharp(\lambdabar)\circ \psisharp(\thetabar) =
\psisharp(\thetabar) \circ \psisharp(\lambdabar) .
$$
  Then in $\Chow 1 {E\times_M E}$ we have
$$\spreadlines{1\jot}\align
\psi_1^*(\thetabar\rel)  \circ\graph{\mu(N)} &= \frac 1{N^2}
  \sum_{b\in\ZZ/N\ZZ)^2} \tgraph{\mu(N)}\circ\tgraph{\tau(b)} \\
  & = \frac 1{N^2} \sum_{b\in\ZZ/N\ZZ)^2} \tgraph{\mu(N)\circ\tau(b)} \\
  & = \tgraph{\mu(N)}.
\endalign
$$
 Now apply $\psisharp\pibar1$ to the relation
$$
[\Delta(E/M)] = \pican0(E/M) +\pican 1(E/M) + \pican 2(E/M) .
$$
 Then from the characterizing property 3.1(1) of the $\pican i$ we get
that
$$
 \tgraph{\mu(-1)}\circ \pican i(E/M) = \pican i(E/M) \circ
\tgraph{\mu(-1)} = (-1)^i \pican i(E/M) .
$$
 It follows that $\psi_1^*(\lambdabar\rel)$ annihilates $\pican0(E/M)$
and $\pican2(E/M)$, whence the same is true of $\psisharp\pibar1$, and
also that $\psi_1^*(\lambdabar\rel)\circ \pican1(E/M) =\pican1(E/M)$.
Hence
$$
 \psisharp\pibar 1 =
\psi_1^*(\thetabar\rel) \circ \pican1(E/M).
$$
 Now multiply both sides of this equation by $N$. Then again using 3.1(1)
we get
$$\spreadlines{1\jot}\align
 N(\psisharp\pibar 1) &=
\psi_1^*(\thetabar\rel) \circ \tgraph{\mu(N)} \circ \pican1(E/M) \\
 &= \tgraph{\mu(N)} \circ \pican1(E/M) \\
 &= N \pican1(E/M).
\endalign
$$
 Therefore $\psisharp\pibar1 = \pican1(E/M)$ in $\Chow 1 {\E\times_M\E}$
as required.

 Now consider $i=0,2$.  From 3.2.3 and 3.2.7 we have
$$
 \psisharp\pibar i = \psisharp \pbar_i - \tfrac12
\psisharp((\phibar\times_\Mbar\phibar)^{*} [\fra] ) .
$$
 As we have already observed (3.2.1) that $\psisharp \pbar_i =\pican
i(E/M)$ for $i=0,2$ \cite{K\"unnemann, 1994, 4.1.2(iv)}, what we need to
show is that $\psisharp((\phibar\times_\Mbar\phibar)^{*} [\fra])=0$.  But
from
the definition of $\psisharp$ and the commutativity of the diagram
$$\CD
 E\times_M E @CC \psi_1 C  \Ebar\times_\Mbar\Ebar \\
 @VVV   @VVV  \\
 M @CC j C \Mbar
\endCD
$$
 it follows that
$$
 \psisharp((\phibar\times_\Mbar\phibar)^{*} [\fra]) = (\phi\times_M\phi)^*
j^*
([\fra]) .
$$
 Therefore it will suffice to prove that $j^* ([\fra])=0$ in $\Chow 1 M$,
or equivalently, that $[\fra]\in\Chow 1\Mbar$ can be supported in
{}~$\Minf$.

  To see this, let $\ebar_0 := \ebar := \alphabar((0,0),\Mbar)$, and let
$\ebar_1 := \alphabar((1,0),\Mbar)$, and $\ebar_2 :=
\alphabar((0,1),\Mbar)$.  Then for distinct $i,j\in \{0,1,2\}$ the
intersection
cycle
$[\ebar_i \cdot \ebar_j] = 0$ in $\Chow 2{\Ebar}$, since these sections
are
distinct in every fibre.  Now let $\eta$ denote the generic point
of ~$M$.  Then $N\ebar_0(\eta)$ and $N\ebar_1(\eta)$ are
$\QQ(\eta)$-rational zero cycles on $\Ebar_\eta$, each summing to
$\ebar_0(\eta)$ on $\Ebar_\eta$, whence by Abel's theorem they are
linearly equivalent on ~$\Ebar_\eta$.  More precisely,
$$
 N\ebar_0(\eta) = N\ebar_1(\eta) + \operatorname{div}(f_\eta)
$$
 for some $f_\eta\in\QQ(\eta)$.  But then as cycles
$$
  N\ebar_0 = N\ebar_1  +  \phibar^*(\frak b) + D + \operatorname{div}(F)
\tag 1
$$
 for some zero-cycle ~$\frak b$ on ~$M$ and some divisor ~$D$ supported in
{}~$\Einf$ and some $F\in \QQ(\Ebar)$ (corresponding to $f_\eta$).  If we
now
intersect both sides of ~(1) with $\ebar_2$ and push the resulting cycle
down
to $\Mbar$ by ~$\phibar_*$, then we find that $\frak b$ is linearly
equivalent
on $\Mbar$ to some zero-cycle $\frak b'$ supported on ~$\Minf$.  Therefore
we
may rewrite ~(1) as
$$
 N\ebar_0 \sim_{\text{lin}} N\ebar_1  +  D'
\tag 2
$$
 on $\Ebar$, with $D'$ a divisor supported in ~$\Einf$.
Now intersecting (2) with $\ebar_0 =\ebar$, it follows that the
self-intersection cycle $N[\ebar\cdot\ebar]$ can be supported in
$\Einf$.  And since $[\fra] =\phibar_*([\ebar\cdot\ebar])$ in $\Chow 1
\Mbar$, it follows that $[\fra]$ can be supported in $\Minf$ and
$j^*([\fra]) =0$ in $\Chow 1 M$, which was what we needed to show.
\Qed
\enddemo

\subsubhead Remark \endsubsubhead
 A similar argument can be used to show that $[\e\cdot\e]=0$ in $\Chow 2
E$.

\subhead 3.2.9. Definition of $\pibarinf$  \endsubhead
 Let
$$
 \pibarf := \sum_{i=0}^2 \pibari \qquad \text{\sl in } \Chow
2{\Ebar\times\Ebar} .
$$
 Then Proposition ~3.2.8 implies that $\psisharp\pibarf =
[\Delta(E/M)]$.  Let
$$
 \pibarinf := [\Delta(\Ebar)] -  \pibarf \qquad\text{\sl in } \Chow
2{\Ebar\times\Ebar} .
$$
 Then it follows from the mutual orthogonality and idempotency of the
$\pibari$, for $i=0,1,2$, that $\pibarf$ and $\pibarinf$ are projectors
as well, and $\pibarinf$ is orthogonal to all the others.   In fact we
can say more, using Proposition~3.2.8 and the geometry of ~$\Ebar$.

\proclaim{Lemma 3.2.10 {\rm (the structure of $\pibarinf$)}}
 In $\Chow 2{\Ebar\times\Ebar}$,
$$
 \pibarinf = \sum_{c\in\Minf} \pibar c ,
$$
 where the $\pibar c$ are mutually orthogonal projectors, orthogonal to
the $\pibari$ for $0\le i\le 2$, and of the form
$$
 \pibar c = \sum_{i,j\in\ZmodN} r_c(i,j) [\theta_c(i)]\times_{\{c\}}
[\theta_c(j)]
$$
 for some rational numbers $r_c(i,j)$.
\endproclaim

\demo{Proof}
 Consider the diagram
$$\spreadmatrixlines{1.0\jot}\matrix\format%
\r&\enspace\c&\enspace\c&\enspace\c&\enspace\c&\enspace\c&\enspace\l\\
\CH_2(\Einf\times_{\Minf}\Einf,\QQ) & \lra
&  \CH_2(\Ebar\times_\Mbar \Ebar),\QQ) & \overset{\psi_1^*}\to\lra &
\CH_2(E\times_M E,\QQ)
&\lra & 0  \\
 & \searrow &   @VV (\psi_2)_* V   @VVV   \\
 & & \CH^2(\Ebar\times\Ebar,\QQ) & \lra & \CH^2(E\times E,\QQ) &\lra & 0
\endmatrix
$$
 whose horizontal rows are exact \cite{Fulton, 1984, 1.8, p.21}.
Then in the notation of ~3.2.7 $\pibarf =(\psi_2)_*\pibarrel f$, and it
follows
from Proposition~3.2.8 that the difference $[\Delta(\Ebar)]\rel -
\pibarrel f$
in $\CH_2(\Ebar\times_\Mbar \Ebar),\QQ)$ maps to zero in $\CH_2(E\times_M
E,\QQ)$.  Therefore $\pibarinf \in \Chow 2{\Ebar\times\Ebar}$ must be in
the
image of
$\CH_2(\Einf\times_{\Minf}\Einf,\QQ)$.  Thus, since
$\Einf\times_{\Minf}\Einf$ is 2@-dimensional, with components of the
form $\theta_c(i) \times_{\{c\}} \theta_c(j)$, we get that
$\pibarinf$ can be written in the form indicated.  On the other hand,
the disjointness of the fibres $\Ebar_c =\Einf_c$ implies that the
distinct $\pibar c$, for $c\in\Minf$, are mutually orthogonal, and
idempotent; and hence as constituents of $\pibarinf$, orthogonal also to
$\pibari$ for $0\le i\le2$.  This proves the lemma.
\Qed
\enddemo

\subhead 3.3. Extending canonical relative projectors to $\Etwotild$
\endsubhead
\nobreak \par\nobreak
\subsubhead 3.3.1. Introduction to the method \endsubsubhead
 The basic idea of the projectors we now define for
$\Etwotild$ is that we would like them to be the tensor products over
$\Mbar$,
in the sense of 1.1.3, of the projectors $\pibar i$
defined above for $\Ebar$.  But since neither $\Ebar$ nor $\Etwotild$ is
smooth over $\Mbar$, and further since after the blowing-up
$\beta:\Etwotild\to\Etwobar$ $\Etwotild$ is no longer a product over
{}~$\Mbar$,
the definition 1.1.3 of the tensor product of
correspondences does not directly apply to our situation. For
this reason we shall define projectors for
$\Etwotild$ directly as combinations of graphs and transposes of graphs
of morphism, as we did for $\pibar i$.  In particular this means that
again we start with explicit representative cycles.

   Firstly we write down correspondences in $\Etwotild\times\Etwotild$
but supported on $\Etwotild\times_\Mbar\Etwotild$ that act on $\Etwotild$
like $\pibar i$ on one fibre factor and identity on the other (in spite
of the fact that the construction of $\Etwotild$ by desingularizing
$\Etwobar$ destroyed the fibre product structure!), for $0\le i\le 2$.
Then in order to get the mutually orthogonal projectors we actually want
from these, we have to show that the correspondence that acts as
$\pibar{i_1}$
on the first factor and identity on the second commutes with the
correspondence that acts as $\pibar{i_2}$ on the second factor and
identity on the first.  Once that is done, we check that the restrictions
of these projectors to $\Etwo\times_M\Etwo$, in a similar sense as 3.2.7
and ~3.2.8, see 3.3.7 to ~3.3.9 below, are indeed tensor products of the
canonical relative projectors.  Finally, similarly as for $\Ebar$ we
define a $\pitildinf(\Etwotild/\Mbar)$ for ~$\Etwotild$.

\subhead 3.3.2. Definition of $\pitild0\sparen j(\Etwotild/\Mbar)$ and
$\pitild2\sparen j(\Etwotild/\Mbar)$ \endsubhead
  Recall from 3.2.2(4) that we wrote $\pbar_2 = \graph{\mubar(0)}$ and
$\pbar_0 = \tgraph{\mubar(0)}$, where $\mubar(0) :=
\alphabar((0,0),\bigdot)\circ \phibar :\Ebar \to \Ebar$ is the projection
onto the zero-section morphism.  Consider now $\mubar(0)\times_\Mbar
\id_{\Ebar} : \Etwobar\to\Etwobar$.  Since the image of this map is
disjoint from the center $\Etwobarinf_0$ of the blowing up ~$\beta :
\Etwotild\to\Etwobar$,  by factoring it through ~$\beta$ it lifts to a
morphism of $\Etwotild$ that respects the fibre structure of
$\Etwotild\to\Mbar$.  More precisely, let
$$
 \mutild(0,1) := (\beta')\inv \circ (\mubar(0)\times_\Mbar
\id_{\Ebar}) \circ \beta :\Etwotild \lra \Etwotild
$$
 where $\beta'$ is the restriction of $\beta$ to $\Etwotild -
\beta\inv(\Etwobarinf_0)$, where it is an isomorphism.  If we define
$\mutild(1,0)$ similarly, then the product in either order
$$
 \mutild(0,0) := \mutild(0,1)\circ \mutild(1,0) =
\mutild(1,0) \circ \mutild(0,1) = \alphatild(\bold
0,\phitwotild(\bigdot)):\Etwotild \lra \Etwotild
$$
 is the projection onto the zero-section of ~$\Etwotild$, where
$\alphatild: (\ZmodN)^4\times\Mbar\to\Etwotild$ is the level-$N$
structure.

  Now let
$$\spreadlines{1\jot} \align
 \pitild0\sparen1(\Etwotild/\Mbar) &:= \tgraph{\mutild(0,1)} -
\tfrac 1 2 \tgraph{\mutild(0,1)}\circ\graph{\mutild(0,1)} \\
 \pitild0\sparen2(\Etwotild/\Mbar) &:= \tgraph{\mutild(1,0)} -
\tfrac 1 2  \tgraph{\mutild(1,0)}\circ\graph{\mutild(1,0)} \\
  \pitild2\sparen1(\Etwotild/\Mbar) &:= \graph{\mutild(0,1)} -
\tfrac 1 2 \tgraph{\mutild(0,1)}\circ\graph{\mutild(0,1)} \\
 \pitild2\sparen2(\Etwotild/\Mbar) &:= \graph{\mutild(1,0)} -
 \tfrac 1 2 \tgraph{\mutild(1,0)}\circ\graph{\mutild(1,0)} ,
\endalign
$$
 where the notation is chosen to suggest that $\pitild i\sparen
j(\Etwotild/\Mbar)$ acts like $\pibar i$ on the $j^{\text{th}}$ fibre
factor and identity on the other.  Then the idempotency of each of these,
and the orthogonality of $\pitild0\sparen j(\Etwotild/\Mbar)$ and
$\pitild2\sparen j(\Etwotild/\Mbar)$ for fixed ~$j$, follows easily from
observing that
$$\spreadlines{1\jot} \align
  \graph{\mutild(0,1)} \circ \tgraph{\mutild(0,1)} &=
(\mutild(0,1)\times \mutild(0,1))_*([\Delta(\Etwotild)]) = 0 \\
 \graph{\mutild(1,0)} \circ \tgraph{\mutild(1,0)} &=
(\mutild(1,0)\times \mutild(1,0))_*([\Delta(\Etwotild)])  = 0 ,
\endalign
$$
 in $\Chow 3{\Etwotild\times\Etwotild}$.

  Similarly as in 3.2.3 now let $\frak a\sparen1$ and $\frak a\sparen2$
be two disjoint zero cycles on $\Mbar$ that both also represent $[\frak
a] = \phibar_*([\ebar\cdot\ebar])$ in $\Chow 1{\Mbar}$.  Then the
correction term $\tfrac12 \tgraph{\mutild(0,1)}\circ\graph{\mutild(0,1)}$
can be represented by a 3-dimensional cycle $\frak b\sparen1$ supported
on $(\phitwotild\times_\Mbar\phitwotild)^{-1}(|\fra\sparen1|)$, where
$|\frak c|$ denotes the support of a zero cycle $\frak c$ on ~$\Mbar$,
and similarly there is a cycle $\frak b\sparen2$ representing $\tfrac 1 2
\tgraph{\mutild(1,0)}\circ\graph{\mutild(1,0)}$ and supported on
$(\phitwotild\times_\Mbar\phitwotild)^{-1}(|\fra\sparen2|)$.  This can be
seen by a direct computation:  In order to have a proper intersection on
$\Etwotild\times\Etwotild\times\Etwotild$ we can move (similarly as in the
proof of~3.2.2) inside the second factor by looking first at the
$\Etwobar$
over which it lies and there moving the zero section ~$\ebar_0$ in the
relevant
factor ~$\Ebar$ to a cycle $\ebar_0'$ such that $\ebar_0$ and $\ebar_0'$
intersect properly and have no common points over the cusps and such that
$\ebar_0'$ does not pass through the crossing points of the components
over the
cusps.  This then gives a corresponding moving for the
$\tgraph{\mutild(1,0)}$
which leads to a proper intersection.  Then we get (at least set
theoretically)
$$\spreadlines{1\jot} \align
 \frak b\sparen 1 = \{(\tilde x, \tilde y) :\ &\tilde x , \tilde y
\in\Etwotild,\ \phitwotild(\tilde x) = \phitwotild(\tilde y) \in \fra
\sparen 1,\\
 &\ \beta(\tilde x) = (x_1,x_2), \ \beta(\tilde y) =(y_1,y_2), \
x_2=y_2\}
\endalign
$$

 For later reference note also that the correction terms $[\frak
b\sparen1]$ and $[\frak b\sparen2]$ are nilpotent of order ~2 in $\Chow
3{\Etwotild\times\Etwotild}$, and orthogonal to each other.

\subhead 3.3.3. Definition of $\pitild1\sparen j(\Etwotild/\Mbar)$
\endsubhead
 Now to define a correspondence that acts as $\pibar 1$ on one fibre
factor and identity on the other, we first observe that, for
$g\in\G$ acting on $\Ebar$, the fibre product morphism $g\times_\Mbar
\id_{\Ebar} :\Etwobar \to \Etwobar$ scheme-theoretically preserves the
center $\Etwobarinf_0$ of the blowing-up $\beta:\Etwotild\to\Etwobar$.
Therefore it lifts uniquely to a morphism, say $\chitild(g,\id) :
\Etwotild
\to \Etwotild$, of ~$\Etwotild$ \cite{Hartshorne, 1977, II.7.15, p.165}.
Similarly $\id_{\Ebar}\times_\Mbar g$ lifts to a morphism, say
$\chitild(\id,g)$, and moreover, for $g_1,g_2\in\G$ we have
$$
\chitild(g_1,\id) \circ \chitild(\id,g_2) = \chitild(\id,g_2) \circ
\chitild(g_1,\id) =: \chitild(g_1,g_2) .
$$
  Thus $\G^2:=\G\times\G$ acts as a group of fibrewise automorphisms on
{}~$\Etwotild$, and this action extends $\QQ$-linearly to give a
homomorphism
$$
 \QQ[\G^2] \lra \Chow 3 {\Etwotild\times\Etwotild} .
$$
  As special cases, for $\boldkey a = (a_1,a_2) \in (\mu_2\times\mu_2)$
we write $\mutild(\boldkey a) :\Etwotild \to \Etwotild$ for the
corresponding morphism, and for $\boldb\in \ZmodNsq\times\ZmodNsq$ we let
$\tautild(\boldb):\Etwotild \to\Etwotild$ denote the corresponding
morphism.  Then analogously as in the definition ~3.2.5 of $\pibar1$, let
$$\spreadlines{1\jot} \alignat2
 \lambdatild\sparen1 &:= \frac 1 2 \big( \graph{\mutild(1,1)} -
\graph{\mutild(-1,1)} \big) \qquad &
 \lambdatild\sparen2 &:= \frac 1 2 \big( \graph{\mutild(1,1)} -
\graph{\mutild(1,-1)} \big) \\
  \thetatild\sparen1 &:= \frac 1 {N^2} \sum_{b\in(\ZZ/N\ZZ)^2}
\graph{\tautild(b,0)}  &
  \thetatild\sparen2 &:= \frac 1 {N^2} \sum_{b\in(\ZZ/N\ZZ)^2}
\graph{\tautild(0,b)} ,
\endalignat
$$
 and then
$$\spreadlines{1\jot} \align
  \pitild1\sparen1(\Etwotild/\Mbar) &:=
\thetatild\sparen1\circ\lambdatild\sparen1 =
\lambdatild\sparen1\circ\thetatild\sparen1 \\
  \pitild1\sparen2(\Etwotild/\Mbar) &:=
\thetatild\sparen2\circ\lambdatild\sparen2 =
\lambdatild\sparen2\circ\thetatild\sparen2 .
\endalign
$$
 As in the definition of $\pibar1$, it follows easily from identities in
the group ring $\QQ[\G\times\G]$ that $\lambdatild\sparen j$ commutes
with $\thetatild \sparen j$, and that all the $\lambdatild\sparen j$ and
$\thetatild\sparen j$ and thus the $\pitild1\sparen j(\Etwotild/\Mbar)$
are idempotent.  Here we also have
$$
 \pitild1\sparen 1(\Etwotild/\Mbar) \circ \pitild1\sparen
2(\Etwotild/\Mbar) = \pitild1\sparen 2(\Etwotild/\Mbar) \circ
\pitild1\sparen 1(\Etwotild/\Mbar) ,
$$
 because the two factors of $\G\times\G$ commute.

 The following lemma should be compared with Proposition~3.2.6.

\proclaim{Lemma 3.3.4}
 For fixed $j=1$ or~2, the $\pitild i \sparen j(\Etwotild/\Mbar)$, for
$i=0,1,2$, are mutually orthogonal idempotents in
$\CH^3(\Etwotild\times\Etwotild,\QQ)$.
\endproclaim

\demo{Proof}
 All that remains to be checked is that $\pitild 1 \sparen
j(\Etwotild/\Mbar)$ is orthogonal to both $\pitild 0 \sparen
j(\Etwotild/\Mbar)$ and $\pitild 2 \sparen j(\Etwotild/\Mbar)$.  But for
this one can argue similarly as for Proposition~3.2.6, that
$\lambdatild\sparen 1$ is orthogonal to both $\graph{\mutild(0,1)}$ and
$\tgraph{\mutild(0,1)}$, and $\lambdatild\sparen 2$ is orthogonal to both
$\graph{\mutild(1,0)}$ and $\tgraph{\mutild(1,0)}$.
\Qed
\enddemo

\subhead 3.3.5. Definition of $\pitild{i_1,i_2}(\Etwotild/\Mbar)$
\endsubhead
 For $0\le i_1,i_2\le2$ define
$$
 \pitildi(\Etwotild/\Mbar) := \pitild{i_1}\sparen1(\Etwotild/\Mbar)
\circ \pitild{i_2}\sparen2(\Etwotild/\Mbar) .
$$

\proclaim{Proposition 3.3.6}
 The $\pitildi(\Etwotild/\Mbar)$, for $0\le i_1,i_2\le2$, are mutually
orthogonal projectors in ~$\Chow 3 {\Etwotild\times\Etwotild}$.
\endproclaim

\demo{Proof}
 This proposition will follow immediately from the Lemma~3.3.4 as soon as
we verify the commutativity relation, that for all $i_1,i_2=0,1,2$,
$$
\pitild{i_1}\sparen1(\Etwotild/\Mbar) \circ \pitild{i_2}
\sparen2(\Etwotild/\Mbar) = \pitild {i_2}\sparen2(\Etwotild/\Mbar)
\circ \pitild {i_1}\sparen1(\Etwotild/\Mbar) .
\tag 1
$$
  We shall verify this case by case.

 {\it Case $i_1 =i_2 =1$.}
 We have already seen in 3.3.3 that (1) holds because the
$\pitild{1}\sparen j(\Etwotild/\Mbar)$, for $j=1,2$, are homomorphic
images of commuting projectors in the group ring $\QQ[\G^2]$.

  {\it Case $i_1=1\ne i_2$ or $i_1\ne 1 =i_2$.}  In this case the
commutativity relation~(1) will follow if we can show that the graph of
$\chitild(g_1,\id)$ commutes with both the graph and the transposed graph
of $\mutild(0,1)$, and similarly that the graph of $\chitild(\id,g_2)$
commutes with both the graph and the transposed graph of $\mutild(1,0)$.
But recalling that $\tgraph{\chitild} =
\graph{\chitild\inv}$ whenever $\chitild$ is an automorphism of
$\Etwotild$, and then using ~1.1.7, the problem reduces to proving that
for any $g\in\G$,
$$\spreadlines{1\jot} \align
\mutild(0,1)\circ\chitild(\id,g) &= \chitild(\id,g)\circ\mutild(0,1), \\
\mutild(1,0)\circ\chitild(g,\id) &= \chitild(g,\id)\circ\mutild(1,0).
\endalign
$$
 as endomorphisms of ~$\Etwotild$.

 To prove the first of these, say, since the argument is the same for
both, first recall that by definition $\mutild(0,1) := (\beta')\inv \circ
(\mubar(0,1) \circ \beta$,  where $\beta'$ is the restriction of $\beta$
to $\Etwotild - \beta\inv(\Etwobarinf_0)$, on which it is an isomorphism,
and $\mubar(0,1) := \mubar(0) \times_\Mbar \id_{\Ebar} :
\Etwobar\to\Etwobar$.  On the other hand, the automorphism
$\chitild(\id,g)$ preserves the exceptional divisor of $\Etwotild$, as
it was lifted to a morphism on $\Etwotild$ from $\chibar(\id,g) :=
\id\times_\Mbar g: \Etwobar\to\Etwobar$, which preserves the center
$(\Etwobarinf)_0$ of the blowing-up.  Therefore by \cite{Hartshorne, 1977,
II.7.15, p.165}
$$
 \beta \circ\chitild(\id,g) = \chibar(\id,g) \circ \beta .
$$
 Combining this with the definition of $\mutild(0,1)$, we get
$$\spreadlines{1\jot} \align
\mutild(0,1) \circ \chitild(\id,g) &= (\beta')\inv \circ \mubar(0,1)
  \circ \beta \circ \chitild(\id,g)  \\
 &= (\beta')\inv \circ \mubar(0,1) \circ \chibar(\id,g) \circ \beta \\
 &= (\beta')\inv \circ \chibar(\id,g) \circ \mubar(0,1) \circ \beta \\
 &= \chitild(\id,g) \circ \mutild(0,1).
\endalign
$$

  {\it Case $i_1 = i_2\ne 1$.}
 First consider the cases $\pitild{0,0}(\Etwotild/\Mbar)$ and
$\pitild{2,2}(\Etwotild/\Mbar)$ which are similar.  Take for instance
$\pitild{2,2}(\Etwotild/\Mbar)$.  Using, as remarked in ~3.3.2, that the
correction terms are orthogonal, we have
$$\spreadlines{1\jot}\aligned
 \pitild{2,2}(\Etwotild/\Mbar) :&= \pitild2\sparen1(\Etwotild/\Mbar)
\circ \pitild2\sparen2(\Etwotild/\Mbar) \\
 & = \graph{\mutild(0,0)} -\tfrac12 \graph{\mutild(0,1)} \circ
\tgraph{\mutild(1,0)}  \circ \graph{\mutild(1,0)} \\
 & \hphantom{= \graph{\mutild(0,0)}}\ - \tfrac12 \tgraph{\mutild(0,1)}
\circ \graph{\mutild(0,1)} \circ \graph{\mutild(1,0)} .
\endaligned
\tag 2
$$
 Thus proving the commutativity relation ~(1) for
$\pitild{0,0}(\Etwotild/\Mbar)$ and $\pitild{2,2}(\Etwotild/\Mbar)$
reduces to proving the relations
$$\spreadlines{1\jot}  \aligned
  \tgraph{\mutild(0,1)} \circ \graph{\mutild(1,0)} &=
 \graph{\mutild(1,0)} \circ \tgraph{\mutild(0,1)} ,\\
 \tgraph{\mutild(1,0)} \circ \graph{\mutild(0,1)} &=
   \graph{\mutild(0,1)} \circ \tgraph{\mutild(1,0)} ,
\endaligned
\tag 3
$$
 which are straightforward to verify by direct computation.

  {\it Case $1\ne i_1 \ne i_2\ne 1$.}
 It remains to consider $\pitild{0,2}(\Etwotild/\Mbar)$ and
$\pitild{2,0}(\Etwotild/\Mbar)$.  Take for instance
$\pitild{0,2}(\Etwotild/\Mbar)$.  Again using the orthogonality of the
correction terms, now we get
$$\spreadlines{1\jot}\aligned
 \pitild{0,2}(\Etwotild/\Mbar) :&= \pitild0\sparen1(\Etwotild/\Mbar)
\circ \pitild2\sparen2(\Etwotild/\Mbar) \\
 & = \tgraph{\mutild(0,1)} \circ \graph{\mutild(1,0)} -\tfrac12
\tgraph{\mutild(0,1)} \circ \tgraph{\mutild(1,0)}  \circ
\graph{\mutild(1,0)} \\
 & \hphantom{= \graph{\mutild(0,0)}}\ - \tfrac12 \tgraph{\mutild(0,1)}
\circ \graph{\mutild(0,1)} \circ \graph{\mutild(1,0)} .
\endaligned
\tag 4
$$
 Then after writing out $\pitild2\sparen2(\Etwotild/\Mbar) \circ
\pitild0\sparen1(\Etwotild/\Mbar)$ we find that the commutativity
relation~(1) also follows in this case from the relations~(3).
 This completes the proof of the proposition.
\Qed
\enddemo

\subhead 3.3.7.  Notations and definitions related to cycles on
$\Etwotild\times_\Mbar \Etwotild$  \endsubhead
  As in 3.2.7, when $\alpha :\Etwotild \to \Etwotild$ is a morphism that
respects the fibre structure of $\Etwotild \to \Mbar$ then the graph
$\Gamma_\alpha$ of $\alpha$ is supported on
$\Etwotild\times_\Mbar\Etwotild$ and we write $\graph\alpha\rel$ and
$\tgraph\alpha\rel$ for its class and the class of its transpose in
$\CH_3(\Etwotild\times\Etwotild,\QQ)$.  Now consider the inclusions
$$
 \lup 2\psi : \Etwo\times_M\Etwo @C \lup2\psi_1 CC
\Etwotild\times_\Mbar\Etwotild @C \lup2\psi_2 CC \Etwotild
\times\Etwotild .
$$
 Then we have $\graph\alpha = (\lup2\psi_2)_* \graph\alpha\rel$, and,
by abuse of notation we define
$$
 \psitwosharp (\graph\alpha) := \lup2\psi_1^* (\graph\alpha\rel) \qquad
\text{in } \Chow 2{\Etwo\times_M\Etwo},
$$
 and similarly for the transpose of the graph.  If $\beta:\Etwotild \to
\Etwotild$ is another morphism respecting the fibre structure of
$\Etwotild \to\Mbar$, then
$$\spreadlines{1\jot} \aligned
 \psitwosharp (\graph{\alpha\circ\beta}) &= \psitwosharp (\graph\alpha)
\circ \psitwosharp(\graph\beta) \\
 \psitwosharp (\tgraph{\alpha\circ\beta}) &= \psitwosharp (\tgraph\beta)
\circ \psitwosharp(\tgraph\alpha) .
\endaligned
\tag 1
$$
 Also as before we extend these definitions by linearity.

  Next we apply these definitions to the explicit cycles in ~3.3.3.
There we defined $\lambdatild\sparen j$ and $\thetatild\sparen j$, for
$j=1,2$, as linear combinations of graphs of automorphisms that respect
the fibre structure of $\Etwotild\to\Mbar$, so $\lambdatild\sparen
j\strut\rel$ and $\thetatild\sparen j\strut\rel$ in
$\CH_3(\Etwotild\times_\Mbar\Etwotild,\QQ)$ and
$\psitwosharp(\lambdatild\sparen j)$ and $\psitwosharp(\thetatild\sparen
j)$ in $\Chow 2 {\Etwo\times_M\Etwo}$ are defined, for $j=1,2$.  If we
write, as we may,
$$\spreadlines{1\jot} \aligned
 \pitild1\sparen1(\Etwotild/\Mbar) &= \frac 1{2N^2}\sum_{g\in\G}
\varepsilon(g)^{-1} \graph{\chitild(g,\id)} \\
 \pitild1\sparen2(\Etwotild/\Mbar) &= \frac 1{2N^2}\sum_{g\in\G}
\varepsilon(g)^{-1} \graph{\chitild(\id,g)}
\endaligned
\tag 2
$$
 with $\varepsilon$ as in 3.2.5, then $\pitild1\sparen
j\strut\rel(\Etwotild/\Mbar)$ and $\psitwosharp(\pitild1\sparen
j(\Etwotild/\Mbar))$ are defined in the obvious way, for $j=1,2$, and
moreover in $\Chow 2 {\Etwo\times_M\Etwo}$
$$
 \psitwosharp(\pitild1\sparen j(\Etwotild/\Mbar)) =
\psitwosharp(\lambdatild\sparen j) \circ \psitwosharp(\thetatild\sparen
j) = \psitwosharp(\thetatild\sparen j) \circ
\psitwosharp(\lambdatild\sparen j) .
$$

  Next we want to apply the definitions above to $\pitild i\sparen j
(\Etwotild/\Mbar)$, for $i=0,2$ and $j=1,2$, as defined in ~3.3.2.
Recall that there we chose explicit cycles $\frak b\sparen j$ supported
on $\Etwotild\times_\Mbar\Etwotild$ and such that
$$
 [\frak b\sparen1] = \tfrac12 \tgraph{\mutild(0,1)}\circ
\graph{\mutild(0,1)} \quad \text{ and }\quad [\frak b\sparen2] =
\tfrac12 \tgraph{\mutild(1,0)}\circ
\graph{\mutild(1,0)} .
$$
 Thus we have elements $[\frak b\sparen j]\rel \in
\CH_3(\Etwotild\times_\Mbar\Etwotild,\QQ)$, and therefore also elements
$\pitild i \sparen j\strut\rel(\Etwotild/\Mbar) \in
\CH_3(\Etwotild\times_\Mbar\Etwotild,\QQ)$ such that
$$
 \pitild i \sparen j(\Etwotild/\Mbar) = (\lup2\psi_2)_*(\pitild i \sparen
j\strut\rel(\Etwotild/\Mbar))
$$
  for $i=0,2$ and $j=1,2$.  Hence we may also define
$$
 \psitwosharp\pitild i\sparen j (\Etwotild/\Mbar) := \lup2\psi_1^*
(\pitild i\sparen j\strut\rel (\Etwotild/\Mbar)) \in
\Chow2{\Etwo\times_M\Etwo} ,
$$
 for $i=0,2$ and $j=1,2$.

\medpagebreak

 Finally, we would like to apply the definitions at the beginning of
this section to $\pitildi(\Etwotild/\Mbar)$, for $0\le i_1,i_2\le 2$.
The following lemma shows how we can do this, even though these
projectors were defined in ~3.3.5 as a composition of cycle classes, i.e.,
$$
 \pitildi(\Etwotild/\Mbar) := \pitild{i_1}\sparen1(\Etwotild/\Mbar)
\circ \pitild{i_2}\sparen2(\Etwotild/\Mbar) .
$$

\proclaim{Lemma 3.3.8}
\nobreak\par\nobreak
\noindent {\rm (1)}\enspace
 There exist $\pitildi\rel(\Etwotild/\Mbar)
\in \CH_3(\Etwotild\times_\Mbar\Etwotild,\QQ)$ such that
$$
 (\lup2\psi_2)_* \pitildi\rel(\Etwotild/\Mbar) =
\pitildi(\Etwotild/\Mbar) ,
$$
 for $0\le i_1,i_2\le 2$.
\par
\noindent {\rm (2)}\enspace  Let
$$
 \psitwosharp \pitildi(\Etwotild/\Mbar) :=
\lup2\psi_1^*(\pitildi\rel(\Etwotild/\Mbar)) .
$$
 Then in $\Chow 2{\Etwo\times_M\Etwo}$ we have
$$
 \psitwosharp \pitildi(\Etwotild/\Mbar) = \psitwosharp
(\pitild{i_1}\sparen1(\Etwotild/\Mbar) \circ \psitwosharp
\pitild{i_2}\sparen2(\Etwotild/\Mbar) ,
$$
 for $0\le i_1,i_2\le 2$.
\endproclaim

\demo{Proof}
 We will prove this lemma case by case, as we did for
Proposition~3.3.6.

 {\it Case $i_1=i_2=1$.}
 Consider the character $\varepsilon_2:\G^2 \to
\{\pm1\}$ defined by $\varepsilon_2(g_1,g_2) :=
\varepsilon(g_1)\varepsilon(g_2)$, where $\varepsilon :\G \to \{\pm1\}$
is the character defined in ~3.2.5.  Then
$$\spreadlines{1\jot}\aligned
 \pitild{1,1}(\Etwotild/\Mbar) :&= \pitild 1\sparen1(\Etwotild/\Mbar)
\circ \pitild 1\sparen2(\Etwotild/\Mbar) \\
 & = \frac 1{4N^4} \sum_{(g_1,g_2)\in\G^2} \varepsilon_2(g_1,g_2)^{-1}
\graph{\chitild(g_1,g_2)} .
\endaligned
\tag 3
$$
 We may use this expression to define
$\pitild{1,1}\rel(\Etwotild/\Mbar)$, proving ~\therosteritem1, and then
\therosteritem2 follows from observing that
$$
 \psitwosharp(\graph{\chitild(g_1,g_2)} =
\psitwosharp(\graph{\chitild(g_1,\id)} \circ
\psitwosharp(\graph{\chitild(\id,g_2)} ,
$$
 see ~3.3.7(1).

 {\it Case $i_1=1\ne i_2$ or $i_1\ne 1 =i_2$.}
  Consider for instance $i_1=1$ and $i_2=2$.  Then
$$\spreadlines{1\jot}\aligned
 \pitild{1,2}(\Etwotild/\Mbar) :&= \pitild 1\sparen1(\Etwotild/\Mbar)
\circ \pitild 2\sparen2(\Etwotild/\Mbar) \\
 &= \frac1{2N^2} \sum_{g\in\G} \big( \graph{\chitild(g,\id)\circ
\mutild(1,0)} -  [(\id_{\Etwotild} \times_\Mbar \chitild(g,\id))_* (\frak
b\sparen2)] \big) ,
\endaligned
\tag 4
$$
 where the second term in each summand comes from ~1.1.7 applied to
$\graph{\chitild(g,\id)} \circ [\frak b\sparen2]$.    Now
$\graph{\chitild(g,\id)\circ \mutild(1,0)}\rel \in
\CH_3(\Etwotild\times_\Mbar\Etwotild,\QQ)$ is defined, because it comes
from the graph of a morphism, and the second term is supported on
$\Etwotild\times_\Mbar\Etwotild$ (indeed even on $(\phitwotild
\times_\Mbar\phitwotild)^{-1}(|\fra\sparen2|)$), as well.  Therefore we
may define $\pitild{1,2}\rel(\Etwotild/\Mbar)$ by the explicit
expression ~(4), and this proves part~\therosteritem1 in this case.

 As for showing that, with the definitions as given here,
$$
 \psitwosharp \pitild{1,2}(\Etwotild/\Mbar) = \psitwosharp
(\pitild{1}\sparen1(\Etwotild/\Mbar) \circ \psitwosharp
\pitild{2}\sparen2(\Etwotild/\Mbar) ,
\tag 5
$$
 first we claim that $\psitwosharp([\frak b\sparen2]) =0$.  From the
explicit computation of $[\frak b\sparen2]\rel$ in~3.3.2 we get
$$ \spreadlines{1\jot} \align
 \psitwosharp([\frak b\sparen2]) &= \lup2\psi_1^*([\frak b\sparen2]\rel)
\\
 &= \tfrac12 [((\phi\times_M\phi)^*j^*(\fra\sparen2))\times_M \Delta(E/M)]
, \endalign
$$
 and we have already seen in the proof of 3.2.8 that
$j^*([\fra\sparen2])=0$.  On the other hand,
$$
 \psitwosharp([(\id_{\Etwotild} \times_\Mbar \chitild(g,\id))_* (\frak
b\sparen2)]\rel) = (\id_{\Etwotild} \times_\Mbar \chitild(g,\id))_*
\big( \lup2\psi_1^*([\frak b\sparen2]\rel) =0 ,
$$
 where the first equality follows because $\lup2\psi_1$ is an open
immersion which is preserved by the action of $(g,\id)\in\G^2$.  Now (5)
follows for $(i_1,i_2)=(1,2)$, and the other cases are similar, except
that when $i_1=0$ or $i_2=0$ we use transposed graphs throughout.

 {\it Case $i_1 = i_2\ne 1$.}
  Take for instance $(i_1,i_2)=(2,2)$, the other case $(i_1,i_2)=(0,0)$
will be similar.  Then
$$\spreadlines{1\jot} \align
\pitild{2,2}(\Etwotild/\Mbar) &= \pitild2\sparen1(\Etwotild/\Mbar) \circ
\pitild2\sparen2(\Etwotild/\Mbar) \\
  &= \graph{\mutild(0,0)} -
\graph{\mutild(0,1)}\circ [\frak b\sparen2] -  [\frak b\sparen1]
\circ \graph{\mutild(1,0)} \\
  &= \graph{\mutild(0,0)} -(\id_{\Etwotild}\times_\Mbar
\mutild(0,1))_*([\frak b\sparen2]) - (\mutild(1,0)\times_\Mbar
\id_{\Etwotild})^*([\frak b\sparen1]) .
\tag 6 \\
\endalign
$$
 This last expression gives us explicit cycles with which to define
$\pitild{2,2}\rel(\Etwotild/\Mbar)$, proving part~\therosteritem1 for
this case.  To prove part~\therosteritem2 we must verify by
straightforward computation that
$\lup2\psi_1^*((\id_{\Etwotild}\times_\Mbar
\mutild(0,1))_*([\frak b\sparen2]\rel) = 0$, and similarly {\sl mutatis
mutandis;} the proofs are similar to the previous ones.

 {\it Case $1\ne i_1\ne i_2 \ne1$.}
  Take for instance $(i_1,i_2) = (0,2)$.  Then similarly as in the
previous case we have
$$ \spreadlines{1\jot} \align
\pitild{1,2}(\Etwotild/\Mbar) &= \pitild1\sparen1(\Etwotild/\Mbar) \circ
\pitild2\sparen2(\Etwotild/\Mbar) \\
  &= \tgraph{\mutild(0,1)} \circ \graph{\mutild(1,0)} -
(\id_{\Etwotild}\times_\Mbar \mutild(0,1))^*([\frak b\sparen2])
\tag 7 \\
 & \hphantom{\tgraph{\mutild(0,1)} \circ \graph{\mutild(1,0)}}\ -
(\mutild(1,0)\times_\Mbar \id_{\Etwotild})^*([\frak b\sparen1])
\endalign
$$
 Now to see that $\tgraph{\mutild(0,1)} \circ \graph{\mutild(1,0)}$ is
or can be supported on $\Etwotild\times_\Mbar\Etwotild$ we can compute
at the level of cycles where
$$
 \tgraph{\mutild(0,1)} \circ \graph{\mutild(1,0)} =
[\pr_{13*}([(\Gamma_{\mutild(1,0)}\times \Etwotild)\cdot
(\Etwotild\times \tr\Gamma_{\mutild(0,1)})]) .
$$
 Then we see that this can be represented by a cycle supported on the set
$$
 \{(\tilde x,\tilde y) : \phitwotild(\tilde x) =\phitwotild(\tilde y),\
\beta(\tilde x) =(0,x_2),\ \beta(\tilde y)=(y_1,0),\ \text{with }
x_2,y_1\in\Ebar \}
\tag 8
$$
 contained in $\Etwotild\times_\Mbar\Etwotild$.  Using this we can
define $\pitild{0,2}\rel(\Etwotild/\Mbar)$ via formula~(7).  For part
{}~\therosteritem2 we use firstly that the correction terms vanish after
applying $\lup2\psi_1^*$, as above, and that if we use a cycle
representative for $\tgraph{\mutild(0,1)} \circ \graph{\mutild(1,0)}$
supported on the set~(8) then with the obvious notation we get
$$
 \psitwosharp(\tgraph{\mutild(0,1)} \circ \graph{\mutild(1,0)} =
\lup2\psi_1^*(\tgraph{\mutild(0,1)}\rel) \circ
\lup2\psi_1^*(\graph{\mutild(1,0)}\rel) .
$$
 Set-theoretically this is immediate, and in order to see that the
intersection multiplicities are correct use \cite{Weil, 1948, VIII.4,
Thm.10,
p.233}.  This completes the proof of the lemma.
\Qed
\enddemo

\proclaim{Proposition 3.3.9}
 In $\CH^2(\Etwo\times_M\Etwo,\QQ)$ we have
$$
\lup2\psisharp \pitildi(\Etwotild/\Mbar) = \pican{i_1}(E/M)\tensor_M
\pican{i_2}(E/M) ,
\tag 1
$$
 for $0\le i_1,i_2\le 2$, and moreover
$$
\lup 2\psisharp\big(\sum_{i_1+i_2 =i} \pitildi(\Etwotild/\Mbar)\;\big) =
\pi_i^{\text{can}}(\Etwo/M),
\tag 2
$$
 for $0\le i\le 2$.
\endproclaim

\demo{Proof}
 Firstly we claim that
$$
 \lup2\psisharp(\pitild{i_1}\sparen1(\Etwotild/\Mbar)) =
\psisharp(\pibar{i_1}) \tensor_M [\Delta(E/M)] ,
$$
 with the tensor product defined as in  ~1.1.3, and similarly
for $\lup2\psisharp(\pitild{i_2}\sparen2(\Etwotild/\Mbar))$.  For $i_1
=1$ this comes immediately from the expression ~3.3.7(2).  If $i_1
=0$, say, then we have seen in the proof of Lemma~3.3.8 that
$\lup2\psi_1^*([\frak b\sparen j]\rel) =0$, from which it follows that
$\pitild{0}\sparen1\strut\rel(\Etwotild/\Mbar) =
\tgraph{\mutild(0,1)}\rel$.  Therefore,
$\lup2\psi_1^*(\pitild{0}\sparen1\strut\rel(\Etwotild/\Mbar)) =
\tgraph{\mu(0} \tensor_M [\Delta(E/M)]$ as claimed.  The argument is the
same
if $i_1=2$ or if $i_1$ is replace by ~$i_2$.  Hence (1) now follows
from Lemma~3.3.8(2) and Proposition~3.2.8.  Then (2) follows from ~(1),
the K\"unneth formula for relative Chow motives over ~$M$, and the
characterizing property ~3.1(1) of the canonical relative projectors for
abelian schemes.
\Qed
\enddemo

\subhead 3.3.10. Definition of $\pitildinf(\Etwotild/\Mbar)$
\endsubhead
  Let
$$
 \pitildf(\Etwotild/\Mbar) := \sum_{i_1, i_2 =0}^2
\pitild{i_1,i_2}(\Etwotild/\Mbar) .
$$
 Then by Proposition~3.3.8
$\lup2\psisharp\pitildf(\Etwotild/\Mbar) = [\Delta(\Etwo/M)]$ in $\Chow
2{\Etwo\times_M\Etwo}$.  Let
$$
 \pitildinf(\Etwotild/\Mbar) := [\Delta(\Etwotild)] -
\pitildf(\Etwotild/\Mbar) \qquad \text{\sl in }
\CH^3(\Etwotild\times\Etwotild,\QQ) .
$$
  Then it is immediate from the orthogonality and idempotency of the
$\pitildi(\Etwotild/\Mbar)$ that $\pitildf(\Etwotild/\Mbar)$ and
$\pitildinf(\Etwotild/\Mbar)$ are mutually orthogonal projectors, and
that $\pitildinf(\Etwotild/\Mbar)$ is orthogonal to all the
$\pitildi(\Etwotild/\Mbar)$, for $0\le i_1,i_2 \le 2$.  Similarly as for
$\pibarinf$, we can say more about the structure of
$\pitildinf(\Etwotild/\Mbar)$.

\proclaim{Lemma 3.3.11 {\rm (structure of $\pitildinf(\Etwotild/\Mbar)$)}}
 For $c\in\Minf$ let $\Theta_c(\m)$ denote the components of the fibre
$\Etwotild_c$ over ~$c$ (as $\m$ runs through pairs of integers and
pairs of half-integers mod~$N\ZZ$, as in ~2.3.4).  Then in
$\CH^3(\Etwotild\times\Etwotild,\QQ)$,
$$
 \pitildinf(\Etwotild/\Mbar) = \pitildinf\sparen2(\Etwotild/\Mbar) +
\pitildinf\sparen4(\Etwotild/\Mbar),
$$
 with
$$\spreadlines{1\jot} \align
 \pitildinf\sparen2(\Etwotild/\Mbar) &:= \sum_{c\in\Minf} \pitild
c\sparen2(\Etwotild/\Mbar) \\
 \pitildinf\sparen4(\Etwotild/\Mbar) &:= \sum_{c\in\Minf}  \pitild
c\sparen4(\Etwotild/\Mbar)  =
\tr{\vphantom{\pi}}\pitildinf\sparen2(\Etwotild/\Mbar)
\endalign
$$
 where
$$\spreadlines{1\jot} \align
\pitild c\sparen2(\Etwotild/\Mbar) &:= \sum_{\m\in I}
[Z_c(\m)\times_{\{c\}} \Theta_c(\m) ] , \\
 \pitild c\sparen4(\Etwotild/\Mbar) &:= \tr{\vphantom{\pi}}\pitild
c\sparen2(\Etwotild/\Mbar)
  = \sum_{\m\in I} [\Theta_c(\m) \times_{\{c\}} Z_c(\m)] , \\
\endalign
$$
 for some $[Z_c(\m)] \in \Chow2{\Etwotild}$ supported in $\Etwotild_c$.
 Moreover, all the $\pitild c\sparen2(\Etwotild/\Mbar)$ and $\pitild
c\sparen4(\Etwotild/\Mbar)$ are projectors, mutually orthogonal, and also
orthogonal to all $\pitildi(\Etwotild/\Mbar)$, for $0\le i_1,i_2 \le 2$.
\endproclaim

\demo{Proof}
Consider the diagram
$$\spreadmatrixlines{1\jot} \matrix\format%
\r&\enspace\c&\enspace\c&\enspace\c&\enspace\c&\enspace\c&\enspace\l\\
\CH_3(\Etwotildinf\times_{\Minf}\Etwotildinf,\QQ) & \to
&  \CH_3(\Etwotild\times_\Mbar \Etwotild,\QQ) &
\overset{\lup2\psi_1^*}\to\to &
\CH_3(\Etwo\times_M \Etwo,\QQ) &\to & 0  \\
 & \searrow &   @VV (\lup2\psi_2)_* V   @VVV   \\
 & & \CH^3(\Etwotild\times\Etwotild,\QQ) & \to & \CH^3(\Etwo\times
\Etwo,\QQ) &\to & 0
\endmatrix
$$
 whose  horizontal rows are exact \cite{Fulton, 1984, 1.8, p.21}.  Then
$\pitildf(\Etwotild/\Mbar) = (\lup2\psi_2)_*
\pitildf\rel(\Etwotild/\Mbar)$, in
the notation of ~3.3.7, and it follows from Proposition~3.3.9 that the
difference $[\Delta(\Etwotild)]\rel - \pitildf\rel(\Etwotild/\Mbar)$ in
$\CH^3(\Etwotild\times\Etwotild,\QQ)$ maps to zero in $\CH_3(\Etwo\times_M
\Etwo,\QQ)$.  Hence $\pitildinf(\Etwotild/\Mbar)\in
\CH^3(\Etwotild\times\Etwotild,\QQ)$ must be in the image of
$\CH_3(\Etwotildinf\times_{\Minf}\Etwotildinf,\QQ)$.
On the other hand,
the components of $\Etwotildinf\times_{\Minf}\Etwotildinf$ are of the
form $\Theta_c(\m) \times_{\{c\}} \Theta_c(\m')$, which by
Proposition~2.3.1 are products of rational surfaces.  Therefore,
for each of these
components, linear equivalence coincides with homological equivalence,
and thus the K\"unneth formula for homology allows us to conclude that
$\CH_3(\Etwotildinf\times_{\Minf}\Etwotildinf,\QQ)$ is generated by
elements of the form $[\Theta_c(\m)] \times_{\{c\}} [Z_c]$ and $[Z_c']
\times_{\{c\}} [\Theta_c(\m)]$, for $c\in\Minf$ and $\m\in I$ and
$[Z_c],\,[Z_c'] \in \CH_1(\Etwotild_c,\QQ)$.  Hence,
$\pitildinf(\Etwotild/\Mbar)$ can be written in the form claimed.  But
in addition, every class of the form $[\Theta_c(\m)] \times_{\{c\}}
[Z_c]$ is orthogonal to every class of the form $[Z_c'] \times_{\{c\}}
[\Theta_c(\m)]$ for reasons of dimension, and cycles which can be
supported over distinct $c\in\Minf$ are orthogonal, as they are
disjoint.  Therefore all the $\pitild c\sparen j (\Etwotild/\Mbar)$ are
mutually orthogonal.  However, they must also be idempotent and
orthogonal to all the $\pitildi(\Etwotild/\Mbar)$, for
$i_1,i_2 =0,1,2$, since this is true for $\pitildinf(\Etwotild/\Mbar)$.
\Qed
\enddemo

\subsubhead 3.3.12. Splitting $\pitild{1,1}(\Etwotild/\Mbar)$ into
symmetric and antisymmetric parts  \endsubsubhead
  Before leaving this section there is one further refinement we need.
First, observe that fibrewise permutation of the fibre factors of
$\Etwobar\to\Mbar$ preserves the center of the blowing-up ~$\beta$
scheme-theoretically, whence it lifts uniquely to a morphism, say
$\sigma: \Etwotild\to\Etwotild$ of ~$\Etwotild$.  Thus we get an action
of the permutation group $\S_2$ on $\Etwotild$, which together with the
action of $\G^2$ gives a group action of the semidirect product
$\G^2\rtimes\S_2$ on ~$\Etwotild$.  (This is the group $\Gamma_2$ of
\cite{Scholl, 1990, 1.1.1}.)

Next, let
$$\spreadlines{1\jot} \aligned
 A_2 &:= \tfrac 1 2 \big([\Delta(\Etwotild)] + [\tGamma_{\sigma}]\big)\\
 S_2 &:= \tfrac 1 2 \big([\Delta(\Etwotild)] - [\tGamma_{\sigma}]\big)
\endaligned
\qquad \text{\sl in } \Chow 3{\Etwotild\times\Etwotild} .
$$
 Then $A_2$ and $S_2$ are mutually orthogonal projectors whose sum is the
identity in $\Chow 3{\Etwotild\times\Etwotild}$.  Moreover, the
restrictions (in the sense of ~3.2.7 and 3.2.8) $\lup2\psisharp A_2$ and
$\lup2\psisharp S_2$ of $A_2$ and $S_2$ respectively to $\Chow 2
{\Etwo\times_M\Etwo}$, in the notation of Proposition~3.3.8, project the
tensor square of a correspondence in $\Chow 1{\E\times_M\E}$ to its
exterior and symmetric square parts, respectively, cf\.
\cite{K\"unnemann, 1994}, \cite{del Ba\~no Rolla, 1995}.

Now we compose these projectors with $\pitild{1,1}(\Etwotild/\Mbar)$,
and write $\alt2\pitild{1,1}(\Etwotild/\Mbar)$ for $A_2\circ
\pitild{1,1}(\Etwotild/\Mbar)$ and $\sym2\pitild{1,1}(\Etwotild/\Mbar)$
for $S_2\circ\pitild{1,1}(\Etwotild/\Mbar)$.  Then it is easy to check
(by looking in $\QQ[\G^2\rtimes\S_2]$) that $A_2$ and $S_2$ commute with
$\lambdatild\sparen1 \circ \lambdatild\sparen2$ as well as with
$\thetatild\sparen1 \circ \thetatild\sparen2$, and therefore with
$\pitild{1,1}(\Etwotild/\Mbar)$.  Thus, in addition to
$$
 \alt2\pitild{1,1} + \sym2\pitild{1,1} =\pitild{1,1}\; ,
$$
 we also have that $\alt2\pitild{1,1}(\Etwotild/\Mbar)$ and
$\sym2\pitild{1,1}(\Etwotild/\Mbar)$ are orthogonal to each other as
well as to all the $\pitildi(\Etwotild/\Mbar)$, for $(i_1,i_2) \ne
(1,1)$.  Furthermore, from the definitions and Proposition~3.3.8,
$$
 \lup2\psisharp(\sym2\pitild{1,1}(\Etwotild/\Mbar)) =
 \Sym^2_M\pican1(\E/\M) ,
$$
 whereas
$$
 \lup2\psisharp(\alt2\pitild{1,1}(\Etwotild/\Mbar)) =
\twedge^2_M\pican1(\E/\M) \simeq \pican2(\Etwo/M),
$$
 as follows from the definitions, Proposition~3.3.8, and the result of
\cite{Shermenev, 1974} and \cite{K\"unneman, 1994, Thm.3.3.1}.

\head \bf 4. Analysis of the Chow motives $h(\Ebar)$ and $h(\Etwotild)$
\endhead
\noindent
  This section is the technical center of the paper, for here we analyze
the Chow motives determined by the projectors defined in section
three in order to identify them up to isomorphism, when we
can, with Chow motives that can be defined in terms of lower-dimensional
varieties.  For example, we view $\LL^d \simeq (\spec K, \id_K, -d)$ as
being supported on a point, and $h(\Mbar) \simeq \1 \oplus \LL \oplus
h^1(\Mbar)$ as consisting of a constituent submotive belonging
essentially to the curve together with two constituent submotives
supported on points; the precise isomorphisms that we prove in this
section are stated in Theorem~4.2, below.  We reserve exploring the
implications of this theorem for Chow-K\"unneth decompositions and
filtrations
on
the Chow groups of ~$\Ebar$ and ~$\Etwotild$ until the next two sections.

\subhead 4.1. Notation \endsubhead
 Let $\W1 := (\Ebar,\pibar1)$ and $\W2 := (\Etwotild,
\sym2\pitild{1,1}(\Etwotild/\Mbar))$.  Then these are the Chow motives
for modular forms constructed in \cite{Scholl, 1990}, for $k=1,2$; and
modulo homological equivalence, they are the motives for modular forms
constructed in \cite{Deligne, 1969}.

 In the statement of the next theorem, a positive integer coefficient on
a motive indicates the multiplicity with which that motive, up
to isomorphism, occurs.

\proclaim{Theorem 4.2}
As Chow motives in $\scM(K)$,
$$\spreadlines{1\jot}\align
 h(\Ebar) \simeq & \ \1 \oplus m\,\LL \oplus \LL^2
\tag 1 \\
  & \oplus h^1(\Mbar) \oplus (h^1(\Mbar)\tensor \LL) \\
  & \oplus \W1
\endalign
$$
 for some positive integer ~$m$, and
$$\spreadlines{1\jot} \align
h(\Etwotild) \simeq &\ \1 \oplus n\,\LL \oplus n\,\LL^2 \oplus \LL^3
\tag 2 \\
  & \oplus h^1(\Mbar) \oplus 3\,(h^1(\Mbar)\tensor\LL) \oplus
    (h^1(\Mbar)\tensor\LL^2) \\
  & \oplus 2\,(\W1) \oplus 2\,(\W1\tensor \LL) \\
  & \oplus \W2
\endalign
$$
 for some positive integer ~$n$.
\endproclaim

\subsubhead Remark 4.2.1 \endsubsubhead
 It will follow from the proof together with 2.3.2 that $m=\frac12N^2(N-1)
\prod_{p|N}(1-p^{-2})$.  Unfortunately, we don't have equally precise
information about ~$n$.

\subsubhead 4.2.2. Organization of the proof \endsubsubhead
The rest of this section is devoted to the proof of Theorem~4.2, and is
divided into five parts.  In the first part we analyze the motives
defined by $\pibar0$ and $\pibar2$ for $\Ebar$, and by
$\pitild{0,0}(\Etwotild/\Mbar)$ and $\pitild{2,2}(\Etwotild/\Mbar)$ for
{}~$\Etwotild$; these are the constituents of lowest and highest weights.
Then we prove a proposition that describes the action of the extended
relative projectors on the components of the cusp fibres; we need this in
the analysis of all the remaining projectors.  Next we look at the
remaining $\pitildi(\Etwotild/\Mbar)$, for it turns out that they can be
treated together.  After that we describe the motives defined by
$\pibarinf$ and $\pitildinf(\Etwotild/\Mbar)$, and then finally we put
everything together to complete the proof of the theorem.

\subhead 4.3. All zeroes or all twos \endsubhead
 We begin with a little lemma to help deal with the nuisance of
the correction terms occurring in the projectors with zeroes or twos.
This lemma may be compared with the lemma of Beilinson on the lifting of
idempotents by a nilpotent ideal \cite{Jannsen, 1994, p.289}.

\proclaim{Lemma 4.3.1}
 When $X$ is a smooth (connected) projective variety, and $p, p' \in
\Chow{\dim X}{X\times X}$ are projectors such that $(p-p')\circ (p-p')
=0$, then as Chow motives $(X,p) \simeq (X,p')$.
\endproclaim

\demo{Proof}
 The identity map, i.e., $[\Delta(X)]$, induces the isomorphism.
 Write $p'=p+n$, with $n\circ n =0$.  Then from $p\circ p =p$ and
$(p+n)\circ (p+n)=(p+n)$ it is elementary to deduce that
$p\circ(p+n)\circ p =p$ and $(p+n) \circ p \circ (p+n) =(p+n)$, as
required.
\Qed
\enddemo

\proclaim{Proposition 4.3.2}
  As Chow motives in $\scM(K)$,
\roster \smallskip
\item $(\Ebar,\pibar0)\simeq h(\Mbar)$.
\smallskip
\item $(\Ebar,\pibar2)\simeq h(\Mbar) \tensor \LL$.
\smallskip
\item $(\Etwotild, \pitild{0,0}(\Etwotild/\Mbar)) \simeq h(\Mbar)$.
\smallskip
\item $(\Etwotild, \pitild{2,2}(\Etwotild/\Mbar)) \simeq
h(\Mbar)\tensor \LL^2$.
\endroster
\endproclaim

\demo{Proof}
 All the isomorphisms are induced by the graphs or transposed graphs of
the structure maps onto $\Mbar$ and the zero-sections.
 We give first the argument for \therosteritem3, as \therosteritem1 is
similar but simpler.  From the lemma it follows that
$(\Etwotild, \pitild{0,0}(\Etwotild/\Mbar)) \simeq (\Etwotild,
\tgraph{\mutild(0,0)})$ since, as we have observed (3.2.3 and 3.3.2), all
the correction terms are nilpotent of order~2. Then to obtain that
$(\Etwotild, \tgraph{\mutild(0,0)}) \simeq (\Mbar,[\Delta(\Mbar)])$, it
suffices to show
$$\spreadlines{1\jot}\align
 \tgraph{\mubar(0,0)} \circ \tgraph{\phitwotild} \circ [\Delta(\Mbar)]
\circ \tgraph{\alphatild(\bold0)} \circ \tgraph{\mubar(0,0)} &=
\tgraph{\mubar(0,0)} , \\
 [\Delta(\Mbar)] \circ \tgraph{\alphatild(\bold0)} \circ
\tgraph{\mubar(0,0)} \circ \tgraph{\phitwotild} \circ [\Delta(\Mbar)] &=
[\Delta(\Mbar)] .
\endalign
$$
  But these follow from the identities
$$\spreadlines{1\jot} \align
\mutild(0,0)\circ \alphatild(\bold0) \circ \phitwotild \circ \mutild(0,0)
 & = \mutild(0,0) , \\
 \phitwotild \circ \mutild(0,0) \circ \alphatild(\bold0) & = \id_\Mbar ,
\endalign
$$
 where $\alphatild$ is the extension of the level-$N$ structure of
{}~$\Etwo$ to
$\Etwotild$, as in ~3.3.2.
 Now transposing everything proves ~\therosteritem4, and likewise the
correspondences that prove ~\therosteritem2 are the transposes of those
that prove ~\therosteritem1.
\Qed
\enddemo

\subhead 4.4. Action of projectors on fibres and components at infinity
\endsubhead
  Next we consider the action of the our projectors on
fibres and the components of the fibres at infinity.  Roughly speaking,
$\pibarf$ and $\pitildf(\Etwotild/\Mbar)$ annihilate the components of
the cusp fibres---indeed, it was so that this would be the case that
$\pibar1$, and consequently the $\pitildi(\Etwotild/\Mbar)$ with ~$i_1$
or ~$i_2=1$, were chosen as they were---and $\pibarinf$ and
$\pitildinf(\Etwotild/\Mbar)$ act as the identity on those components,
but there are some nuances involving the identity components; the next
proposition gives a precise statement.  As a matter of notation, for any
$t\in\Mbar$ we let $\Ebar_t := \phibar\inv(t)$ and $\Etwotild_t :=
\phitwotild\inv(t)$.  Further, for any cusp $c\in\Minf$ we let
$\theta_c(0)$ be the identity component of $\Ebar_c$, i.e., the
component containing $\alphabar((0,0),c)$, and similarly let
$\Theta_c(\bold0)$ be the identity component of $\Etwotild_c$, the
component containing $\alphatild(\bold0,c)$.

\proclaim{Proposition 4.4.1}
\roster
\item For all $t\in\Mbar$, in $\Chow1{\Ebar}$
$$\spreadlines{1\jot} \align
 \pibar0([\Ebar_t]) &= [\Ebar_t]  \\
 \pibari([\Ebar_t]) &= 0 \quad\text{\sl for } i\ne 0 .
\endalign
$$
\item For all $t\in\Mbar$, in $\Chow1{\Etwotild}$
$$\spreadlines{1\jot} \align
\pitild{0,0}(\Etwotild/\Mbar)([\Etwotild_t]) &= [\Etwotild_t] \\
\pitildi(\Etwotild/\Mbar)([\Etwotild_t]) &= 0  \quad\text{\sl for }
(i_1,i_2)\ne (0,0) .
\endalign
$$
\item  For $c\in\Minf$,  in $\Chow1{\Ebar}$
$$\spreadlines{1\jot} \align
\pibari([\theta_c(m)]) &=0 \quad\text{\sl unless } m=0 \text{ and } i=0,
\\
 \pibar c([\theta_c(m)]) &= [\theta_c(m)] \quad\text{\sl for } m\ne 0, \\
 \pibar0([\theta_c(0)]) &= [\Ebar_c] .
\endalign
$$
\item For $c\in\Minf$,  in $\Chow1{\Etwotild}$
$$\spreadlines{1\jot} \align
\pitildi(\Etwotild/\Mbar)([\Theta_c(\m)]) &= 0 \quad\text{\sl unless }
\m=\0 \text{ and } (i_1,i_2)=(0,0), \\
\pitild c\sparen2(\Etwotild/\Mbar)([\Theta_c(\m)]) &= [\Theta_c(\m)]
\quad\text{\sl for } \m\ne \0 ,\\
\pitild{0,0}(\Etwotild/\Mbar)([\Theta_c(\0)]) &= [\Etwotild_c] .
\endalign
$$
\endroster
\endproclaim

\demo{Proof}
  To begin, we can write $\pibar0 = \tgraph{\mubar(0)} - \frac12
(\psi_2)_*(\phibar\times_\Mbar\phibar)^*(\frak a)$, for a certain class
$\frak
a\in\Chow1{\Mbar}$, see Lemma~3.2.2.  Thus, on any fibre, or any
component of a fibre, $\pibar0$ acts as $\mubar(0)^*$, which acts by
mapping (the class of) the identity component of a fibre to (the class
of) that entire fibre.  By orthogonality, we also get that the other
projectors
defined in section three annihilate the class of an entire fibre.
Similarly
$\pitild{0,0}(\Etwotild/\Mbar)$ acts on fibres or components of fibres as
$\mutild(0,0)^*$, likewise mapping (the class of) the identity component
of any
fibre to (the class of) that entire fibre.  And again, by orthogonality,
we
also get that the other projectors annihilate the class of an entire
fibre.
This proves parts \therosteritem1 and \therosteritem2, and also the
statements
about the action of $\pibar0$ or $\pitild{0,0}(\Etwotild/\Mbar)$ in parts
\therosteritem3 and ~\therosteritem4.

  Next consider $\pibar2 = \graph{\mubar(0)}$ plus a vertical correction
term.  This acts on (the class of) any component of any fibre as
$\mubar(0)_*$, thereby annihilating (the class of) that component.
Similarly $\pitild 2\sparen1(\Etwotild/\Mbar)$ and $\pitild
2\sparen2(\Etwotild/\Mbar)$ act on vertical cycles as $\mutild(0,1)_*$
and $\mutild(1,0)_*$, respectively, from which it follows that they
annihilate vertical two-dimensional cycles, in particular (classes of)
components of fibres.  Thus $\pitildi(\Etwotild/\Mbar)([\Theta_c(\m)]) =
0$ whenever $i_1=2$ or $i_2=2$, for as we saw in the proof of Proposition
{}~3.3.6 we may write $\pitild{i,2}(\Etwotild/\Mbar) = \pitild
i\sparen1(\Etwotild/\Mbar)\circ \pitild 2\sparen2(\Etwotild/\Mbar)$ and
$\pitild{2,i}(\Etwotild/\Mbar) = \pitild i\sparen2(\Etwotild/\Mbar)\circ
\pitild 2\sparen1(\Etwotild/\Mbar)$.

  Now consider $\pibar1 = \lambdabar\circ\thetabar$, as in ~3.2.5.  Then
$\thetabar$ acts on a component $\theta_c(m)$ of $\Einf$ by
$\thetabar([\theta_c(m)]) = \frac 1 N [\Ebar_c]$, while
$\lambdabar([\theta_c(m)]) = \frac12([\theta_c(m)] -[\theta_c(-m)])$, as
follows from 2.2.1 and 3.2.5.  So it is easy to see that their combined
effect is to annihilate any ~$[\theta_c(m)]$.

  Finally we consider $\pitildi(\Etwotild/\Mbar)$ with $i_1=1$ or
$i_2=1$; as above, we will be finished if we can show that
$\pitild 1\sparen j(\Etwotild/\Mbar) ([\Theta_c(\m)]) =0$ for any
component $\Theta_c(\m)$ of $\Etwotildinf$, for $j=1$ or~2.  For
definiteness, suppose for the moment that $j=2$, and write $\pitild
1\sparen 2(\Etwotild/\Mbar) = \lambdatild\sparen2 \circ
\thetatild\sparen2$, as in ~3.3.3.  Then letting $\Theta_c(m,n)$
represent the components of $\Etwotild_c$, with the indexing described
in ~2.3.4, we find that $\lambdatild\sparen2 ([\Theta_c(m,n)]) =
\frac 12 ([\Theta_c(m,n)]-[\Theta_c(m,-n)])$, and $\thetatild\sparen2
([\Theta_c(m,n)]) = \frac 1 N \sum_{n\in\ZmodN} [\Theta_c(m,n)]$.  Thus
the combined effect of the two is to annihilate $[\Theta_c(m,n)]$, as
required.  Since the argument is the same $j=1$, this completes the proof.
\Qed
\enddemo

\subhead 4.5. Isomorphisms between submotives of $h(\Etwotild)$ and
submotives of $h(\Ebar)$ \endsubhead
   The next proposition identifies several of the motivic constituents of
$h(\Etwotild)$ defined by the projectors defined in section three with
motives
supported on lower dimensional varieties.  Although some of these can be
supported on $\Mbar$ or $\spec K$, what we actually verify is that some
of the submotives of $\Etwotild$ are isomorphic to submotives of $\Ebar$,
so we state the proposition this way and defer further reduction until
the last part of the section.

\proclaim{Proposition 4.5.1}
  As Chow motives in $\scM(K)$,
\roster
\item  $(\Etwotild,\pitild{0,1}(\Etwotild/\Mbar)) \simeq
(\Etwotild,\pitild{1,0}(\Etwotild/\Mbar)) \simeq (\Ebar,\pibar1)$;
\medskip
\item $(\Etwotild,\pitild{0,2}(\Etwotild/\Mbar)) \simeq
(\Etwotild,\pitild{2,0}(\Etwotild/\Mbar)) \simeq (\Ebar,\pibar2)$;
\medskip
\item $(\Etwotild,\alt2\pitild{1,1}(\Etwotild/\Mbar)) \simeq
(\Ebar,\pibar2)$;
\medskip
\item  $(\Etwotild,\pitild{1,2}(\Etwotild/\Mbar)) \simeq
(\Etwotild,\pitild{2,1}(\Etwotild/\Mbar)) \simeq (\Ebar,\pibar1,-1)$.
\endroster
\endproclaim

\demo{Proof}
  Since the proofs of these isomorphisms between a submotive of
$\Etwotild$
and a submotive of $\Ebar$ all follow a similar pattern, when an
argument applies generally we will use $\tildepi$ to represent any of the
seven projectors for $\Etwotild$ above, and $\barpi$ for the
corresponding projector on $\Ebar$, and $m$ for the corresponding Tate
twist (when present); but when the differences in detail require it,
we will refer to the specific cases \therosteritem1--\therosteritem4.

 With this notation and that in Propositions ~3.2.8 and ~3.3.9, the
first observation is that in each case there is an isomorphism
$(\Etwo,\lup2\psisharp\tildepi) \simeq (\E,\psisharp\barpi, m)$ of
relative Chow motives over~$M$.  For parts ~\therosteritem1,
\therosteritem2 and \therosteritem4 this follows more or less formally
from the tensor structure of the category $\scM(M)$, as in ~1.1.2 and
{}~1.1.3; whereas for part ~\therosteritem3 it follows from the theorem
of \cite{Shermenev, 1974} and \cite{K\"unnemann, 1994, Thm.3.3.1}, as
mentioned in ~3.3.12.  So let $\alpha$ on $\Etwo\times_M\E$ and $\beta$
on $\E\times_M\Etwo$ be cycles inducing this isomorphism in each
direction, and let $\alphatild$ and $\betatild$ denote their closures in
$\Etwotild \times \Ebar$ and $\Ebar\times\Etwotild$, respectively.  Then
we claim that $[\alphatild] \in \Corr^{m}(\Etwotild, \Ebar)$ and
$[\betatild] \in \Corr^{-m}(\Ebar,\Etwotild)$ induce inverse isomorphisms
between $(\Etwotild,\tildepi)$ and $(\Ebar,\barpi,m)$.  To verify this we
must show that
$$\spreadlines{1\jot}\align
 \barpi \circ [\alphatild] \circ \tildepi \circ [\betatild] \circ
\barpi &= \barpi
\tag 5
\\
 \tildepi \circ [\betatild] \circ \barpi  \circ [\alphatild] \circ
\tildepi &= \tildepi .
\tag 6
\endalign
$$

  Consider (5) first.  What we already know is that the
correspondences on both sides of the equation can be supported on
$\Ebar\times_\Mbar\Ebar$, and that their restrictions (in the sense of
3.2.7 and~3.2.8) to $\E\times_M\E$
$$
 \psi_1^*\big( \barpi \circ [\alphatild] \circ \tildepi \circ [\betatild]
\circ \barpi \big) =  \psi_1^*\barpi
$$
 coincide.  Thus the exactness of the sequence
$$
 \CH_2(\Ebarinf\times_{\Minf}\Ebarinf,\QQ) \to
\CH_2(\Ebar\times_\Mbar\Ebar,\QQ) \to \CH_2(\E\times_M\E,\QQ) \to 0
$$
  implies that the difference
$$
(\barpi \circ [\alphatild] \circ \tildepi \circ [\betatild]
\circ \barpi) - \barpi = \sum \Sb{c\in\Minf} \\ {m,n\in\ZmodNsq} \endSb
a_c(m,n) [\theta_c(m)]\times_{\{c\}} [\theta_c(n)] ,
$$
 since it lies in the image of
$\CH_2(\Ebarinf\times_{\Minf}\Ebarinf,\QQ)$ in
$\CH_2(\Ebar\times_\Mbar\Ebar,\QQ)$.  Then composing with $\barpi$ on the
left and right gives
$$
\spreadlines{1\jot} \align
(\barpi \circ [\alphatild] \circ \tildepi \circ [\betatild]
\circ \barpi) - \barpi & = \sum \Sb{c\in\Minf} \\{m,n\in\ZmodNsq} \endSb
a_c(m,n) \tr{\vphantom{\pi}}\barpi([\theta_c(m)]) \times_{\{c\}}
\barpi([\theta_c(n)])
\\ & =0
\endalign
$$
 by applying Proposition ~4.4.1.  This proves (5).

 The argument for (6), using the right-exact sequence
$$
 \CH_3(\Etwotildinf\times_{\Minf}\Etwotildinf,\QQ) \to
\CH_3(\Etwotild\times_\Mbar\Etwotild,\QQ) \to
\CH_3(\Etwo\times_M\Etwo,\QQ) \to 0 ,
$$
 runs in a completely parallel manner up to the point where
$$\spreadlines{1\jot} \align
 &(\tildepi \circ [\betatild] \circ \barpi  \circ [\alphatild] \circ
\tildepi) - \tildepi \\
 =\ &\sum \Sb{c\in\Minf}\\{\m\in I} \endSb
 \Big( a_c(\m) \big([Z'_c(\m)]\times_{\{c\}} [\Theta_c(\m)]\big) +
b_c(\m) \big([\Theta_c(\m)] \times_{\{c\}} [Z''_c(\m)]\big)\Big ) ,
\endalign
$$
 for some one-cycles $Z'_c(\m),\,Z''_c(\m)$ on $\Etwotild_c$ and rational
numbers $a_c(\m),\,b_c(\m)$, with $\Theta_c(\m)$ running over the
components of ~$\Etwotildinf$ and $I$ the indexing described in ~2.3.4.
Now composing with $\tildepi$ on both left and right leaves $(\tildepi
\circ [\betatild] \circ \barpi \circ [\alphatild] \circ \tildepi) -
\tildepi$ fixed, but on the other terms,
$$\spreadlines{1\jot} \alignat2
 \tildepi\circ\big([Z'_c(\m)]\times_{\{c\}} [\Theta_c(\m)]\big)\circ
\tildepi &= \tr{\vphantom{pi}}\tildepi([Z'_c(\m)]) \times_{\{c\}}
\tildepi([\Theta_c(\m)]) &&=0 \\
 \tildepi\circ\big([\Theta_c(\m)] \times_{\{c\}} [Z''_c(\m)]\big)\circ
\tildepi &= \tr{\vphantom{pi}}\tildepi([\Theta_c(\m)]) \times_{\{c\}}
\tildepi([Z''_c(\m)]) &&=0 ,
\endalignat
$$
 since $\tr{\vphantom{pi}}\tildepi = \tildepi$ and
$\pitild{0,0}(\Etwotild/\Mbar) \ne \tildepi \ne
\pitild{2,2}(\Etwotild/\Mbar)$, so that Proposition ~4.4.1 applies.  This
proves (6), and concludes the proof of the proposition.
\Qed
\enddemo

\subhead 4.6. The motives defined by $\pibarinf$ and
$\pitildinf(\Etwotild/\Mbar)$ \endsubhead
  Finally we must analyze the motives defined by $\pibarinf$ and
$\pitildinf(\Etwotild/\Mbar)$.  Since these were each defined as the
difference between the diagonal and the sum of the $\pibari$ or
$\pitildi(\Etwotild/\Mbar)$
respectively, it requires some care to get a good grip on them.  However,
in the end the motives themselves have a rather simple form, as a sum
of powers of Lefschetz motives, essentially because all the components
of the cusp fibres supporting these projectors are rational varieties.

\proclaim{Proposition 4.6.1}
 As Chow motives in $\scM(K)$,
\roster
\item $(\Ebar,\pibar c) \simeq (N-1) \LL$;
\medskip
\item $(\Etwotild,\pitild c\sparen2(\Etwotild/\Mbar)) \simeq s\LL$;
\medskip
\item $(\Etwotild,\pitild c\sparen4(\Etwotild/\Mbar)) \simeq s\LL^2$
\endroster
 for $c\in\Minf$ and some $0<s\in\ZZ$.
\endproclaim

\demo{Proof}
  We give first the proof for \therosteritem2 and \therosteritem3, which
come together, and comment at the end on ~\therosteritem1, since it can
be proved similarly, and even more easily.  For convenient reference,
recall that $\pitild c\sparen2(\Etwotild/\Mbar)$ and $\pitild
c\sparen4(\Etwotild/\Mbar)$ respectively have the form
$$
 \pitild c\sparen2(\Etwotild/\Mbar) = \sum_{\m\in I}
[Z_c(\m)\times_{\{c\}} \Theta_c(\m) ] =
\tr{\vphantom{\pi}}\pitild c\sparen4(\Etwotild/\Mbar)
$$
 where the $Z_c(\m)$ are some one-cycles supported on $\Etwotild_c$,
about which {\sl\`a~priori\/} we know nothing else, and $I$ is the
indexing described in ~2.3.4.  The proof will proceed in several steps.

 {\it Step one.}
 Firstly, we claim that $\m=\0$, if it occurs, can be eliminated from
the expression for $\pitild c\sparen2(\Etwotild/\Mbar)$ and $\pitild
c\sparen4(\Etwotild/\Mbar)$, where (as in the proof of Proposition
{}~4.4.1) $\Theta_c(\0)$ denotes the identity component of $\Etwotild_c$.
For observe that in $\Chow2{\Etwotild}$ the class of the fibre over
$c\in\Minf$ can be written as
$$
 [\Etwotild_c] = [\Theta_c(\0)] +
\sum_{\m\ne\0} [\Theta_c(\m)] ,
$$
 where the sum runs over all components of $\Etwotild_c$ other than the
identity component.  Then using this to give an alternate expression for
$[\Theta_c(\0)]$, we rewrite
$$
 \pitild c\sparen2(\Etwotild/\Mbar) = [Z_c(\0)]\times_{\{c\}}
[\Etwotild_c] +
\sum_{\m\ne\0} [Z'_c(\m)]\times_{\{c\}} [\Theta_c(\m) ] ,
\tag 4
$$
 for suitable one-cycles $Z'_c(\m)$ supported on ~$\Etwotild_c$.  Next,
there exists $\frak d(c) \in \Chow1{\Mbar}$ rationally equivalent to
{}~$[c]$ but with support disjoint from ~$\Minf$.  Since from Proposition
{}~4.4.1 we know that $\pitildinf(\Etwotild/\Mbar)$ annihilates
$[\Etwotild_x]$ for any $x\in M$, it follows that
$$
 \pitildinf\sparen2(\Etwotild/\Mbar)([\Etwotild_c]) =
\pitildinf\sparen2(\Etwotild/\Mbar)(\phitwotild^*(\frak d(c))) = 0 .
\tag 5
$$
  We also know from Proposition ~4.4.1(4) that
$$
\pitild c\sparen2(\Etwotild/\Mbar) ([\Theta_c(\m)]) = [\Theta_c(\m)]
 \qquad\text{for } \m\ne\0 .
\tag 6
$$
 Now we compose both sides of ~(4) with $\pitild
c\sparen2(\Etwotild/\Mbar)$ on
the left.  Since $\pitild c\sparen2(\Etwotild/\Mbar)$ is idempotent, the
left-hand side is unchanged.  As for the right-hand side, from (5), (6)
and the
general observation that the composition of a correspondence $\pi$ with a
correspondence of the form $[Z]\times[T]$ is $\pi \circ ([Z]\times[T]) =
[Z]\times \pi([T])$, we conclude that
$$
 \pitild c\sparen2(\Etwotild/\Mbar)  = \sum_{\m\ne\0}
[Z'_c(\m)\times_{\{c\}} \Theta_c(\m) ] = \tr{\vphantom{\pi}}\pitild
c\sparen4(\Etwotild/\Mbar) ,
$$
 with $Z'_c(\m)$ as in~(4).  Thus we have an expression for $\pitild
c\sparen2(\Etwotild/\Mbar)$ with no $\m =\0$ term, as claimed.  Indeed, by
comparing this with ~(4), it follows that $[Z_c(\0)] =0$.

{\it Step two.}
 Next, we claim that without loss of generality, we can replace the
one-cycles $[Z'_c(\m)]$ by one-cycles $[Z''_c(\m)]$ with the property
that
$$
\pitild c\sparen4(\Etwotild/\Mbar)([Z''_c(\m)]) = [Z''_c(\m)] ,
$$
 where, by virtue of step one, $\m\ne\0$.  In fact, if we replace
$[Z'_c(\m)]$ by
$$
 [Z''_c(\m)] := [Z'_c(\m)] - \pitildf(\Etwotild/\Mbar)([Z'_c(\m)])
$$
 in the last expression for $\pitild c\sparen2(\Etwotild/\Mbar)$, then
$\pitild c\sparen2(\Etwotild/\Mbar)$ remains unchanged.  For using the
orthogonality of
$\tr{\vphantom{\pi}}\pitildf(\Etwotild/\Mbar) =
\pitildf(\Etwotild/\Mbar)$ with $\pitild c\sparen2(\Etwotild/\Mbar)$, we
can write
$$\spreadlines{1\jot} \align
 &\ \pitild c\sparen2(\Etwotild/\Mbar) = \pitild
c\sparen2(\Etwotild/\Mbar) - \pitild c\sparen2(\Etwotild/\Mbar)\circ
\tr{\pitildf}(\Etwotild/\Mbar) \\
 = & \sum_{\m\ne\0} [Z'_c(\m)] \times_{\{c\}} [\Theta_c(\m)] -
\sum_{\m\ne\0} \pitildf(\Etwotild/\Mbar)([Z'_c(\m)])\times_{\{c\}}
[\Theta_c(\m)] ,
\endalign
$$
 from which it follows that
$$
 \pitild c\sparen2(\Etwotild/\Mbar) =  \sum_{\m\ne\0} [Z''_c(\m)]
\times_{\{c\}} [\Theta_c(\m)] = \tr{\vphantom{\pi}}\pitild c
\sparen4(\Etwotild/\Mbar).
$$
 Furthermore, it's clear that $\pitildf(\Etwotild/\Mbar)([Z''_c(\m)]) =
0$, from which it follows by orthogonality that $\pitild
c\sparen4(\Etwotild/\Mbar)([Z''_c(\m)]) = [Z''_c(\m)]$.

 {\it Step three.}
  Next we claim that the Chow groups of the motives defined by $\pitild
c\sparen2(\Etwotild/\Mbar)$ and $\pitild c\sparen4(\Etwotild/\Mbar)$
respectively are
$$\spreadlines{1\jot} \align
 \CH ((\Etwotild,\pitild c\sparen2(\Etwotild/\Mbar)),\QQ) &= \Span_\QQ\{
[\Theta_c(\m)]\mid \m\ne\0 \} , \\
 \CH ((\Etwotild,\pitild c\sparen4(\Etwotild/\Mbar)),\QQ) &= \Span_\QQ\{
[Z''_c(\m)]\mid \m\ne\0 \} ,
\endalign
$$
 and thus, in particular, these are finite-dimensional vector spaces.
 For the right-hand sides are contained in the left-hand sides because
by Proposition ~4.4.1 $\pitild c\sparen2(\Etwotild/\Mbar)$ acts on
$[\Theta_c(\m)]$ as the identity for $\m\ne\0$, and similarly, by step
two above $\pitild c\sparen4(\Etwotild/\Mbar)$ acts on $[Z''_c(\m)]$ as
the identity,
$\m\ne\0$.  On the other hand, $\CH ((\Etwotild,\pitild
c\sparen2(\Etwotild/\Mbar)),\QQ)$ is contained in the span of the
$[\Theta_c(\m)]$ other than the identity component because for any cycle
{}~$\xi$ on ~$\Etwotild$,
$$\spreadlines{1\jot} \align
 \pitild c\sparen2(\Etwotild/\Mbar)(\xi) &=
\pr_{2*}\big[(\xi\times\Etwotild) \cdot \big(\sum_{\m\ne\0}  Z''_c(\m)
\times  \Theta_c(\m) \big)\big] \\
 &= \sum_{\m\ne\0} (\xi\cdot Z''_c(\m))\,[\Theta_c(\m)] .
\endalign
$$
 The inclusion of $\CH ((\Etwotild,\pitild
c\sparen4(\Etwotild/\Mbar)),\QQ)$ in the span of the $[Z''_c(\m)]$ for
$\m\ne\0$ follows similarly.

 {\it Step four.}
 We claim that the intersection pairing on
{}~$\Etwotild$ restricts nondegenerately to a pairing
$$
 \CH^1((\Etwotild,\pitild c\sparen2(\Etwotild/\Mbar)),\QQ) \tensor
\CH^2((\Etwotild,\pitild c\sparen4(\Etwotild/\Mbar)),\QQ) \lra
\CH^3(\Etwotild,\QQ) \simeq \QQ .
$$
  For any $[\Theta] \in \CH^1((\Etwotild,\pitild
c\sparen2(\Etwotild/\Mbar)),\QQ)$, consider
$$\spreadlines{1\jot}\align
 [\Theta] &= \pitild c\sparen2(\Etwotild/\Mbar) ([\Theta]) \\
  & = \pr_{2*}\big[(\Theta\times\Etwotild) \cdot
\big(\sum_{\m\ne\0}  [Z''_c(\m)  \times  \Theta_c(\m)]\big)\big] \\
 &= \sum_{\m\ne\0} (\Theta\cdot Z''_c(\m))\,[\Theta_c(\m)] .
\endalign
$$
 Thus, unless it is already zero, $[\Theta]$ cannot be
orthogonal to all $[Z''_c(\m)]$ for $\m\ne\0$.
Similarly, no $[Z] \in \CH^2((\Etwotild,\pitild
c\sparen4(\Etwotild/\Mbar)),\QQ)$ can be orthogonal to
all $[\Theta_c(\m)]$ for $\m\ne\0$.

 It also follows that $\CH^1((\Etwotild,\pitild
c\sparen2(\Etwotild/\Mbar)),\QQ)$ and $\CH^2((\Etwotild,\pitild
c\sparen4(\Etwotild/\Mbar)),\QQ)$ must have the same dimension.

 {\it Conclusion of the proof for parts \therosteritem2 and
{}~\therosteritem3.}
  Now choose any convenient basis for $\CH^1((\Etwotild,\pitild
c\sparen2(\Etwotild/\Mbar)),\QQ)$, say $\{\omega_l\mid l=1,\dots,s\}$,
for some~$s$, and replace each $[\Theta_c(\m)]$ in the last expression
for $\pitild c\sparen2(\Etwotild/\Mbar)$ by a linear combination of these
{}~$\omega_l$.  The outcome is then
$$
 \pitild c\sparen2(\Etwotild/\Mbar) = \sum_{l=1}^s \zeta_l\times\omega_l
= \tr{\vphantom{\pi}}\pitild c\sparen4(\Etwotild/\Mbar)
$$
 for some $\zeta_l \in \CH^2((\Etwotild,\pitild
c\sparen4(\Etwotild/\Mbar)),\QQ)$.  Then for
$1\le l_0\le m$ we have, similarly as above,
$$\spreadlines{1\jot}\align
 \omega_{l_0} &= \pitild c\sparen2(\Etwotild/\Mbar) (\omega_{l_0}) \\
  & = \pr_{2*}\big((\omega_{l_0}\times\Etwotild) \cdot
\big(\sum_{l=1}^s \zeta_l \times \omega_l\big)\big) \\
 &= \sum_{l=1}^s (\omega_{l_0}\cdot \zeta_l)\,\omega_l .
\endalign
$$
 But since $\{\omega_l,\, 1\le l\le s\}$ is a basis of
$\CH^1((\Etwotild,\pitild c\sparen2(\Etwotild/\Mbar)),\QQ)$, the
intersection multiplicity
$$
(\omega_{l_0}\cdot \zeta_l) = \cases 1 &\text{when } l=l_0, \\
 0 &\text{when } l\ne l_0 . \endcases
$$
 This means that $\{\zeta_l,\, 1\le l\le s\}$ is the dual basis of
$\CH^2((\Etwotild,\pitild c\sparen4(\Etwotild/\Mbar)),\QQ)$, and that the
individual terms in the expression above for $\pitild
c\sparen2(\Etwotild/\Mbar)$ and $\pitild c\sparen4(\Etwotild/\Mbar)$ are
mutually orthogonal idempotents.  And as we saw in ~1.1.2(c), projectors
of this form define powers of Lefschetz motives.  Thus the motives
defined by $\pitild c\sparen2(\Etwotild/\Mbar)$ and $\pitild
c\sparen4(\Etwotild/\Mbar)$ have the form asserted.

 {\it Proof of part \therosteritem1.}
 The proof of part \therosteritem1 can be carried out in the same way,
with a few small differences and simplifications.  Starting with the
expression
$$
 \pibar c = \sum_{m,n\in\ZmodN} r_c(m,n) [\theta_c(m)]\times_{\{c\}}
[\theta_c(n)],
$$
 for some $r_c(m,n)\in\QQ$, the same argument as step one applied twice
leads to
$$
 \pibar c = \sum_{m\ne 0,\,n\ne 0} s_c(m,n) [\theta_c(m)]\times_{\{c\}}
[\theta_c(n)],
$$
 for some $s_c(m,n)\in\QQ$.  Then steps three and four are replaced and
made more precise by \cite{Shioda, 1972, Thm.1.1 and Lemma~1.3}, which
imply that $\{[\theta_c(m)]\mid 0\ne m\in\ZmodN\}$ is already
algebraically independent and has a nondegenerate intersection matrix
$((\theta_c(m)\cdot\theta_c(n)))$, i.e., of rank ~$(N-1)$, see
remark~2.2.2.  From this it follows that $(s_c(m,n))$ is the inverse of
the intersection matrix.  Then if we rewrite
$$
 \pibar c = \sum_{0\ne m\in\ZmodN}  [\theta_c(m)]\times_{\{c\}}
\big(\sum_{n\ne 0} s_c(m,n)[\theta_c(n)]\big),
$$
 we see $\pibar c$ as the sum of $(N-1)$ mutually orthogonal projectors
of the form $[A]\times[B]$ with $(A\cdot B)=1$.  This proves part
{}~\therosteritem1, and concludes the proof of the proposition.
\Qed
\enddemo

\bigpagebreak
\subhead 4.7. Proof of Theorem~4.2 \endsubhead
  Now we prove Theorem ~4.2.  Consider first
$\Ebar$:  From Propositions ~4.3.2 and ~4.6.1 we get
$$\spreadlines{1\jot} \align
h(\Ebar) &\simeq (\Ebar,\pibar0) \oplus (\Ebar,\pibar1) \oplus (\Ebar,
\pibar2) \oplus (\Ebar,\pibarinf) \\
  & \simeq h(\Mbar) \oplus \W1 \oplus (h(\Mbar)\tensor\LL) \oplus r\LL ,
\endalign
$$
  where it follows from ~4.6.1 that $r=(N-1)\cdot\#(\Minf)$.  Then by
using
that
$$
 h(\Mbar) \simeq \1 \oplus \LL \oplus h^1(\Mbar) ,
$$
 the decomposition asserted in the statement of the theorem follows.
The argument for $h(\Etwotild)$, using in addition Proposition ~4.5.1, is
entirely similar.
\Qed


\head \bf 5. Chow-K\"unneth decompositions and the cohomology of $\Ebar$
and $\Etwotild$
\endhead
\noindent
  We can now give two proofs of the existence of Chow-K\"unneth
decompositions for $\Ebar$ and $\Etwotild$.  The first proof very quickly
deduces the existence and a description of the Chow-K\"unneth
decompositions for $\Ebar$ and $\Etwotild$ from Theorem ~4.2 using
\cite{Scholl, 1990, Thm.1.2.1} to tell us the cohomology of $\W1$ and
{}~$\W2$.  The second proof also starts with Theorem ~4.2, but then uses a
description of the total cohomology spaces $\Hdot (\Ebar,\Qdot)$ and
$\Hdot (\Etwotild,\Qdot)$ to obtain the Chow-K\"unneth decompositions
for $\Ebar$ and ~$\Etwotild$, and at the same time compute the cohomology
of ~$\W1$ and ~$\W2$, i.e., the cases $k=1$ and $k=2$ of
\cite{Scholl, 1990, Thm.1.2.1}.

 Recall that a positive integer coefficient on a motive indicates the
multiplicity with which that motive, up to isomorphism, occurs.

\proclaim{Theorem 5.1}
 With $m$ and $n$ as in Theorem~4.2,
\nobreak\par\nobreak
\noindent {\rm(1)}\enspace  $\Ebar$ has a Chow-K\"unneth decomposition,
with
$$\spreadlines{1\jot} \alignat2
h^0(\Ebar) & \simeq \1 & h^4(\Ebar) &\simeq \LL^2 \\
h^1(\Ebar) & \simeq h^1(\Mbar) &  h^3(\Ebar) &\simeq h^1(\Mbar)\tensor
\LL \\
h^2(\Ebar) &\simeq m\LL \oplus \W1 \qquad&&
\endalignat
$$
\par
\noindent {\rm(2)}\enspace  $\Etwotild$ has a Chow-K\"unneth
decomposition, with
$$\spreadlines{1\jot} \alignat2
 h^0(\Etwotild) &\simeq \1 &  h^6(\Etwotild) &\simeq \LL^3 \\
 h^1(\Etwotild) &\simeq h^1(\Mbar) & h^5(\Etwotild) &\simeq
h^1(\Mbar)\tensor \LL^2 \\
 h^2(\Etwotild) &\simeq n \LL \oplus 2(\W1) & h^4(\Etwotild) &\simeq
n \LL^2 \oplus 2(\W1\tensor \LL) \\
 h^3(\Etwotild) &\simeq 3(h^1(\Mbar)\tensor\LL) \oplus \W2 \qquad&&
\endalignat
$$
\endproclaim

\subsubhead Remark 5.1.1 \endsubsubhead
 The existence of Chow-K\"unneth decompositions for surfaces in general
is proved in \cite{Murre, 1990}.  Proposition~5.1 describes what it
looks like specifically for ~$\Ebar$, and also gives a more refined
decomposition for this surface.

\demo{5.1.2. The first proof}
 After Lemma ~1.2.5, it is only necessary to verify that all of the
submotives given by theorem~4.2 have \CKd s.  It is clear that $\LL^d$
has a \CKd, and easy to see that $h(\Mbar)\tensor\LL^d$ does, as well.
But $\W1$ and $\W2$ also have \CKd s, for by \cite{Scholl, 1990,
Thm.1.2.1},
$$\spreadlines{1\jot} \align
 \Hdot (\W1,\Qdot) &\simeq \Hdot^1(\Mbar, j_*R^1\phi_*\Qdot) \subset
\Hdot^2(\Ebar,\Qdot) \\
\Hdot (\W2,\Qdot) &\simeq \Hdot^1(\Mbar, j_*\Sym^2R^1\phi_*\Qdot)
\subset \Hdot^3(\Etwotild,\Qdot),
\endalign
$$
 which means in particular that the cohomology of ~$\W1$ is purely of
weight ~2, so $\id_{(\W1)} = \pi_2(\W1)$ is a \CKd\space for ~$\W1$,
and similarly the cohomology of ~$\W2$ is purely of weight ~3, so
$\id_{(\W2)} = \pi_3(\W2)$ is a \CKd\space for ~$\W2$.  Thus Theorem
{}~4.2 gives $h(\Ebar)$ and $h(\Etwotild)$ respectively as
direct sums of motives with \CKd s, therefore by Lemma ~1.2.5, both
$\Ebar$ and $\Etwotild$ have \CKd s.  By collecting together the
components of each given weight, we get the \CKd s for $\Ebar$ and
$\Etwotild$ as claimed.
\Qed
\enddemo

\subhead 5.2. The cohomology of $\Ebar$ and $\Etwotild$ \endsubhead
  In the proof just given, the nontrivial cohomology computations
were already taken care of by \cite{Scholl, 1990, Thm.1.2.1}.  But we
can also prove the existence of \CKd s for $\Ebar$ and $\Etwotild$
independently of that theorem, while at the same time computing the
cohomology of ~$\W1$ and $\W2$, which are the cases $k=1$ and
$k=2$ of \cite{Scholl, 1990, Thm.1.2.1}.  Toward this end, we recall some
facts about the cohomology of $\Ebar$ and $\Etwotild$.

\proclaim{Proposition 5.2.1}
$$\spreadlines{\medskipamount} \align
 \Hdot (\Ebar,\Qdot) \simeq & \bigoplus_{p=0}^2 \big( \Hdot^p(\Mbar,
\Qdot) \oplus \Hdot^p(\Mbar,\Qdot(-1))\big) \tag 1 \\
 & \oplus \Hdot^1(\Mbar, j_*R^1\phi_*\Qdot) \oplus
\Hdot^2\strut_{\Minf}(\Mbar, \U_\infty)
\endalign
$$
 where $\U_\infty$ is a skyscraper sheaf supported over $\Minf$ that
contributes to cohomology only in degree ~2.  Moreover, the intersection
form on $\Ebar$ induces perfect pairings
$$\spreadlines{1\jot}\align
 \Hdot^p(\Mbar,\Qdot(j)) &\tensor \Hdot^{2-p}(\Mbar,\Qdot(-(j+1))) \quad
\text{for } 0\le p\le 2 \\
 \Hdot^1(\Mbar, j_*R^1\phi_*\Qdot) &\tensor \Hdot^1(\Mbar,
j_*R^1\phi_*\Qdot) \\
 \Hdot^2\strut_{\Minf}(\Mbar, \U_\infty) &\tensor
\Hdot^2\strut_{\Minf}(\Mbar, \U_\infty)
\endalign
$$
 into $\Hdot^4(\Ebar,\Qdot)\simeq \Qdot(-2)$.
$$\spreadlines{\medskipamount} \align
 \Hdot (\Etwotild,\Qdot) \simeq & \bigoplus_{p=0}^2 \big(
\Hdot^p(\Mbar,\Qdot) \oplus 3\Hdot^p(\Mbar,\Qdot(-1)) \oplus
\Hdot^p(\Mbar,\Qdot(-2)) \big)  \tag 2 \\
 & \oplus \ 2\big( \Hdot^1(\Mbar,j_*R^1\phi_*\Qdot) \oplus
\Hdot^1(\Mbar,j_*R^1\phi_*\Qdot(-1)) \big) \\
 & \oplus \ \Hdot^2\strut_{\Minf}(\Mbar,\U_\infty\sparen2) \oplus
\Hdot^4\strut_{\Minf}(\Mbar,\U_\infty\sparen4)
\endalign
$$
 where $\U_\infty\sparen j$ is a skyscraper sheaf supported over $\Minf$
that contributes to cohomology only in degree~$j$, for $j=2,4$.
Moreover, the intersection form on $\Etwotild$ induces perfect pairings
into $\Hdot^6(\Etwotild,\Qdot)\simeq \Qdot(-3)$ on the isotypic
components corresponding to
$$\spreadlines{1\jot} \align
 \Hdot^p(\Mbar,\Qdot(j)) &\tensor \Hdot^{2-p}(\Mbar,\Qdot(-(j+2))) \\
 \Hdot^1(\Mbar, j_*R^1\phi_*\Qdot(j)) &\tensor \Hdot^1(\Mbar,
j_*R^1\phi_*\Qdot(-(j+1))) \\
 \Hdot^2\strut_{\Minf}(\Mbar,\U_\infty\sparen2) &\tensor
\Hdot^4\strut_{\Minf}(\Mbar,\U_\infty\sparen4) .
\endalign
$$
\endproclaim

\demo{Proof}
 All of this is well-known, but as we do not know of a convenient
reference, we sketch the argument for $\Etwotild$, the argument for
$\Ebar$ being similar.  Firstly, the decomposition theorem of
\cite{Beilinson et al., 1983} implies that
$$ 
 \Hdot (\Etwotild,\Qdot) \simeq \bigoplus_{p=0}^2 \bigoplus_{q=0}^4
\Hdot^p(\Mbar, j_*R^q(\phitwo)_*\Qdot) \oplus \bigoplus_{s=1}^2
\Hdot^{2s}\strut_{\Minf}(\Mbar, \U_\infty\sparen{2s}) ,
$$
 where $\U_\infty\sparen{2s}$ is a skyscraper sheaf supported on $\Minf$
contributing in degree ~$2s$, as well as the Poincar\'e duality pairings
$$ \spreadlines{1\jot} \align
 \Hdot^p(\Mbar, j_*R^q(\phitwo)_*\Qdot) &\tensor \Hdot^{2-p}(\Mbar,
j_*R^{4-q}(\phitwo)_*\Qdot) \\
 \Hdot^{2}\strut_{\Minf}(\Mbar, \U_\infty\sparen{2}) &\tensor
\Hdot^{4}\strut_{\Minf}(\Mbar, \U_\infty\sparen{4}) .
\endalign
$$
 The next observation is that as a sheaf on $\M$,
$$
 R^q(\phitwo)_*\Qdot \simeq \bigoplus_{r=0}^2 m(2,q,r) \Sym^r
R^1\phi_*\Qdot(\tfrac{r-q}2) ,
$$
where
$$
 m(2,q,r) := \binom 2 {\frac{q-r}2} \binom 2 {\frac{q+r}2} - \binom 2
{\frac{q-r}2 -1}\binom 2 {\frac{q+r}2 +1} ,
$$
 with the convention that any of these binomial coefficients vanish if
its argument is negative or non-integral.  This is easily computed by
observing that $R^q(\phitwo)_*\Qdot$ is the locally constant sheaf
associated to the action of the fundamental group of $\M$ on
$\Hdot^q(\Etwotild_t,\Qdot)$, for general $t\in M$, and that the
fundamental group of $M$ is a form of $\SL(2)$.  Via this last
identification, $\Sym^r R^1\phi_*\Qdot$ is the locally constant sheaf
associated to the symmetric tensor representation of $\SL(2)$ of degree
{}~$r$.  When $r>0$ this is an irreducible representation of dimension
greater than ~1, so in particular there are no invariants or
coinvariants.  Therefore $\Hdot^p(\Mbar, j_*\Sym^rR^1\phi_*\Qdot)$
vanishes when $r>0$ and $p=0$ or~2.  Furthermore, Schur's lemma implies
that $j_*\Sym^rR^1\phi_*\Qdot$ can only be Poincar\'e dual to a Tate
twist of itself, and this completes the proof.
\Qed
\enddemo

\subhead 5.3. The second derivation of the \CKd s of $\Ebar$ and
$\Etwotild$, and computation of the cohomology of $\W1$ and~$\W2$
\endsubhead
  Using Proposition ~5.2.1 we derive the \CKd s of $\Ebar$ and
$\Etwotild$ without using the result of \cite{Scholl, 1990, Thm.1.2.1},
and determine the cohomology of $\W1$ and ~$\W2$.

\subsubhead 5.3.1. The proof for $\Ebar$ and $\W1$ \endsubsubhead
We consider first the Chow motive decomposition of $\Ebar$ given by
Theorem ~4.2, and begin by matching the cohomology groups of the
constituent motives whose cohomology we know with the constituents of
$\Hdot (\Ebar,\Qdot)$ as given in Proposition ~5.2.1.  By matching
weights also, we obtain
$$ \spreadlines{1\jot} \alignat 2
 \Hdot (\1,\Qdot) &\simeq \Hdot^0(\Mbar,\Qdot) &&\simeq
\Hdot^0(\Ebar,\Qdot)\\
 \Hdot (h^1(\Mbar),\Qdot) &\simeq \Hdot ^1(\Mbar,\Qdot) && \simeq
\Hdot^1(\Ebar,\Qdot) \\
 \Hdot (h^1(\Mbar)\tensor\LL,\Qdot) &\simeq \Hdot ^1(\Mbar,\Qdot(-1))
&& \simeq \Hdot^3(\Ebar,\Qdot) \\
  \Hdot (\LL^2,\Qdot) &\simeq \Hdot^2(\Mbar,\Qdot(-1)) &&\simeq
\Hdot^4(\Ebar,\Qdot) .
\endalignat
$$
 It therefore follows that the motives $\W1$ and $(\Ebar,\pibarinf)$
have cohomology purely of weight ~2, even if we had not already computed
that $(\Ebar,\pibarinf)$ is isomorphic to a sum of Lefschetz motives.
This already proves the existence of a \CKd\space for $\Ebar$, and that
the cohomology of the sum $\W1 \oplus (\Ebar, \pibarinf)$ must be
isomorphic to the sum $\Hdot^1(\Mbar, j_*R^1\phi_*\Qdot) \oplus
\Hdot^2\strut_{\Minf}(\Mbar, \U_\infty)$.  Then to compute the cohomology
of $\W1 = (\Ebar,\pibar1)$, and of $(\Ebar,\pibarinf)$ as a constituent
of $\Hdot (\Ebar,\Qdot)$, we observe first that
$\pibar1(\Hdot^2\strut_{\Minf}(\Mbar, \U_\infty)) =0$, since $\U_\infty$
is supported over $\Minf$ and $\pibar1$ acts as zero on all components
of $\Einf$, by Proposition ~4.4.1.  Therefore
$$
 \Hdot (\W1,\Qdot) \subseteq \Hdot^1(\Mbar, j_*R^1\phi_*\Qdot).
$$
  Conversely, it follows
from Lemma ~3.2.10 that $\Hdot ((\Ebar,\pibarinf),\Qdot)$ is generated by
the classes of some $\theta_c(m)$ (modulo homological equivalence), and
thus consists entirely of algebraic cohomology classes in
$\Hdot^2(\Ebar,\Qdot(1))$.  On the other hand, by virtue of the
$\Gal(K^{\text{sep}}/K)$-module structure of
$H_{\text{\'et}}^1(\Mbar \tensor
K^{\text{alg}},j_*R^1\phi_*\QQ_\ell)$
\cite{Deligne, 1969}, or the Hodge structure of
$H_B^1(\Mbar(\CC)^{\text{an}}, j_*R^1\phi_*\QQ)$ \cite{Shioda, 1972}
\cite{\v Sokurov, 1976 and 1981} \cite{Zucker, 1979}, $\Hdot^1(\Mbar,
j_*R^1\phi_*\Qdot)$ cannot contain any algebraic cohomology classes.
Therefore the only possibility is that $\Hdot (\W1,\Qdot) \simeq
\Hdot^1(\Mbar, j_*R^1\phi_*\Qdot)$ and $\Hdot ((\Ebar,\pibarinf),\Qdot)
\simeq \Hdot^2\strut_{\Minf}(\Mbar, \U_\infty)$.

\subsubhead 5.3.2. The proof for $\Etwotild$ and $\W2$ \endsubsubhead
 The argument computing the \CKd\space of $\Etwotild$ and the cohomology
of $\W2$ follows similar lines.  From the Chow motive computations in
section four, using known cohomology groups and matching weights we get
$$ \spreadlines{1\jot} \alignat2
 \Hdot (\1,\Qdot) &\simeq \Hdot^0(\Mbar,\Qdot) &&\simeq
\Hdot^0(\Etwotild,\Qdot) \\
 \Hdot(h^1(\Mbar),\Qdot) &\simeq \Hdot^1(\Mbar,\Qdot) &&\simeq
\Hdot^1(\Etwotild,\Qdot) \\
 \Hdot(h^2(\Mbar),\Qdot) &\simeq \Hdot^2(\Mbar,\Qdot) &&\subset
\Hdot^2(\Etwotild,\Qdot) \\
 \Hdot(2(\W1),\Qdot) &\simeq \Hdot^1(\Mbar, j_*R^1\phi_*\Qdot)^{\oplus 2}
&&\subset \Hdot^2(\Etwotild,\Qdot) \\
 \Hdot(3(h^0(\Mbar)\tensor\LL),\Qdot) &\simeq
\Hdot^0(\Mbar,\Qdot(-1))^{\oplus3} &&\subset \Hdot^2(\Etwotild,\Qdot) \\
 \Hdot(3(h^1(\Mbar)\tensor\LL),\Qdot) &\simeq
\Hdot^1(\Mbar,\Qdot(-1))^{\oplus3}
&&\subset \Hdot^3(\Etwotild,\Qdot) \\
 \Hdot(3(h^2(\Mbar)\tensor\LL),\Qdot) &\simeq
\Hdot^2(\Mbar,\Qdot(-1))^{\oplus3}
&&\subset \Hdot^4(\Etwotild,\Qdot) \\
 \Hdot(2(\W1\tensor\LL),\Qdot) &\simeq \Hdot^1(\Mbar,
j_*R^1\phi_*\Qdot(-1))^{\oplus2} &&\subset \Hdot^4(\Etwotild,\Qdot) \\
  \Hdot(h^0(\Mbar)\tensor\LL^2,\Qdot) &\simeq  \Hdot^0(\Mbar,\Qdot(-2))
&&\subset \Hdot^4(\Etwotild,\Qdot) \\
 \Hdot(h^1(\Mbar)\tensor\LL^2,\Qdot) &\simeq  \Hdot^1(\Mbar,\Qdot(-2))
&&\simeq
\Hdot^5(\Etwotild,\Qdot) \\
 \Hdot(h^2(\Mbar)\tensor\LL^2,\Qdot) &\simeq  \Hdot^2(\Mbar,\Qdot(-2))
&&\simeq
\Hdot^6(\Etwotild,\Qdot) .
\endalignat
$$
 Therefore the cohomology of the sum of motives $\W2 \oplus
(\Etwotild, \pitildinf\sparen2(\Etwotild/\Mbar)) \oplus
(\Etwotild, \pitildinf\sparen4(\Etwotild/\Mbar))$ is the sum of the
cohomology groups $\Hdot^1(\Mbar, j_*\Sym^2\phi_*\Qdot) \oplus
\Hdot^{2}\strut_{\Minf}(\Mbar, \U_\infty\sparen{2}) \oplus
\Hdot^{4}\strut_{\Minf}(\Mbar, \U_\infty\sparen{4})$.  Then by
Proposition ~4.4.1
$\sym2\pitild{1,1}(\Etwotild/\Mbar)$ annihilates (the classes of) the
components of $\Etwotildinf$, which means that $\Hdot (\W2,\Qdot)$ is
disjoint from $\Hdot^{2}\strut_{\Minf}(\Mbar, \U_\infty\sparen{2})$.
Therefore we have $\Hdot ((\Etwotild,\pitildinf\sparen
2(\Etwotild/\Mbar)),\Qdot) \simeq \Hdot^{2}\strut_{\Minf}(\Mbar,
\U_\infty\sparen{2})$.  But we also know not only that
$\Hdot^{2}\strut_{\Minf}(\Mbar, \U_\infty\sparen{2})$ pairs
nondegenerately with $\Hdot^{4}\strut_{\Minf}(\Mbar,
\U_\infty\sparen{4})$, but also that the Chow groups of
$(\Etwotild,\pitildinf\sparen 2(\Etwotild/\Mbar))$ and
$(\Etwotild,\pitildinf\sparen 4(\Etwotild/\Mbar))$ pair nondegenerately,
by step four in the proof of Proposition ~4.6.1, so we must also have
that $\Hdot ((\Etwotild,\pitildinf\sparen 4(\Etwotild/\Mbar)),\Qdot)
\simeq \Hdot^{4}\strut_{\Minf}(\Mbar, \U_\infty\sparen{4})$.  Therefore
the only remaining possibility is that
$$
 \Hdot (\W2,\Qdot) \simeq \Hdot^1(\Mbar, j_*\Sym^2\phi_*\Qdot),
$$
 as claimed, from which it also follows that $\Etwotild$ has a \CKd.
\Qed

\head \bf 6. The filtration on the Chow groups of $\Ebar$ and
$\Etwotild$ \endhead
\nobreak\noindent
 Recall that  Conjecture~A predicts the
existence of a Chow-K\"unneth decomposition; for $\Ebar$ and
$\Etwotild$ this is proved in Theorem~5.1 (see also Theorem~4.2).  In
this section we start with those \CKd s, and then for $\Ebar$ and
$\Etwotild$ we prove Conjectures~B, that $\Chow j {h^i(X)} =0$ for $i<j$
or $i>2j$, and~D, that $F^1\Chow j X =\Chowhom j X$, and a large part of
Conjecture~C, that the filtration on the Chow groups induced by these
\CKd s is the natural one.  Although the conjectures have been proved
for surfaces in general \cite{Murre, 1990}, here we give a different
proof for ~$\Ebar$, using the extra structure that Theorems~5.1 and ~4.2
reveal.  In particular, we find the Chow groups of $\W1$ (see Theorem~6.2
below), which we then use in the proof of Conjectures~B and ~D for the
threefold ~$\Etwotild$.  As for proving Conjecture~C for $\Etwotild$,
precise statements are given in Theorem~6.2 below, but our results may be
summarized by observing first that it is trivially true for $\Chow
0{\Etwotild}$, and it is equivalent to Conjecture~D, which we prove, for
$\Chow 1{\Etwotild}$; but then we also prove that $F^1\Chow j {\Etwotild}
= \Chowhom j{\Etwotild}$ for $j=2,3$, and that $F^2\Chow3{\Etwotild}
=\Chowalb3{\Etwotild}$.  So what's missing is $F^2\Chow 2{\Etwotild}$,
which is contained in the kernel of an Abel-Jacobi map defined on
$\Chowhom2{\W2}$ (Proposition~6.5.6), and $F^3\Chow 3{\Etwotild}$,
which we show equals $\Chow 3{\W2}$ (Proposition ~6.6.1).

\subhead 6.1. Notation \endsubhead
 With the present state of knowledge about Chow groups we can at best
prove the naturality of a step in the filtration when there is a
clear, geometrically described candidate for it.  If there are such
natural
candidates and if the filtration is this natural one, then by abuse of
language
we will say that Conjecture~C is true.  For a
smooth projective variety ~$X$ over a field ~$k$, we have
$$
 \Chowhom j X := \Ker\{\gamma:\Chow j X \to \Hdot^{2j}(X,\Qdot(j))\},
$$
 where $\gamma$ is the cycle class map.  Further, let
$$
\Chowalb d X := \Ker\{ \Alb: \Chowhom d X \to
\Alb(X)\tensor\QQ\},
$$
 where $\Alb(X)$ is the Albanese of $X$ and $d=\dim X$.
 Finally, supposing for simplicity that $\operatorname{char}k= 0$, let
$$
 \Chowaj j X := \Ker\{\AJ : \Chowhom j X \to J^j(X)\tensor \QQ\},
$$
 where $\AJ$ is the Abel-Jacobi map to the $j$@-th intermediate Jacobian
{}~$J^j(X)$.

\proclaim{Theorem 6.2}
\nobreak\par\nobreak\noindent {\rm (1)} For the \CKd\space of $\Ebar$
described in Theorem ~5.1(1) we have
\roster
\item"(i)" $\Chow j{h^i(\Ebar)} =0$ for $i<j$ or $i>2j$, i.e.,
Conjecture~B is true for ~$\Ebar$;
\smallskip
\item"(ii)" $F^1\Chow j {\Ebar} = \Chowhom j{\Ebar}$ for $1\le j\le 2$,
i.e., Conjecture~D is true for ~$\Ebar$;
\smallskip
\item"(iii)" $F^2\Chow 2{\Ebar} = \Chowalb 2{\Ebar}$, and therefore the
filtration is independent of the choice of Chow-K\"unneth projectors
$\pi_i(\Ebar)$, i.e., Conjecture ~C is true for ~$\Ebar$.
\endroster
 In particular, for the Chow groups of $\W1$ we have
$$\spreadlines{1\jot} \gather
 \Chow 0{\W1} = \Chow 1{\W1} =0 \\
 \Chow 2{\W1} = F^2\Chow 2{\W1} = \Chowalb 2{\Ebar} .
\endgather
$$

\noindent {\rm (2)} For the \CKd\space of $\Etwotild$ described in
Theorem ~5.1(2) we have
\roster
\item"(i)" $\Chow j{h^i(\Etwotild)} =0$ for $i<j$ or $i>2j$, i.e.,
Conjecture~B is true for ~$\Etwotild$;
\smallskip
\item"(ii)" $F^1\Chow j {\Etwotild} = \Chowhom j{\Etwotild}$ for $1\le
j\le 3$, i.e., Conjecture~D is true for ~$\Etwotild$.
\smallskip
\item"(iii)" Towards Conjecture ~C we also have
\smallskip
\itemitem{\rm(a)}
$F^2\Chow 2{\Etwotild} \subseteq \Chowaj2{\Etwotild}$, when
$\operatorname{char}K= 0$.
\smallskip
\itemitem{\rm(b)} $F^2\Chow3{\Etwotild} = \Chowalb3{\Etwotild}$.
\endroster
 In particular, for the Chow groups of $\W2$ we have
$$\spreadlines{1\jot} \gather
 \Chow 0{\W2} = \Chow 1{\W2} =0 \\
 \Chow 2{\W2} = F^1\Chow 2{\W2} = \Chowhom 2{\W2} \\
 \Chow 3{\W2} = F^3\Chow 3{\W2} = F^3\Chow 3{\Etwotild} .
\endgather
$$
\endproclaim

\subhead 6.3. Preliminaries to the proof of Theorem~6.2 \endsubhead
 Before getting into the proof of Theorem~6.2, we begin with some
elementary
but useful observations.

\subsubhead 6.3.1. The conjectures for $\Chow 0 X$ \endsubsubhead
 For any smooth projective $X$ with a Chow-K\"unneth decomposition,
$\Chow 0 X$ trivially satisfies Conjectures~B, C and~D.  For $\Chow 0 X
\tensor\Qdot = \Hdot^0(X,\Qdot)$, from which it follows that $\pi_0(X)$
is the identity on $\Chow 0 X$.  Then by orthogonality, $\pi_i(X)(\Chow 0
X) = 0$ for $i>0$.

\subsubhead 6.3.2. The Chow groups of a motive \endsubsubhead
 Recall from the definitions in section~1 that for any Chow motive $M_0$
we have
$$
 \Chow j {M_0\tensor \LL^m} = \Chow{j-m} {M_0} .
$$

\subsubhead 6.3.3. The Chow groups of $\spec K$ and $h^1(\Mbar)$
\endsubsubhead
 Two special cases of 6.3.2 which we will use in the proof of Theorem~6.2
are
$$
\CH^j(\LL^m,\QQ) \simeq \cases \QQ , & \text{if } j=m, \\
    0 , &\text{otherwise}
 \endcases
$$
 and
$$
\CH^j(h^1(\Mbar)\tensor\LL^m),\QQ) \simeq \cases
     \Jac(\Mbar)\tensor\QQ, &\text{if }  j=m+1, \\
    0 , &\text{otherwise} .
\endcases
$$
 Moreover, these motives satisfy Conjectures ~B and ~D, in an obvious
sense (cf.\space Proposition ~1.2.5), and they satisfy Conjecture~C in
the sense that the filtrations on their Chow groups are the natural ones.

\subsubhead 6.3.4. The motives of $\W1$ and $\W2$ \endsubsubhead
 Given the \CKd s in Theorem~5.1, the previous paragraph together with
Lemma~1.2.5 imply that to prove Conjectures~B and ~D for $\Ebar$ and
$\Etwotild$ it would suffice to prove them for $\W1$ and $\W2$ (with the
obvious understanding of what that means).  However, as in section~5.3,
the reality is that it works the other way around:  Anything nontrivial
that we are able to say about the Chow groups of $\W1$ and $\W2$ comes
indirectly, via analyzing the Chow groups of $\Ebar$ and ~$\Etwotild$.
As Chow motives, we have
$$
 \W1 \simeq h^2(\W1) \qquad  \W1 \tensor \LL \simeq h^4(\W1\tensor\LL),
\qquad  \W2 \simeq h^3(\W2) ,
$$
 since these motives have cohomology only in these degrees, see ~5.1
and ~5.3.  Then by applying 6.3.1 for $\Ebar$ and$\Etwotild$ we find
that
$$
 \Chow 0{\W1} = \Chow 0{\W2} =0 .
$$

\subsubhead 6.3.5. Organization of the proof \endsubsubhead
 The rest of this section is devoted to the proof of Theorem~6.2.  In the
next subsection we consider $\Chow 1{\Ebar}$ and $\Chow 1 {\Etwotild}$,
in the following, $\Chow 2{\Ebar}$ and $\Chow 2 {\Etwotild}$, and in the
last, $\Chow 3{\Etwotild}$.

\subhead 6.4.  Analysis of $\Chow 1{\Ebar}$ and $\Chow 1{\Etwotild}$
\endsubhead
 As the proof of Conjectures~B, C and ~D is the same for both $\Chow
1{\Ebar}$ and $\Chow 1{\Etwotild}$, the details are written out only for
{}~$\Etwotild$.   The proof is based on two lemmas, the first of which
describes a general approach to verifying the conjectures for $\Chow
1{M_0}$ for any Chow motive ~$M_0$, while the second identifies the
Picard (as well as the Albanese) variety of an elliptic modular variety
with the Jacobian of the elliptic modular curve over which it lies.

\proclaim{Lemma 6.4.1}
 Let $M_0$ be a Chow motive in $\scM(k)$, and assume $M_0$ has a
Chow-K\"unneth decomposition.  Suppose that $\pi_1(M_0)(\xi) =
\xi$ for all $\xi\in \CH^1_{\text{hom}}(M_0,\QQ)$.  Then for $\xi \in
\Chow1{M_0}$,
\roster
\smallskip
\item $\xi=\pi_1(M_0)(\xi) +\pi_2(M_0)(\xi)$.
\smallskip
\item $\pi_i(M_0)(\xi) =0$ for $i\ne 1,2$.
\smallskip
\item $\Ker\pi_2(M_0) = \CH^1_{\text{hom}}(M_0,\QQ)$.
\endroster
\endproclaim

\demo{Proof}
 To begin with,
$$
\xi-\pi_2(M_0)(\xi) \in \Ker \pi_2(M_0) \subseteq
\CH^1_{\text{hom}}(M_0,\QQ),
$$
  see ~1.2.3.  Then applying the hypothesis,
$$
\xi-\pi_2(M_0)(\xi) = \pi_1(M_0)(\xi-\pi_2(M_0)(\xi)) = \pi_1(M_0)(\xi) ,
$$
 where the second equality follows from the orthogonality of
$\pi_1(M_0)$ and $\pi_i(M_0)$. This proves part \therosteritem1, and part
\therosteritem2 follows from the mutual orthogonality of all the
Chow-K\"unneth projectors.  To prove part \therosteritem3, if
$\xi \in \CH^1_{\text{hom}}(M_0,\QQ)$, then $\xi=\pi_1(M_0)(\xi)$ by
assumption, and therefore $\pi_2(M_0)(\xi) =0$, once more by
orthogonality.
\Qed
\enddemo

\medpagebreak
 A special case of the following lemma already occurs in \cite{Shioda,
1972, p.24}.

\proclaim{Lemma 6.4.2}
\roster
\item $\Pic^0(\Ebar) \simeq \Pic^0(\Etwotild) \simeq \Pic^0(\Mbar) =
\Jac(\Mbar)$
\item $\Alb(\Ebar) \simeq \Alb(\Etwotild) \simeq \Alb(\Mbar) =
\Jac(\Mbar)$
\endroster
\endproclaim

\demo{Proof}
 Consider for instance $\Etwotild$.  Letting $\alphatild(\0) : \Mbar\to
\Etwotild$ denote the extended identity section, then
$$
 \alphatild(\0)^* \circ (\phitwotild)^* : \Jac(\Mbar) \to
\Pic^0(\Etwotild) \to \Jac(\Mbar)
$$
 is the identity map.  Then \therosteritem1 follows from the fact that
$\dim
\Hdot^1(\Etwotild, \Qdot ) = \dim \Hdot^1(\Mbar, \Qdot )$, see Proposition
{}~5.2.1.  By duality, \therosteritem2 follows as well.  The argument is
the same for $\Ebar$.
\Qed
\enddemo

\proclaim{Proposition 6.4.3}
 Conjectures ~B, C and ~D are true for $\Chow 1{\Ebar}$ and
\break
 $\Chow 1{\Etwotild}$.
\endproclaim

\demo{Proof}
  Consider for instance $\Etwotild$.  From Lemma~6.4.2 we get the
following commutative diagram.
$$\CD
 \CH^1_{\text{hom}}(\Etwotild,\QQ) @> \pi_1(\Etwotild) >>
\CH^1_{\text{hom}}(\Etwotild,\QQ) \\  @|  @| \\
\Pic^0(\Etwotild)\tensor\QQ @> \pi_1(\Etwotild) >>
\Pic^0(\Etwotild)\tensor\QQ \\
 @V (\phitwotild)^* V \sim V  @V (\phitwotild)^* V \sim V \\
\Jac(\Mbar)\tensor\QQ  @>\sim > \pi_1(\Mbar) > \Jac(\Mbar)\tensor\QQ
\endCD
$$
 Therefore we can apply Lemma~6.4.1, with $M_0 =h(\Etwotild)$, and the
conclusions of that lemma give the  Conjectures ~B, C and ~D for
$\CH^1(\Etwotild)$.  The same argument works for ~$\Ebar$.
\Qed
\enddemo

\proclaim{Corollary 6.4.4}
\roster
\item $\Chow 1{\W1} =0$.
\smallskip
\item $\Chow 1{\W2} =0$.
\endroster
\endproclaim

\demo{Proof}
 Consider the cycle map $\gamma: \Chow 1{\W1} \to \Hdot^2(\Ebar,
\Qdot)$.  Since from ~5.3.1 we know that $\W1$ has no algebraic
cohomology, $\Chow 1{\W1} \subset \Chowhom 1{\Ebar}$.  But
$\pi_2(\Ebar)$, which acts as the identity on $\W1$, also acts as zero on
$\Chowhom 1{\Ebar}$.  Hence $\Chow 1 {\W1}=0$.

 Similarly but even easier, $\Chow 1{\W2} =0$ because the identity of
$\W2$ is a part of $\pi_3(\Etwotild)$, which acts as zero on $\Chow 1
{\Etwotild}$.
\Qed
\enddemo

\subhead 6.5.  Analysis of $\Chow 2{\Ebar}$ and $\Chow 2{\Etwotild}$
\endsubhead
 \nobreak\par\nobreak
\subsubhead 6.5.1. The Albanese kernel \endsubsubhead
   For any smooth projective variety $X$ of dimensional~$d$ the
Chow-K\"unneth projector $\pi_{2d-1}(X)$ acts as the identity on the
Albanese variety ~$\Alb(X)$ \cite{Murre, 1990} (this is proved by looking
at the torsion points).  Hence from the commutative diagram
$$\CD
\CH^d_{\text{hom}}(X,\QQ) @> \pi_{2d-1} >> \CH^d_{\text{hom}}(X,\QQ) \\
@V \Alb VV  @VV \Alb V \\
\Alb(X)\tensor\QQ @> \sim > \pi_{2d-1} > \Alb(X) \tensor \QQ
\endCD
$$
 it follows that $\Ker(\pi_{2d-1}(X))\subseteq \Chowalb d {X}$; this may
be compared with ~1.2.3.

\proclaim{Proposition 6.5.2}
 Conjectures ~B, C and ~D are true for $\Chow 2{\Ebar}$.  Moreover
$$
 \Chow{}{\W1} = \Chow 2{\W1} = \Chowalb 2{\Ebar} .
$$
\endproclaim

\demo{Proof}
 Consider the \CKd\space of $\Ebar$ in Theorem~5.1(1) (see also
Theorem~4.2):  Other than $\W1$, all the submotives of $h(\Ebar)$ are of
the form $\LL^m$ or $h^1(\Mbar)\tensor\LL^m$, and thus satisfy the
conjectures, as in ~6.3.3.  Thus, by Lemma~1.2.5, to prove Conjecture~B
for $\Chow 2{\Ebar}$ it suffices to verify that $\pi_0(\Ebar)$ and
$\pi_1(\Ebar)$ act as zero on $\Chow 2{\W1}$.  But this is immediate,
since by
our construction of the Chow-K\"unneth decomposition in Theorem~5.1 we
have that
$\id_{(\W1)}$ is orthogonal to $\pi_0(\Ebar)$ and $\pi_1(\Ebar)$.  Thus
Conjecture~B follows.  Moreover, as $\id_{(\W1)}$ is part of
$\pi_2(\Ebar)$ we
have
$$\spreadlines{1\jot} \alignat3
 \Chow 2{\Ebar} &\simeq \Chow 2{h^2(\Ebar)} && \oplus \Chow2{h^3(\Ebar)}
&& \oplus \Chow2{h^4(\Ebar)} \\
 &\simeq \Chow2{\W1} &&\oplus \Jac(\Mbar)\tensor\QQ && \oplus \QQ .
\endalignat
$$
 To prove Conjecture~D, observe that an element $\alpha \in\Chow2{\Ebar}$
is contained in $\Chowhom2{\Ebar}$ if and only if the cycle class map
acts on the component of $\alpha$ in $\Chow2{h^4(\Ebar)}$ as zero.  But
on $\Chow 2{h^4(\Ebar)}$ the cycle class map is the degree map, and thus
an isomorphism with ~$\QQ$.  Therefore $\Ker(\pi_4(\Ebar)) = \Chowhom
2{\Ebar}$, i.e., Conjecture~D is true.  To prove Conjecture~C, from 6.5.1
we know that $\Ker(\pi_3(\Ebar))
\subseteq \Chowalb 2{\Ebar}$.  Then to see that this inclusion is
actually an equality, the essential fact, from Theorem~5.1, is that
$h^3(\Ebar) \simeq h^1(\Mbar) \tensor\LL$, for together with 6.3.3 and
6.4.2 this implies that
$$
\Chowhom 2{\Ebar}/\Ker(\pi_3(\Ebar)) \simeq \Jac(\Mbar)\tensor\QQ \simeq
\Alb(\Ebar)\tensor\QQ .
$$
 Thus $F^2\Chow 2{\Ebar}= \Chowalb 2{\Ebar}$, i.e., Conjecture~C is
true.  Finally, it now follows directly from Theorem~5.1 that $\Chowalb
2{\Ebar} = \Chow 2{\W1}$, and this is the entire Chow ring $\CH(\W1,\QQ)$
by 6.3.1 and Corollary~6.4.4.
\Qed
\enddemo

 This finishes the proof of part~(1) of Theorem~6.2.

\proclaim{Proposition 6.5.3}
 Conjectures ~B and ~D are true for $\Chow 2{\Etwotild}$.
\endproclaim

\demo{Proof}
 Consider the \CKd\space of $\Etwotild$ described in Theorem~5.1(2) (see
also Theorem~4.2):  By the previous proposition and 6.3.3, all the
submotives that occur except possibly $\W1\tensor\LL$ or $\W2$ satisfy
the conjectures.  Thus, applying Lemma ~1.2.5, to prove Conjecture~B for
$\Chow2{\Etwotild}$ it suffices to check that $\pi_i(\Etwotild)$ acts as
zero on both $\Chow2{\W1\tensor\LL}$ and $\Chow 2{\W2}$, for $i<2$ or
$i>4$.
But this is true, since by our Chow-K\"unneth decomposition
$\id_{(\W1\tensor\LL)}$ and $\id_{(\W2)}$ are both orthogonal to these
{}~$\pi_i(\Etwotild)$.  Moreover, $\id_{(\W1\tensor\LL)}$ is part of
$\pi_4(\Etwotild)$ and $\id_{(\W2)}$ is part of ~$\pi_3(\Etwotild)$.

  To prove Conjecture~D we must show that equality holds in the inclusion
$\Ker(\pi_4(\Etwotild))\subseteq \Chowhom 2{\Etwotild}$; it suffices to
see that the cycle class map ~$\gamma$ is injective on $\Chow
2{h^4(\Etwotild)}$.  But from Theorem~5.1 we know that $h^4(\Etwotild) =
2(\W1\tensor\LL) \oplus n\LL^2$.  Then from 6.3.2 and Corollary~6.4.4
we find that $\Chow2{\W1\tensor\LL}= \Chow1 {\W1} =0$, whereas from
{}~6.3.3 and the definitions we get that $\Chow2 {\LL^2} = \Chow0{\spec
K}$, on which $\gamma$ is injective.  Conjecture~D follows.
\Qed
\enddemo

\subsubhead 6.5.4. The Abel-Jacobi kernel \endsubsubhead
 Let $X$ be a smooth projective threefold over a field ~$k$, and assume
for simplicity $\operatorname{char} k =0$.  Then when $X$ has a
Chow-K\"unneth decomposition that satisfies Conjectures~B and ~D, there
is a commutative diagram
$$\CD
\Chowhom 2 X @> \pi_3(X) >> \Chowhom 2 X \\ @V \AJ_X VV  @VV \AJ_X V \\
J^2(X)\tensor \QQ @> \sim > \pi_3(X) > J^2(X)\tensor\QQ
\endCD
$$
 where $J^2(X)$ is the intermediate Jacobian; the lower homomorphism is
an isomorphism because algebraic correspondences respect Hodge structure
and $\pi_3(X)$ is an isomorphism on $H_B^3(X,\QQ_B)$ (which is the
starting point for the construction of $J^2(X)$).  From the diagram it
follows that
$$
\Ker(\pi_3(X)) \subseteq \Ker(\AJ_X) ,
$$
 or, equivalently,
$$
 F^2\Chow2 X \subseteq \Chowaj 2 X ;
$$
 this may be compared with 1.2.3 and 6.5.1.

\proclaim{Conjecture 6.5.5}
 When $X$ is a smooth projective threefold over a field $k$ of
characteristic zero, and there exists a Chow-K\"unneth decomposition for
$X$ such that $\Chow j{h^i(X)} =0$ for $i<j$ or $i>2j$, and $F^1\Chow2 X
=\Chowhom 2 X$, then $F^2\Chow 2 X = \CH^2_{\AJ}(X,\QQ)$ (or
equivalently, $\Ker(\pi_3(X)) = \Ker(\AJ_X)$).
\endproclaim

\proclaim{Proposition 6.5.6}
\roster
\item  $\Chow 2{\W2} = \Chowhom 2{\W2}$.
\smallskip
\item  Assume $\operatorname{char} K = 0$. Then there is a map
$$
 \AJ_{(\W2)} : \Chow2{\W2} \lra J^2(\Etwotild)\tensor\QQ
$$
 and $F^2\Chow 2{\Etwotild} = \Chowaj 2{\Etwotild}$ if and only if
$\AJ_{(\W2)}$ is injective.
\endroster
\endproclaim

\demo{Proof}
 The first statement follows directly from the definitions and the fact
that $\id_{(\W2)} = \pi_3(\W2)$, as observed in ~6.3.4.  Then the
existence of $\AJ_{(\W2)}$ comes by composing $\AJ_{(\Etwotild)}$ with
$\sym2\pitild{1,1}=\id_{(\W2)}$.  To prove the last statement of
part~\therosteritem2, we first note that $h^3(\Etwotild) = \W2 \oplus
3(h^1(\Mbar)\tensor\LL)$, by Theorem~5.2. and next that
$\AJ_{(\Etwotild)}$ is
injective on the summand
$$
 {\Chow 2{h^1(\Mbar)\tensor \LL}}^{\oplus3} \simeq {\Chow
1{h^1(\Mbar)}}^{\oplus3} \simeq (\Jac(\Mbar)\tensor\QQ)^{\oplus3},
$$
 since it coincides with (three copies of) the usual map from divisors on
a
curve to the Jacobian.
\Qed
\enddemo

\subhead 6.6.  Analysis of $\Chow3{\Etwotild}$ \endsubhead
\nobreak\par\nobreak
\proclaim{Proposition 6.6.1}
\roster
\item  Conjectures ~B and ~D are true for $\Chow 3{\Etwotild}$.
\smallskip
\item  $F^2\Chow 3{\Etwotild} = \Chowalb 3{\Etwotild}$.
\smallskip
\item  $\Chow 3{\W2} = F^3\Chow 3{\Etwotild}$.
\endroster
\endproclaim

\demo{Proof}
 For part~\therosteritem1, consider the \CKd\space of $\Etwotild$ in
Theorem~5.1(2) (see also Theorem~4.2):  In view of 6.3.3 and
Proposition~1.2.5, to prove Conjecture~B for $\Chow 3{\Etwotild}$, we
need only verify it for $\W1$, $\W1\tensor\LL$ and $\W2$.  But $\Chow 3
{\W1} =0$.  Thus for Conjecture~B to be true we must have
$\pi_i(\W1\tensor\LL)(\Chow3{\W1\tensor\LL})=0$ for $0\le i\le 2$, which
is the case since $\id_{(\W1\tensor\LL)}$ is orthogonal to $\pitild i
(\Etwotild)$ for $i<3$ and moreover is part of $\pi_4(\Etwotild)$.  We
must also have that $\pi_i(\W2)(\Chow3{\W2})=0$ for $0\le i\le 2$, which
is the case since $\id_{(\W2)}$ is orthogonal to $\pitild i (\Etwotild)$
for
$i<3$ and moreover is part of $\pi_3(\Etwotild)$.  Conjecture~D follows
for
$\Chow 3{\Etwotild}$ similarly as for $\Chow 2{\Ebar}$:
$\Ker(\pi_6(\Etwotild)) \subseteq \Chowhom 3{\Etwotild}$, and
$\Chow 3{h^6(\Etwotild)} =\CH^3(\LL^3,\QQ) \simto \QQ$ is the degree map.

To prove \therosteritem2, first observe that by 6.5.1
$\Ker(\pi_5(\Etwotild)) \subseteq \Chowalb3{\Etwotild}$.  But here we
have equality because of the commutative diagram
$$
\CD
\CH^3_{\Alb}(h^5(\Etwotild),\QQ) @> \sim >>
\CH^3(h^1(\Mbar)\tensor\LL^2,\QQ)\simeq \CH^1(h^1(\Mbar),\QQ) \\ @V \Alb
VV  @VV \Alb V  \\
\Alb(\Etwotild)\tensor\QQ @> \sim >> J(\Mbar)\tensor\QQ
\endCD
$$
 where the the top row is isomorphism from Theorem~5.1 and the bottom row
is an isomorphism by Lemma~6.4.2(2).

Part~\therosteritem3 follows from observing that
$$ F^3\Chow 3{\Etwotild} = \Chow 3{h^3(\Etwotild)} = \Chow
3{3(h^1(\Mbar)\oplus\W2},
$$
 and
$$
\CH^3(h^1(\Mbar)\tensor\LL,\QQ) = \CH^2(h^1(\Mbar),\QQ) =0 .
$$
\Qed
\enddemo

\subsubhead Remark 6.6.2 \endsubsubhead
 We remark that
$$
 F^2\Chow3{\Etwotild} \big/ F^3\Chow3{\Etwotild} = \Chow3{h^4(\Etwotild)}
\simeq {\Chowalb 2{\Ebar}}^{\oplus2} .
$$
 For by Theorem~5.2 (see also Theorem~4.2)
$$
 h^4(\Etwotild) \simeq 2 h^4(\W1\tensor\LL) \oplus h^4(n\LL^2) .
$$
 Then by 6.3.2,
$$\spreadlines{1\jot} \align
 \Chow3{h^4(\Etwotild)} &\simeq {\Chow 3{h^4(\W1\tensor\LL)}}^{\oplus2}
\oplus \Chow3{h^4(\LL^2)}^{\oplus n} \\
  & \simeq \CH^2(\W1,\QQ)^{\oplus2} \\
  & \simeq {\Chowalb 2{\Ebar}}^{\oplus2}
\endalign
$$
 by Proposition~6.5.2.


\bigpagebreak


\Refs
\refindentwd=3.0pc 
\parskip=\smallskipamount

\ref \by del Angel, P.L. and M\"uller-Stach, S.
\paper Motives of uniruled 3-folds
\paperinfo preprint
\yr 1996
\endref

\ref
\by  Ash, A., Mumford,  D., Rapoport, M., Y. Tai, Y.
\book  Smooth Compactifications of Locally Symmetric Spaces
\publ  Math Sci Press
\publaddr  Brookline
\yr  1975
\endref

\ref \by del Ba\~no Rolla
\paper Localisation of $\Gm$-actions and $\lambda$-structure in the
   category of Chow motives
\paperinfo preprint
\yr 1995
\endref

\ref
\by  Beilinson,~A.A.
\paper  Height pairings between algebraic cycles
\inbook $K$-Theory and Arithmetic
\bookinfo Lect. Notes Math. 1289
\pages 1--26
\publaddr New York, Berlin
\publ Springer Verlag
\yr  1987
\moreref 
\inbook Current Trends in Arithmetical Algebraic Geometry
\bookinfo Contemp. Math. {\bf 67}
\pages 1--26
\publ  Amer. Math. Soc.
\publaddr Providence
\yr  1987
\endref

\ref  \bysame, 
 Bernstein,~J.N. and Deligne,~P.
\paper  Faisceaux pervers
\paperinfo Analyse et topologie sur les espaces singulier, I
\jour  Ast\'erisque
\vol  100
\yr  1983
\endref

\ref \by Bloch, S.
\book Lectures on Algebraic Cycles
\bookinfo Duke Univ. Math. Series IV
\publ Duke Univ. Press
\publaddr Durham
\yr 1980
\endref

\ref
\by Deligne, P.
\paper Formes modulaires et repr\'esentations $\ell$-adiques
\inbook  S\'eminaire Bourbaki
\bookinfo Lecture Notes in Mathematics 179
\publ  Springer-Verlag
\publaddr  New York
\yr  1969
\endref

\ref \bysame 
\paper  Hodge cycles on abelian varieties (Notes by J.S.~Milne)
\inbook Hodge Cycles, Motives, and Shimura Varieties
\bookinfo Lecture Notes in Mathematics 900
\publ  Springer-Verlag
\publaddr  New York
\yr  1982
\pages  9--100
\endref

\ref \bysame 
  and Rapoport, M.
\paper  Les sch\'emas de modules de courbes elliptiques
\inbook Modular Functions of One Variable II
\bookinfo  Lect. Notes. Math. 349
\eds P. Deligne, W. Kuyk
\publ Springer
\publaddr Berlin etc.
\pages 143--316
\yr 1973
\endref


\ref \by Deninger, C. and Murre, J.P.
\paper Motivic decomposition of abelian schemes and the Fourier transform
\jour J. reine angew. Math.
\vol 422
\yr 1991
\pages 201--219
\endref

\ref \by Fulton, W.
\book Intersection Theory
\bookinfo Ergeb. Math. Grenzgeb. (3) {\bf 2}
\publ Springer
\publaddr Berlin
\yr 1984
\endref

\ref \by Gordon, B.B.
\paper Algebraic cycles and the Hodge structure of a Kuga fiber variety
\jour Trans. Amer. Math. Soc.
\vol 336
\yr 1993
\pages 933--947
\endref

\ref \by  Hartshorne, R.
\book  Algebraic Geometry
\publ  Springer-Verlag
\publaddr  New York
\yr  1977
\endref

\ref \by Jannsen,~U.
\paper Motives, numerical equivalence, and semi-simplicity
\jour Invent. Math.
\vol 107
\yr 1992
\pages 447--452
\endref

\ref \bysame 
\paper Motivic sheaves and filtrations on Chow groups
\inbook Motives
\bookinfo Proc. Symp. Pure Math. 55, vol 1
\eds Jannsen, Kleiman, Serre
\publ American Mathematical Society
\publaddr Providence
\yr 1994
\pages 245--302
\endref


\ref \by Katz, N.M. and Mazur, B.
\book Arithmetic Moduli of Elliptic Curves
\bookinfo Ann. Math. Studies 108
\publ Princeton U. Press
\publaddr Princeton
\yr 1985
\endref 

\ref \by Kleiman, S.L.
\paper Algebraic cycles and the Weil conjectures
\inbook  Dix Exposes sur la Cohomologie des Schemas
\publ  North-Hollandhoff
\publaddr  Amsterdam
\yr  1968
\pages 359--386
\endref

\ref \bysame 
\paper Motives
\inbook  Geometry, Oslo 1970
\bookinfo Proc. Fifth Nordic Summer School in Math, Oslo 1970
\publ  Walters-Noordhoff
\publaddr  Groningen
\yr  1972
\pages 53--82
\endref


\ref \bysame 
\paper  The standard conjectures
\inbook Motives
\bookinfo Proc. Symp. Pure Math. 55, part~1
\eds U.~Jannsen, S.~Kleiman, J.-P.~Serre
\publ Amer. Math. Soc.
\publaddr Providence
\pages 3--20
\yr 1994
\endref

\ref \by  Kodaira, K.
\paper  On compact complex analytic surfaces
\jour  Ann. of Math.
\vol  71
\yr  1960
\pages  111--152
\moreref
\paper \ II
\jour ibid.
\vol 77
\yr 1963
\pages 563--626
\moreref
\paper \ III
\jour ibid.
\vol 78
\yr 1963
\pages 1--40
\endref


\ref \by K\"unnemann, K.
\paper On the Chow motive of an abelian scheme
\inbook Motives, Seattle 1991
\bookinfo Proc. Symp. Pure Math. 55, vol 1
\eds Jannsen, Kleiman, Serre
\publ American Mathematical Society
\publaddr Providence
\yr 1994
\pages 189--205
\endref

\ref \by Manin,~Yu.
\paper Correspondences, motifs and monoidal transformations
\jour Mat. Sbornik
\vol 77
\yr 1968
\pages 475--507
\moreref
\jour AMS Transl.
\yr 1970
\endref

\ref
\by  Miyake, T.
\book  Modular Forms
\publ  Springer
\publaddr  Berlin, etc.
\yr  1989
\endref

\ref  \by Murre, J.P.
\paper  On the motive of an algebraic surface
\jour  J. reine angew. Math.
\vol  409
\yr 1990
\pages 190--204
\endref

\ref  \bysame 
\paper  On a conjectural filtration of the Chow groups of an algebraic
variety
\jour  Indag. Math. (New Series)
\vol  4
\yr 1993
\pages 177--188
\moreref \ II
\pages 189--203
\endref

\ref \by Schoen, C.
\paper Complex multiplication cycles on elliptic modular threefolds
\jour Duke Math. J.
\vol 53
\yr1986
\pages 771--794
\endref

\ref  \by Scholl, A.J.
\paper Motives for modular forms
\jour Invent. Math.
\vol 100
\yr 1990
\pages 419--430
\endref

\ref  \bysame 
\paper Classical motives
\inbook Motives, Seattle 1991
\bookinfo Proc. Symp. Pure Math. 55, vol 1
\eds Jannsen, Kleiman, Serre
\publ American Mathematical Society
\publaddr Providence
\yr 1994
\pages 163--187
\endref

\ref \by Shermenev
\paper The motive of an abelian variety
\jour Funct. Anal.
\vol 8
\yr 1974
\pages 55--61
\endref

\ref \by Shioda, T.
\paper  On elliptic modular surfaces
\jour  J. Math. Soc. Japan
\vol  24
\yr  1972
\pages  20--59
\endref

\ref \by \v Sokurov, V.V.
\paper  Holomorphic differential forms of higher degree on Kuga's modular
          varieties
\jour  Math. USSR Sbornik
\vol  30
\yr  1976
\pages  1199--142
\endref

\ref \bysame 
\paper  The study of the homology of Kuga varieties
\jour  Math USSR Izvestija
\vol  16
\yr  1981
\pages  399--418
\endref

\ref \by Weil, A.
\book Foundations of Algebraic Geometry
\publ American Mathematical Society
\publaddr Providence
\yr 1948
\endref 

\ref \by  S. Zucker
\paper  Hodge theory with degenerating coefficients: L$_2$-cohomology in
the
        Poincar\'e metric
\jour  Ann. of Math.
\vol  109
\yr  1979
\pages  415--476
\endref


\endRefs
\enddocument